\documentclass[10pt,dvipsnames,reqno]{amsproc}

\pdfoutput=1 
\DeclareUnicodeCharacter{0301}{\'{a}}
\DeclareUnicodeCharacter{0301}{\'{e}}
\DeclareUnicodeCharacter{0301}{\'{i}}
\DeclareUnicodeCharacter{0304}{****}

\usepackage{filecontents}

\usepackage[style=numeric, backend=biber, backref=true, bibencoding=utf8, datamodel=mrnumber, maxbibnames=6, sorting=nyt]{biblatex}

\DeclareFieldFormat{mrnumber}{%
  \ifhyperref
    {\href{http://www.ams.org/mathscinet-getitem?mr=1#1}{MR#1}}
    {MR#1}}

\DeclareFieldFormat{mrclass}{#1}

\DeclareNameAlias{bymrreviewer}{byeditor}

\newbibmacro*{mrinfo}{%
  \printfield{mrnumber}%
  \iffieldundef{mrclass}
    {\setunit*{\addcomma\space}}
    {\setunit*{\addspace}}%
  \printfield{mrclass}%
  \setunit*{\addcomma\space}%
  \ifnameundef{mrreviewer}
    {}
    {\bibstring{byreviewer}%
     \setunit{\addspace}%
     \printnames[bymrreviewer]{mrreviewer}}}

\newtoggle{bbx:mrinfo}
\DeclareBibliographyOption[boolean]{mrinfo}[true]{\settoggle{bbx:mrinfo}{#1}}
\renewbibmacro*{doi+eprint+url}{%
  \iftoggle{bbx:doi}
    {\printfield{doi}}
    {}%
  \newunit\newblock
  \iftoggle{bbx:mrinfo}
    {\usebibmacro{mrinfo}}
    {}%
  \newunit\newblock
  \iftoggle{bbx:eprint}
    {\usebibmacro{eprint}}
    {}%
    }

\usepackage{amsbsy, amsmath, amssymb, amsthm,
  caption, colortbl, diagbox, enumitem, graphicx, mathtools,
  mdframed, soul, subcaption, xspace, url}

\usepackage{symbols}
\usepackage[utf8]{inputenc}
\usepackage{xcolor}
\usepackage{hyperref}
\hypersetup{
     colorlinks=true,
     linkcolor=blue,
     filecolor=blue,
     citecolor=blue,
     urlcolor=blue,
     }
\usepackage{threeparttable}
\usepackage{multirow}
\definecolor{col}{RGB}{240,180,60}
\usepackage{float}
\restylefloat{table}

\allowdisplaybreaks
\usepackage[text={13.8cm,25.2cm}, centering]{geometry}

\makeatletter
\def\section{\@startsection{section}{1}%
  \z@{.36\linespacing\@plus\linespacing}{.36em}%
  {\large\scshape\MakeUppercase}}

\def\subsection{\@startsection{subsection}{2}%
  \z@{.36\linespacing\@plus.36\linespacing}{-.36em}%
  {\normalfont\scshape}}

\def\subsubsection{\@startsection{subsubsection}{3}%
  \z@{0.276\baselineskip}{-.36em}%
  {\normalfont\itshape}}

\makeatother

\theoremstyle{plain}

\newtheorem{thm}{Theorem}[subsection]

\theoremstyle{definition}

\newtheorem{df}[thm]{Definition}
\newtheorem{ex}[thm]{Example}
\newtheorem{rem}[thm]{Remark}

\title[Closure and closed morphisms]{Internal Neighbourhood Structures II: Closure and closed morphisms}
\author{Partha Pratim Ghosh}
\address{Department of Mathematical Sciences \\
  University of South Africa \\
  Unisa Science Campus \\
  corner of Christian de Wet \& Pioneer Avenue \\
  Florida 1709 \\
  Johannesburg, Gauteng \\
  South Africa}
\address{National Institute for Theoretical and Computational Sciences (NITheCS), South Africa}

\email{ghoshpp@unisa.ac.za}
\thanks{Support received from {\em European Union Horizon 2020 MCSA Irses project 731143} is thankfully acknowledged.
}

\begin{document}

\setlength{\abovedisplayskip}{0pt}
\setlength{\abovedisplayshortskip}{0pt}
\setlength{\belowdisplayskip}{0pt}
\setlength{\belowdisplayshortskip}{0pt}

\begin{abstract}
  Internal preneighbourhood spaces inside any finitely complete category with finite coproducts and proper factorisation structure were first introduced in my earlier paper. This paper proposes a closure operation on internal preneighbourhood spaces and investigates closed morphisms and its close allies. Consequently it introduces  analogues of several well known classes of topological spaces for preneighbourhood spaces. Some preliminary properties of these spaces are established in this paper. The results of this paper exhibit preneighbourhood systems are more general than closure operators and conveniently allows identifying properties of classes of morphisms independent of \emph{continuity} of morphisms with respect to induced closure operators.
\end{abstract}

\makeatletter \@namedef{subjclassname@2020}{\textup{2020} Mathematics
  Subject Classification} \makeatother


\maketitle

\section{Introduction}
\label{sec:introduction}

The notion of an internal preneighbourhood space was first considered
in \cite[][]{2020}. The present paper introduces a \emph{closure
  operator} on an internal preneighbourhood space (see Definition
\ref{df:closure-closed}{}). The \emph{closure operator} entails in discussing
\emph{closed morphisms} (see Definition \ref{df:closed-morphisms}{}). The
rest of the paper discuss notions closely aligned with \emph{closed
  morphisms} --- \emph{dense morphisms} (see Definition
\ref{df:dense-mor}{}), \emph{proper morphisms} (see Definition
\ref{df:proper-mor}{}), \emph{separated morphisms} (see Definition
\ref{df:separated-morphism}{}) and \emph{perfect morphisms} (see
Definition \ref{df:perfect-mor}). Alongside morphisms special
classes of internal preneighbourhood spaces are introduced:
\emph{compact} spaces (see Definition \ref{df:compact-obj}{}), \emph{Hausdorff spaces} (see Definition \ref{df:hausdorff}{}), \emph{compact
  Hausdorff spaces} (see Definition
\ref{df:some-spaces}{\ref{item:cpt-hausdorff}}), \emph{Tychonoff spaces} (see Definition
\ref{df:some-spaces}\ref{item:tychonoff}{}) and \emph{absolutely
  closed spaces} (see Definition
\ref{df:some-spaces}\ref{item:abs-closed}{}). Detailed investigation
on the special classes of internal preneighbourhood spaces shall be
done in later papers. A quick perusal of Table \ref{tab:mor-types-comp-table}{} provide a glimpse of results achieved in this paper as well as helps to compare  similar results in literature, e.g., in \cite{ClementinoGiuliTholen2004}. The table clearly exhibit the extent to which continuity of morphisms with respect to induced closure operations are essential in achieving these properties.

The paper is organised as follows:
\begin{enumerate}[label=(\alph*),wide,labelindent=0ex]
\item In \S\ref{sec:preliminaries}{} notions necessary for the paper are briefly introduced; in the process some seemingly new observations have been listed. In this paper a monotonic, extensional and grounded endomap on a poset is called a closure operator (see \S\ref{sssec:closure-operation}{} for terminology).
  Given a complete lattice $L$, $\mathtt{EGM}(L)$ denotes the complete lattice of closure operators on $L$, $\mathtt{CBSMSL}(L)$ denotes the complete lattice of complete bounded sub-$\wedge$-semilattices of $L$. Proposition \ref{prop:cbsmsl-dually-refl-in-egm}{} shows \opp{\mathtt{CBSMSL}(L)}{} is reflectively embedded in $\mathtt{EGM}(L)$ as the idempotent closure operations.

\item \S\ref{sec:clos-from-pnbd}{} introduce \emph{closure operation} on preneighbourhood spaces.
  \begin{enumerate}[label=(\roman*),wide,labelindent=2ex]

  \item \label{item:sec-3.2}{}
    Each preneighbourhood system $\mu$ on an object $X$, for each $p \in \Sub{X}{\mathsf{M}}$, partitions the subobjects of $X$ in four subsets. The first partition consists of subobjects which are \emph{far away} from $p$ (see equation \eqref{eq:pnbd-induce-far-subobj}{}); such subobjects are incompatible with $p$. The second partition consist of subobjects which are incompatible with $p$ and not \emph{far away} from $p$; the third partition consists of  subobjects $x > p$ and the fourth of subobjects $x \leq p$. The first partition is a down-set of \Sub{X}{\mathsf{M}}\footnote{A subset $P \subseteq L$ of a lattice $L$ is a \emph{down-set} if it is non-empty and $x \leq y \in P \Rightarrow x \in P$. A down-set of the form $\downarrow p = 
  \bigl\{
  x \in L: x \leq p
  \bigr\}$ is a \emph{principal down-set} and is the smallest down-set containing $p$.} (see Lemma \ref{lem:far-prop}{} for details); in the special instance when \Sub{X}{\mathsf{M}}{} is a frame and $\mu$ is a neighbourhood system (see Definition \ref{df:ins}{}) the first partition is a principal down-set.

The set of subobjects in the fourth and second partitions should be the ones of concern in defining  \emph{closure} \cls{p}{\mu}{} of $p$, see Definition \ref{df:closure-closed}; the fixed subobjects of \cls{}{\mu}{} are \emph{$\mu$-closed subobjects}. The assignment $p \mapsto \cls{p}{\mu}$ so defined is a closure operation, its initial properties discussed in Theorem \ref{thm:clos-in-pscompl}{.} The closure \cls{}{\mu}{} is not additive in general (unless when \Sub{X}{\mathsf{M}}{} is atom generated and distributive, Theorem \ref{thm:clos-in-pscompl}\ref{item:dist-atomic-give-additive-cls}{}), nor idempotent (unless \Sub{X}{\mathsf{M}}{} is atom generated and $\mu$, in particular, a neighbourhood system, see Theorem \ref{thm:clos-in-pscompl}\ref{item:pnbd-open-gen-atom-gen-give-idemp-cls} for details). In case when \Sub{X}{\mathsf{M}}{} is pseudocomplemented then for a neighbourhood system $\mu$ on $X$, there is a Galois connection between the semilattices of closed subobjects and the open subobjects yielding a dual equivalence between regular closed and regular open subobjects (see Proposition \ref{prop:intr-clos-reciprocal}{} and Remark \ref{rem:reg-clos-reg-open-dual-equiv}{}).

\item Continuing from \S\ref{sssec:closure-operation} and specialising to the complete lattice \Sub{X}{\mathsf{M}}{,}  in \S\ref{ssec:closure-vs-pnbd}{} it is shown that the complete lattice of preneighbourhood systems on $X$ dually contains a coreflective copy of closure operations on \Sub{X}{\mathsf{M}} (see Theorem \ref{thm:closures-in-pnbd}{}).  This exhibits the generality of the approach via preneighbourhood systems in comparison with closure operations.

\item In absence of idempotence for \cls{}{\mu}{,} there exists its \emph{idempotent hull}  $\clsm{}{\mu} \geq \cls{}{\mu}$ (see Proposition \ref{prop:cbsmsl-dually-refl-in-egm}{} and Remark \ref{rem:idemp-hull-largest}) having the same set of \mc{\mu}{} of closed subobjects. The notion of continuity with respect to \cls{}{\mu}{} as well \clsm{}{\mu}{} are discussed in \S\ref{ssec:closure-hereditary}. Proposition \ref{prop:cont-wrt-cl-idem-hull}{} shows continuity with respect to \cls{}{\mu}{} implies continuity with respect to \clsm{}{\mu}{;} continuity with respect to \clsm{}{\mu} (respectively, \cls{}{\mu}{}) is called \emph{$\mu$-$\phi$ continuity} (respectively, \emph{$\mu$-$\phi$ continuity with respect to closures}) or \emph{continuity} (respectively, \emph{continuity with respect to closures}) in short. Theorem \ref{thm:closure-continuous}{} shows every admissible monomorphism is continuous with respect to closures;  equation \eqref{eq:closure-idemp-hereditary}{} provides correspondence between closed subobjects of a closed subobject and  closed subobjects of  whole space (also Remark \ref{rem:closed-embed}{}).

  Note: for a morphism \Arr{f}{X}{Y}{} of the base category, its property of being \emph{continuous with respect to preneighbourhood system $\mu$ on $X$ and $\phi$ on $Y$} is precisely the definition of it being a preneighbourhood morphism \Arr{f}{\opair{X}{\mu}}{\opair{Y}{\phi}} (see Definition \ref{df:ins}\ref{item:nbdmaps}{}). On the other hand, each preneighbourhood system $\mu$ on $X$ induce closure operations \cls{}{\mu}, \clsm{}{\mu}{,} and \emph{continuity with respect to these closure operations}  is separate and is not affected by the presence or absence of \emph{continuity with respect to preneighbourhood systems} in general.

\item A major obstruction to effecting continuity with respect to closures of a morphism is the inclusion of subobjects in the fourth partition in \ref{item:sec-3.2} while computing the join. An antidote to this occur when \Sub{X}{\mathsf{M}}{} is \emph{atom generated}: firstly, the closure has a simpler description (Remark \ref{rem:closure-when-atom-generated}){,} consequently, continuity with respect to closures for a large class of morphisms is obtained (Corollary \ref{cor:cont-examples}{\ref{item:refl-0-atomic-atom-gen-pres-atom-is-cont}{}}). However in general, continuity of every preneighbourhood morphism is ensured once every dense preneighbourhood morphism is continuous (Proposition \ref{prop:cont-suffices-weak} and Corollary  \ref{cor:cont-iff-dense-are-cont}{}).

\item \S\ref{ssec:closure-ex}{} illustrate notions in some specific contexts. Notable amidst them are closures with respect to functorial neighbourhood systems \cite[see][Theorem 3.38 and Definition 4.3]{2020} on locales, groups and commutative rings without identity. On locales it is known from \cite{2020}{} that the $T$-neighbourhood systems (see equation \eqref{eq:T-nbd-system}, \cite{DubeIghedo2016b,DubeIghedo2016c}) are functorial; it is shown here that the closure with respect to the $T$-neighbourhood system is precisely the usual closure of a sublocale \cite[see][\S III.8{}]{PicadoPultr2012}; furthermore every localic map \Arr{f}{X}{Y}{} is continuous with respect to any preneighbourhood system $\mu$ on $X$, $\phi$ on $Y$ if $\mu$ is larger than the $T$-neighbourhood system on $X$ and $\phi$ is smaller than the $T$-neighbourhood system on $Y$ (see \S{}\ref{sssec:locales}).
  Example \ref{ex:groups}{} shows the normal closure induces a functorial neighbourhood system $\nu_{X}$ on each group $X$. \S\ref{sssec:groups}{} show $\cls{A}{\nu_X} = 
  \bigl\{
  x \in X: \text{ normal closure of } x \text{ meets }A\bigr\}$; furthermore every group homomorphism \Arr{f}{X}{Y}{} is continuous with respect to any preneighbourhood system $\mu \geq \nu_{X}$ on $X$ and $\phi \leq \nu_Y$.
  Similarly, Example \ref{ex:rngs}{} shows the ideal closure of a subring induces a functorial neighbourhood system $\iota_X$ on ring $X$. \S\ref{sssec:rngs}{} show $
  \cls{A}{\iota_X} =
  \bigl\{
  x \in X: (\exists r \in X)(rx \in A)
  \bigr\}$; furthermore every ring homomorphism \Arr{f}{X}{Y}{} is continuous with respect to any preneighbourhood system $\mu \geq \iota_{X}$ on $X$ and $\phi \leq \iota_Y$.
\end{enumerate}

\item \S\ref{ssec:closed-morphisms}{} discuss \emph{closed morphism}, i.e., given the preneighbourhood spaces \opair{X}{\mu}{,} \opair{Y}{\phi}{,} morphisms \Arr{f}{X}{Y}{} which preserve closed subobjects (see equation \eqref{eq:closed-morphisms}{}). Theorem \ref{thm:closed-morphisms}{} provide properties of closed morphisms, while Theorem  \ref{thm:reflzero-fs-imply-closed}{} provide sufficient examples of closed preneighbourhood morphisms.

\item \S\ref{ssec:dense-mor}{} introduce \emph{dense morphisms} and their properties are discussed in Theorem \ref{thm:dense-mor-prop}. It is shown in Theorem \ref{thm:dense-mor-prop}\ref{item:dense-clemb-pfs}{} that if every preneighbourhood morphism is continuous then every preneighbourhood morphism factors as a dense preneighbourhood morphism followed by a closed embedding. This factorisation system is a proper factorisation system for the full subcategory of internal Hausdorff spaces (see Remark \ref{rem:dense-closed-emb-proper}{}).

\item \S\ref{ssec:stably-clos-morph}{} discuss stably closed morphisms, called \emph{proper morphisms}. Theorem \ref{thm:proper-mor-prop}{} discuss properties of proper morphisms. In \S\ref{ssec:compact-obj} \emph{compact preneighbourhood spaces} are introduced as preneighbourhood spaces \opair{X}{\mu}{} for which the unique morphism \Arr{\terma{X}}{\opair{X}{\mu}}{\opair{\termo}{\nabla_{\termo}}}{} is proper. The full subcategory \Comp{\pNHD{\Bb{A}}{}} of compact preneighbourhood spaces is shown to be finitely productive, closed hereditary if every preneighbourhood morphism is continuous (Theorem \ref{thm:compact-prop}\ref{item:cpt-fin-productive-closed-her}{}).

\item \S\ref{ssec:separated-morphisms}{} discuss \emph{separated morphisms} (Definition \ref{df:separated-morphism}{}), Theorem \ref{thm:sep-map-prop}{} discuss properties of separated morphisms. Internal \emph{Hausdorff spaces} are introduced in Definition \ref{df:hausdorff} as those preneighbourhood spaces \opair{X}{\mu}{} for which \Arr{\terma{X}}{\opair{X}{\mu}}{\opair{\termo}{\nabla_{\termo}}}{} is a separated morphism, alternate characterisations are provided in Theorem \ref{thm:ihs-alt}, the full subcategory \Int{\Haus}{\pNHD{\Bb{A}}}{} of internal Hausdorff spaces is shown to be finitely complete, closed under subobjects and images of preneighbourhood morphisms stably continuous and stably in \emph{E} (Corollary \ref{cor:hausisfincomplete}{}).

\item Finally in \S\ref{ssec:perfect-morphisms}{} \emph{perfect morphisms} are discussed --- they are preneighbourhood morphisms which are both proper and separated (Definition \ref{df:perfect-mor}{}). The properties of perfect morphisms are discussed in Theorem \ref{thm:perfect-map-prop}{.} The paper concludes introducing \emph{compact Hausdorff spaces}, \emph{Tychonoff spaces} and \emph{absolutely closed spaces}, Definition \ref{df:some-spaces}{.}

\end{enumerate}

The effort of internalising the notion of space has been pursued in different ways, at least in the references below as well as citations within them:
\begin{enumerate}[label=\roman*., wide, labelindent=0ex]
\item using closure operators as in \cite{DikranjanGiuli1985,Castellini1986,CastelliniPhdThesis,DikranjanGiuli1987,DikranjanGiuliTozzi1988,DikranjanGiuliTholen1989,DikranjanGiuli1989,Castellini1988,CastelliniStrecker1990,Castellini1990,Dikranjan1992,CastelliniKoslowskiStrecker1992c,Castellini1992,CastelliniKoslowskiStrecker1993,CastelliniKoslowskiStrecker1994b,CastelliniHajek1994,ClementinoGiuliTholen1996,Castellini2000,Castellini2001b,CastelliniGiuli2001,ClementinoGiuliTholen2001,Castellini2003,CastelliniHolgate2003b,DikranjanTholenWatson2004,MCDHWT2004,CastelliniGiuli2005,Castellini2008,CastelliniHolgate2010,DikranjanTholen2015},

\item using interior operators and neighbourhood operators \cite[see][page 2, for definition]{2020} as in \cite{Vorster2000,CastelliniRamos2010,Castellini2011,Razafindrakato2012,Razafindrakoto2012b,CastelliniMurcia2013,HolgateRazafindrakato2014,Castellini2015,Castellini2016,HolgateRazafindrakato2017}, 

\item using convergence structures as in \cite{EklundGahler1992,Lowen1979,DikranjanGiuli1993,Slapal2001,Slapal2002,Slapal2005,GiuliSlapal2005,Slapal2008,TomSlapal2010,Slapal2011a,Slapal2011c,Slapal2012,Slapal2016a,Slapal2016b}

\item using a set of axioms for closed morphisms as in  \cite{ClementinoGiuliTholen2004}

\item using a set of axioms for proper morphisms as in \cite{HofmannTholen2012}. 
\end{enumerate}

In all of them the aspect of \emph{continuity} of morphisms is built inside the axioms, or is an easy consequence of the axioms --- for instance \cite[see ][\S 11.1, (F6) and its consequences]{ClementinoGiuliTholen2004}. The approach in present work is transversal: firstly a category with \emph{nice}
properties is shown to have a structure of categorical neighbourhood
system. An object endowed with a neighbourhood system allows the formulation of a closure operation. In several convenient cases the closure operator possess good familiar properties. The continuity of morphisms with respect to induced closure operation in general is not immediate and has to be checked. However, in presence of continuity with respect to closure operations nicer properties are ensured as summarised by Table \ref{tab:mor-types-comp-table}{;} however, several other properties do not need the presence of continuity with respect to the induced closure operation. Thus apart from the generalisation that the method allows it also reveals the extent to which the condition of continuity (with respect to closure operations) is required in obtaining familiar properties of well known classes of morphisms. It is essential to emphasise the generalisation obtained herein is conservative.
   
  The notation and terminology adopted in this paper are largely in
  line with the usage in \cite{Maclane1998} or \cite[][]{Borceux1994-01}. Finally, for any set $I$, the symbol $2^{I}$ (respectively, \subfin{I}) denotes the set of all subsets (respectively, finite subsets) of $I$ and a set $A$ is \emph{small} if it is a member of some set.


\section{Preliminaries}
\label{sec:preliminaries}
This section recalls facts relevant for this paper. In the process some observations are seemingly new.

\subsection{}
\label{ssec:preliminaries:posets}
This section establish notations and terms with regards to posets as used in this paper.

\begin{subsubsection}{}
  \label{sssec:closure-operation}{}
  Given a poset $P$ with a smallest element \textsf{0} and a largest element \textsf{1}, an order preserving endomap \Arr{f}{P}{P}{} is called \emph{extensional} if $x \leq f(x)$ ($x \in P$), \emph{grounded} if $f(\mathsf{0}) = \mathsf{0}$. A grounded and extensional order preserving endomap on $P$ is called a \emph{closure operation} on $P$ and $\mathsf{EGM}(P)$ is the set of all closure operations on $P$. Evidently $\mathsf{EGM}(P)$ is ordered pointwise, i.e., $c \leq d$ if $c(x) \leq d(x)$, for each $x \in P$, where $c, d \in \mathsf{EGM}(P)$. The poset $\mathsf{EGM}(P)$ has smallest closure operation \id{P} and largest closure operation \Arr{\lambda}{P}{P}{}, where $\lambda(x) =
  \begin{cases}
    \mathsf{0}, & \text{ if }x = \mathsf{0} \\
    \mathsf{1}, & \text{ otherwise}
  \end{cases}
  $. A closure operation $c \in \mathsf{EGM}(P)$ is \emph{idempotent} if $\comp{c}{c} = c$. For each $c \in \mathsf{EGM}(P)$, $\Fix{c} = 
  \bigl\{
  x \in P: c(x) = x
  \bigr\}$ is the set of \emph{fixed points} of $c$.

  Given a complete lattice $L$, for each $c \in \mathsf{EGM}(L)$, $\mathsf{0}, \mathsf{1} \in \Fix{c}$ and \Fix{c}{} is closed under arbitrary meets, i.e., \Fix{c}{} is
a complete bounded sub-$\wedge$-semilattice of $L$. Define:
\begin{equation}
  \label{eq:idemp-majorise}
  \hat{c}(y) = \bigwedge
  \bigl\{
  x \in \Fix{c}: y \leq x
  \bigr\}, \qquad\text{ for }y \in L.
\end{equation}
Let $\mathtt{CBSMSL}(L)$ be the set of all complete bounded
sub-$\wedge$-semilattices of the complete lattice $L$. Since an
intersection of complete bounded sub-$\wedge$-semilattices is a complete
bounded sub-$\wedge$-semilattice, the set $\mathtt{CBSMSL}(L)$ is a
complete lattice with $\{\mathsf{0}, \mathsf{1}\}$ as the smallest
complete bounded sub-$\wedge$-semilattice of $L$, $L$ the largest complete
bounded sub-$\wedge$-semilattice of $L$ and intersection of 
complete bounded sub-$\wedge$-semilattices as meet. Finally, for each
$P \in \mathtt{CBSMSL}(L)$ define:
\begin{equation}
  \label{eq:cbsmsl-to-egm-op}
  \nu_P(y) = \bigwedge
  \bigl\{
  x \in P: y \leq x
  \bigr\}, \qquad\text{ for }\in L.
\end{equation}

\begin{Prop}
  \label{prop:cbsmsl-dually-refl-in-egm}
  For every complete lattice $L$,
  $\xymatrix{ {\mathtt{EGM}(L)} \ar@<1ex>[r]^-{\mathtt{Fix}}
    \ar@<-1ex>@{<-}[r]_-{\nu} \ar@{}[r]|-{\bot} &
    {\opp{\mathtt{CBSMSL(L)}}} }$ in the category of partially ordered
  sets and order preserving maps with
  $\comp{\mathtt{Fix}}{\nu} = \id{\opp{\mathtt{CBSMSL(L)}}}$.

  Furthermore, for any $c \in \mathtt{EGM}(L)$,
  $\hat{c} = \nu_{\Fix{c}}$ and for any ordinal $\alpha$ if:
  \begin{equation}
    \label{eq:c-ord-chain}
    c^{\alpha} =
    \begin{cases}
      \id{L}, & \text{ if } \alpha = 0 \\
      \comp{c}{c^{\beta}}, & \text{ if }\alpha = \beta + 1 \text{ is a non-limit ordinal} \\
      \bigvee_{\beta<\alpha}c^{\beta}, & \text{ if } \alpha > 0 \text{ is a limit ordinal}
    \end{cases},
  \end{equation}
  then $c^{\alpha} \in \mathsf{EGM}(L)$, $c \leq c^{\alpha} \leq \hat{c}$,
  $\comp{\hat{c}}{c^{\alpha}} = \hat{c} = \comp{c^{\alpha}}{\hat{c}}$ and
  $\widehat{c^{\alpha}} = \hat{c}$, for all $\alpha \geq 1$.
\end{Prop}

\begin{proof}
  If $c, d \in \mathtt{EGM}(L)$, $P, Q \in \mathtt{CBSMSL}(L)$ then:
  \begin{enumerate}[label=\roman*.]

  \item for any $y \in L$, since $P$ is a complete
    sub-$\wedge$-semilattice of $L$, $\nu_P(y) \in P$, proving
    $\nu_P \in \mathtt{EGM}(L)$ is idempotent (i.e.,
    $\comp{\nu_P}{\nu_P}=\nu_P$) and $\Fix{\nu_P} = P$;
    
  \item $c \leq d$ and $x \in \Fix{d}$ imply $c(x) \leq d(x) = x$, proving
    \Arr{\mathtt{Fix}}{\mathtt{EGM}(L)}{\opp{\mathtt{CBSMSL}(L)}{}} is
    order preserving;
  \item if $P \subseteq Q$ then for each $x \in L$, then
    $\nu_Q(x) = \bigwedge \bigl\{ t \in Q: x \leq t \bigr\} \leq \bigwedge \bigl\{ t \in P: x \leq t
    \bigr\} = \nu_P(x)$, proving
    \Arr{\nu}{\opp{\mathtt{CBSMSL}(L)}}{\mathtt{EGM}(L)}{} is order
    preserving;
    
  \item
    $c \leq \nu_P \Leftrightarrow \biggl( t \in P \Rightarrow \bigl( x \leq t \Rightarrow c(x) \leq t \bigr) \biggr)
    \Leftrightarrow P \subseteq \Fix{c}$, proving \adjt{\mathtt{Fix}}{\nu};
  \end{enumerate}
  proving the first part of the statement.  For the second part, since
  $c^{0} = \id{L} \leq c^{1} = c \leq c^{2} = \comp{c}{c}$, transfinite
  induction implies the first two conditions; hence for any
  $c \leq d \leq \hat{c}$, $\Fix{c} = \Fix{d} = \Fix{\hat{c}}$,
  $\hat{d} = \hat{c}$ and
  $\comp{\hat{c}}{d} = \hat{c} = \comp{d}{\hat{c}}$, completing the
  proof.
\end{proof}

\begin{rem}\label{rem:idemp-hull-largest}{}
  The idempotent closure operation $\hat{c} = \nu_{\Fix{c}}$ is the
  smallest idempotent closure operation larger than $c$; it is called
  the \emph{idempotent hull of $c$} \cite[see][\S4.6, for more
  properties of $\hat{c}$]{DikranjanTholen1995}. 
\end{rem}

\begin{rem}\label{rem:idemp-hull-as-small-join-of-succ-clos-operations}{}
  In the context of \eqref{eq:c-ord-chain}{,} if there exists an
  ordinal $\alpha$ such that $c^{\alpha} = c^{\alpha+1}$, then
  $c^{\alpha}$ is idempotent and hence $c^{\alpha} = \hat{c}$
  \cite[see][\S{4.6}]{DikranjanTholen1995}. Thus, if the underlying set of the lattice $L$ is a small set then for each $c \in \mathsf{EGM}(L)$, $\hat{c} = c^{\alpha}$ for some ordinal $\alpha$.
\end{rem}

\begin{rem}
  The adjunction provides a formula for the join in
  $\mathtt{CBSMSL}(L)$:
  \begin{equation}
    \label{eq:join-in-cbsmsl}
    \bigvee\CAL{P} = \Fix{\bigwedge_{P\in\CAL{P}}\nu_P}, \qquad\text{ for all }\CAL{P} \subseteq
    \mathtt{CBSMSL}(L).
  \end{equation}
\end{rem}

\begin{rem}
  The complete lattice $\mathtt{CBSMSL}(L)$ of all complete bounded
  sub-$\wedge$-semilattices of $L$ is dually reflectively embedded inside
  the complete lattice $\mathtt{EGM}(L)$ of grounded closure
  operations on $L$ as the idempotent closure operations.
\end{rem}

\end{subsubsection}

\begin{subsubsection}{}
  In a lattice $L$ for any $a \in L$ the following two formulae are equivalent:
  \begin{align*}
    a \neq \mathsf{0} & \text{ and } (\forall x \in L)
    \bigl(
    x \leq a \Rightarrow x = \mathsf{0}\text{ or }x = a
    \bigr), \\
    a \neq \mathsf{0} & \text{ and } (\forall x \in L)
    \bigl(
    a \leq x \text{ or } a \wedge x = \mathsf{0}
    \bigr),
  \end{align*}
  and the element $a$ defined by such is an \emph{atom} of $L$; the set of atoms of
$L$ is denoted by \atom{L}{.} A lattice $L$ is \emph{atomic} if for each $x \in L$, there exists
  $a \in \atom{L}$ with $a \leq x$, and \emph{atom
    generated} if
  $x = \bigvee \bigl\{ a \in \atom{L}: a \leq x \bigr\}$\footnote{Usages differ, e.g., in
    \cite[][Chapter IV]{Birkhoff1979} the terms \emph{atomic} and
    \emph{atom generated} are not distinguished, as in this paper, and
    the term \emph{atomic} is used to mean atom generated.}.
  Evidently, every atom generated lattice is atomic, but the converse
  need not be true --- e.g., consider the lattice of divisors of a positive
  natural number. In case of the lattice \Sub{X}{\mathsf{M}} (see \S\ref{sssec:pfs}{}){,} the
  symbol \atom{X}{} abbreviate
  \atom{\Sub{X}{\mathsf{M}}}{.}
  
\end{subsubsection}

\begin{subsubsection}{}
Recall from \cite{Gratzer2003}: in a complete lattice $L$, a
\emph{pseudocomplement} of $a \in L$ is the element $a^{*} \in L$ such
that $x \leq a^{*} \Leftrightarrow x \wedge a = \mathsf{0}$, i.e.,
$a^{*} = \max\bigl\{x \in L: a \wedge x = \mathsf{0}\bigr\}$; a complete
lattice is \emph{pseudocomplemented} if every element has a
pseudocomplement. Evidently, $\mathsf{0}^{*} = \mathsf{1}$,
$\mathsf{1}^* = \mathsf{0}$, $p \leq q \Rightarrow q^{*} \leq p^{*}$, the assignment
$x \mapsto x^{**}$ is an idempotent closure operation on $L$,
$(p \wedge q)^{**} = p^{**} \wedge q^{**}$ and $(\bigvee S)^{*} = \bigwedge_{s \in S}s^{*}$ ($S \subseteq L$). Clearly, $p \parallel p^{*} \Leftrightarrow p, p^{*} \neq \mathsf{0}$, where $x \parallel y$ means $x$ and $y$ are incompatible.

An element $p \in L$ is said to be \emph{implicative} if the order preserving endomap \Arr{p\wedge{-}}{L}{L}{} has a right adjoint \Arr{\implyr{p}{-}}{L}{L}. Evidently, $p \in L$ is implicative if and only if $p \wedge -$ preserve arbitrary joins. Note: $p^{*} = \implyr{p}{\mathsf{0}}$. 
\end{subsubsection}

\subsubsection{}
\label{sssec:downsets}
For a lattice $L$, a \emph{down-set} is a subset $D \subseteq L$ such that $x \leq y \in D \Rightarrow x \in D$, and a \emph{principal down-set} is a set of the form $\downarrow p = 
\bigl\{
x \in L: x \leq p
\bigr\}$ for some $p \in L$. Evidently, a principal down-set $\downarrow p$ is the smallest down-set containing $p$, and for every down-set $D$, $D = \bigcup_{p \in D}(\downarrow p)$. The set of down-sets of $L$ is $\mathtt{Dn}(L)$ is a complete lattice with intersections being the meet and unions as join, the smallest down-set being $\downarrow \mathsf{0} = \{\mathsf{0}\}$ and $\downarrow \mathsf{1} = L$ being the largest down-set.

An \emph{up-set} of $L$ is a down-set of \opp{L}{,} and a \emph{filter} in $L$ is an up-set closed under finite meets.
\subsection{}
\label{ssec:preliminaries:pnbd}

Internal preneighbourhood spaces were considered in
\cite{2020}.

\subsubsection{}
\label{sssec:pfs}

Let \Bb{A}{} be a finitely complete category with finite coproducts.

A morphism $f$ of \Bb{A}{} is said to be \emph{orthogonal} to a morphism
$g$, written \down{f}{g}{} if there exists a unique morphism $w$ such
that $v = \comp{g}{w}$ and $u = \comp{w}{f}$ whenever
$\comp{v}{f} = \comp{g}{u}$. For $\CAL{H} \subseteq \Bb{A}_1$ let
$\up{\CAL{H}} = \bigl\{ f \in \Bb{A}_1: h \in \CAL{H} \Rightarrow \down{f}{h}
\bigr\}$,
$\dn{\CAL{H}} = \bigl\{ f \in \Bb{A}_1: h \in \CAL{H} \Rightarrow \down{h}{f}
\bigr\}$; then $\xymatrix{ {2^{\Bb{A}_1}} \ar@<1ex>[r]^-{\up{-}} \ar@<-1ex>@{<-}[r]_-{\dn{-}} \ar@{}[r]|-{\bot} & {
    \left(
      2^{\Bb{A}_1}
    \right)^{\mathsf{op}}}  }$, a pair \opair{\CAL{A}}{\CAL{B}}{} of subsets of
\low{\Bb{A}}{1}{} is called a \emph{prefactorisation system} if
$\CAL{B} = \dn{\CAL{A}}$ and $\CAL{A} = \up{\CAL{B}}$; a prefactorisation system \opair{\CAL{A}}{\CAL{B}}{} is a \emph{factorisation system} if every morphism factors as a \CAL{A}{-}morphism followed by a \CAL{B}{-}morphism and a factorisation system \opair{\CAL{A}}{\CAL{B}}{} is \emph{proper} if $\CAL{A} \subseteq \Epi{\Bb{A}}$, $\CAL{B} \subseteq \Mono{\Bb{A}}$. The (possibly large) set of prefactorisation systems on \Bb{A}{} is a complete poset with $\opair{\CAL{A}}{\CAL{B}} \leq \opair{\CAL{A'}}{\CAL{B'}}$ if $\CAL{A} \supseteq \CAL{A'} \Leftrightarrow \CAL{B} \subseteq \CAL{B'}$ \cite[see][\S{2}, for details]{Janel1997b}. Given a factorisation system \opair{\CAL{A}}{\CAL{B}}{,} $\CAL{A} \subseteq \Epi{\Bb{A}}$ if and only if for each object $X$ the \emph{diagonal} $d_{X} \in \CAL{B}$, where $\xymatrixcolsep{7.8em}\xymatrix{ {X} \ar@{>->}[]!<2.4ex,0ex>;[r]^-{d_X=\opair{\id{X}}{\id{X}}} & {X \times X}  }$,
\cite[see][Proposition 14.11]{AdamekHerrlichStecker1990}. Hence a factorisation system \opair{\CAL{A}}{\CAL{B}}{} is proper if and only if for every object $X$, $d_{X} \in \CAL{B}$ and the \emph{codiagonal} $c_{X} \in \CAL{A}$, where \Arr{c_X}{X+X}{X}{} is the unique morphism such that $\comp{c_X}{\iota_1} = \id{X} = \comp{c_X}{\iota_2}$, $\iota_1,\iota_2$ are the coproduct injections, and equivalently $\ExtEpi{\Bb{A}} \subseteq \CAL{A} \subseteq \Epi{\Bb{A}}$ and $\ExtMon{\Bb{A}} \subseteq \CAL{B} \subseteq \Mono{\Bb{A}}$.

Given a proper \fact{\textsf{E}}{\textsf{M}} system, a \textsf{M}-subobject of an object $X$, also called an \emph{admissible subobject} of $X$, is a $m \in \mathsf{M}$ with codomain $X$, any two equivalent admissible subobjects of $X$ considered equal. The set of admissible subobjects of $X$ is denoted by \Sub{X}{\mathsf{M}}{.} In this paper the
  morphisms of \textsf{E} are depicted with arrows like \fepi{}{}{}{}
  while the morphisms of \textsf{M} are depicted with arrows like
  $\xymatrix{ {} \ar@{>->}[r] & {}}$.  If \Arr{f}{X}{Y}{} be a
  morphism, then $\xymatrixcolsep{4.8em}\xymatrix{ {X} \ar@{->>}[r]^-{f^{\mathsf{E}}} & {\Img{f}} \ar@{>->}[]!<3.6ex,0ex>;[r]^-{f^{\mathsf{M}}} & {Y}  }$ is the \opair{\mathsf{E}}{\mathsf{M}} factorisation of $f$; more generally, if $m \in \Sub{X}{\mathsf{M}}$ (respectively, $n \in \Sub{Y}{\mathsf{M}}$) then
  the \emph{image} of $m$ (respectively, \emph{preimage} of $n$) under $f$ is
  \img{f}{m} (respectively, \finv{f}{n}{}){,} where
  $\comp{f}{m} = \comp{(\img{f}{m})}{\rest{f}{m}}$ (respectively,
  $\comp{f}{(\finv{f}{n})} = \comp{n}{f_n}$) is the
  \fact{\textsf{E}}{\textsf{M}} of \comp{f}{m}{} (respectively,
  pullback of $n$ along $f$), \rest{f}{m}{} is the \emph{restriction
    of $f$ on $m$} (respectively, $f_{n}$ is the \emph{corestriction}
  of $f$ on $n$); obviously $f^{\mathsf{E}} = \rest{f}{\id{X}}$ and $f^{\mathsf{M}} = \img{f}{\id{X}}$. In presence of finite limits and finite coproducts, \Sub{X}{\mathsf{M}}{} is a lattice, the largest subobject is \id{X}{} and the smallest object  is $\sigma_{X}$, where  $\xymatrix{ {\inito} \ar@{->>}[r]^-{\inita{\emptyset_X}} & {\emptyset_X} \ar@{>->}[]!<12pt,0pt>;[r]^-{\sigma_{X}} & {X} }$ is the \opair{\mathsf{E}}{\mathsf{M}}{-}factorisation of the unique morphism \inita{X}{} from the initial object to $X$. The image and preimage induce adjunction $\xymatrix{ {\Sub{X}{\mathsf{M}}} \ar@<1ex>[r]^-{\img{f}{}} \ar@{<-}@<-1ex>[r]_-{\finv{f}{}} \ar@{}[r]|-{\bot} & {\Sub{Y}{\mathsf{M}}}  }$ for each morphism \Arr{f}{X}{Y}{} of \Bb{A}{.}
  A filter $F$ on $X$ is a filter in  \Sub{X}{\mathsf{M}}{}; the (possibly large) set of filters on $X$ is \Fil{X}, which is a complete algebraic lattice, distributive if and only if \Sub{X}{\mathsf{M}}{} is distributive (\cite[see][Theorem 1.2]{IberkleidMcGovern2009} or \cite[][Proposition 2.7, Corollary 2.8]{2020}{}), with compact elements $\uparrow p = 
  \bigl\{
  x \in \Sub{X}{\mathsf{M}}: x \geq p
  \bigr\}$ ($p \in \Sub{X}{\mathsf{M}}$). For each morphism \Arr{f}{X}{Y}{} the adjunction \adjt{\img{f}{}}{\finv{f}{}}{} induce adjunctions $\xymatrix{ {\Fil{X}} \ar@<1ex>[r]^-{\imgfil{f}{}} \ar@{<-}@<-1ex>[r]_-{\invfil{f}{}} \ar@{}[r]|-{\top} & {\Fil{Y}} }$, where:
  \begin{align}
    \label{eq:imgfil}
    \imgfil{f}{A}
    & = 
      \bigl\{
      y \in \Sub{Y}{\mathsf{M}}: \finv{f}{y} \in A
      \bigr\},
    & \text{ for }A \in \Fil{X}, \\ \intertext{and}
    \label{eq:invfil}
    \invfil{f}{B}
    & = 
      \bigl\{
      x \in \Sub{X}{\mathsf{M}}: (\exists b \in B)(\finv{f}{b} \leq x)
      \bigr\},
    & \text{ for }B \in \Fil{Y}.
  \end{align}

  \begin{subsubsection}{}
    \label{sssec:smallest-element}{}
    A connection between the smallest subobjects of objects need to be highlighted.

    \begin{Prop}
      \label{prop:smallest-subobjects-connected}{}
      For every \Arr{g}{Z}{Y}{,} there exists a unique $\xymatrix{ {\inito_Z} \ar@{->>}[r]^{\omega_{Z,Y}} & {\inito_Y}  }$ such that $\comp{g}{\sigma_Z} = \comp{\sigma_Y}{\omega_{Z,Y}}$ is the \opair{\textsf{E}}{\textsf{M}}-factorisation of \comp{g}{\sigma_Z}{.}

      In particular: $g \in \mathsf{M} \Rightarrow \omega_{Z,Y} \in \Iso{\Bb{A}}$, if $\Homo{\Bb{A}}{Y}{Z} \neq \emptyset$ as well then $\omega_{Y,Z} = \inv{\omega_{Z,Y}}$, and hence $\omega_{Y,Y} = \id{\inito_Y}$.
    \end{Prop}

    \begin{proof}
      Given a morphism \Arr{g}{Z}{Y}{,} consider the diagram:
    \begin{equation}\label{eq:smallest-subobjects-connected}
      \xymatrix{
        & {\inito} \ar@{->>}[dl]_-{\inita{\inito_Z}} \ar@/_{18ex}/[ddl]_-{\inita{Z}}
        \ar@{->>}[dr]^-{\inita{\inito_{Y}}} \ar@/^{18ex}/[ddr]^-{\inita{Y}} \\
        {\inito_Z} \ar@{>->}[]!<0ex,-2.4ex>;[d]_-{\sigma_Z} \ar@{.>}[rr]^-{!\,w_{g}}
        & & {\inito_Y} \ar@{>->}[]!<0ex,-2.4ex>;[d]^-{\sigma_Y} \\
        {Z} \ar[rr]_-g & & {Y}        
      }.
    \end{equation}
    Since $\comp{g}{\comp{\sigma_Z}{\inita{\inito_Z}}} = \comp{g}{\inita{Z}} = \inita{Y} = \comp{\sigma_Y}{\inita{\inito_Y}}$, there exists the unique morphism $w_{g}$ such that $\comp{w_{g}}{\inita{\inito_Z}} = \inita{\inito_Y}$ and $\comp{g}{\inita{Z}} = \comp{\inita{Y}}{w_{g}}$; evidently, $w_{g} \in \mathsf{E}$. If \Arr{g'}{Z}{Y}{} is any other morphism, then $\comp{w_{g'}}{\inita{\inito{Z}}} = \inita{\inito_Y} = \comp{w_g}{\inita{\inito_Z}}$ implies $w_{g } = w_{g'}$, i.e., $w_{g}$ is independent of the choice of $g$; furthermore, it is evident that $w_{g} = \rest{g}{\sigma_Z}$ and the square represents the \opair{\textsf{E}}{\textsf{M}}-factorisation of \comp{g}{\sigma_Z}{.} Taking $\omega_{Y,Z} = w_{g}$ completes the proof.
    \end{proof}

    Thus: $\inito_{\termo}$ is an \textsf{E}-image of each $\inito_Y$ and if $\inita{\termo}$ is an admissible monomorphism then each $\inito_Y \simeq \inito$ (see Theorem \ref{thm:reflzero-prop}\ref{item:morph-refl-zero=str-initial}, Remark \ref{rem:adm-quasi-pt}{}).
    
  \end{subsubsection}
  
\subsubsection{}
\label{sssec:pnbds}

Preneighbourhood systems can now be defined.

\begin{Df}
  \label{df:ins}
  \begin{enumerate}[label=(\alph*),wide,labelindent=0ex]

  \item \label{item:context} A \emph{context} is $\CAL{A} = (\Bb{A}, \mathsf{E}, \mathsf{M})$, where \Bb{A} is a finitely complete category with finite coproducts and a proper factorisation system \opair{\mathsf{E}}{\mathsf{M}} such that for each object $X$, \Sub{X}{\mathsf{M}}{} is a complete lattice.
    
  \item \label{item:nbdsystems}{}
    An order preserving map
    \Arr{\mu}{\opp{\Sub{X}{\mathsf{M}}}}{\Fil{X}} is a
    \emph{preneighbourhood system} if $\mu(\sigma_{X}) = \Sub{X}{\mathsf{M}}$ and 
    $p \in \mu(m) \Rightarrow m \leq p$; if further,
    $p \in \mu(m) \Rightarrow (\exists q \in \mu(m))(p \in \mu(q))$ then
    $\mu$ is a \emph{weak neighbourhood system}; moreover if
    $\mu(\bigvee S) = \bigcap_{s\in S}\mu(s)$ ($S \subseteq \Sub{X}{\mathsf{M}}$) then
    $\mu$ is a \emph{neighbourhood system}. A pair \opair{X}{\mu}{}, where
    $X$ is an object of \Bb{A} and $\mu$ is a preneighbourhood system on
    $X$ is called an \emph{internal preneighbourhood} space.  Likewise
    for \emph{internal weak neighbourhood space} and \emph{internal
      neighbourhood space}. 

  \item \label{item:nbdmaps}{}
    If \opair{X}{\mu}{} and \opair{Y}{\phi}{} are
    internal preneighbourhood spaces then a morphism \Arr{f}{X}{Y}{}
    is a \emph{preneighbourhood morphism} if
    $p \in \phi(u) \Rightarrow \finv{f}{p} \in \mu(\finv{f}{u})$; if \opair{X}{\mu} and
    \opair{Y}{\phi}{} are internal neighbourhood spaces and \finv{f}{}{}
    preserve joins then it is a \emph{neighbourhood morphism}.  The category of internal
    preneighbourhood spaces and preneighbourhood morphisms is \pNHD{\Bb{A}}{};
    \wNHD{\Bb{A}}{} is the full subcategory of internal weak
    neighbourhood spaces and \NHD{\Bb{A}}{} is the subcategory of
    internal neighbourhood spaces and neighbourhood morphisms.

  \item \label{item:open}{}
    Given a preneighbourhood system $\mu$ on $X$, a subobject $p \in \Sub{X}{\mathsf{M}}$ is \emph{$\mu$-open} if $p \in \mu(p)$; the (possibly large) set of all $\mu$-open sets is \mo{\mu}{.}

  \item \label{item:topology}{} A neighbourhood system $\mu$ on $X$ is a \emph{topology} on $X$ if \mo{\mu}{} is a frame in the partial order of \Sub{X}{\mathsf{M}}{.} If $\mu$ is a topology on $X$ then \opair{X}{\mu}{}  is an \emph{internal topological space}, \Int{\Top}{\Bb{A}}{} is the full subcategory of \NHD{\Bb{A}}{} of all internal topological spaces.
    
  \item \label{item:fs}{}
    A morphism \Arr{f}{X}{Y}{} is \emph{formally surjective} (or, also
  referred to in literature as \emph{semistable}, e.g., in \cite[][]{Vermeulen2001}) if for each
  $y \in \Sub{Y}{\mathsf{M}}$ there exists a
  $x \in \Sub{X}{\mathsf{M}}$ such that $y = \img{f}{x}$, or
  equivalently for every $y \in \Sub{Y}{\mathsf{M}}$ the corestriction
  $f_{y}$ is in \textsf{E}.

\item \label{item:FR}{}
  A morphism \Arr{f}{X}{Y}{} is a \emph{Frobenius morphism} if for each $x \in \Sub{X}{\mathsf{M}}$ and $y \in \Sub{Y}{\mathsf{M}}$, $\img{f}{(x\wedge\finv{f}{y})} = y \wedge \img{f}{x}$.

\item \label{item:refl0}{}
  A morphism \Arr{f}{X}{Y}{} is said to \emph{reflect zero} if $\finv{f}{\sigma_Y} = \sigma_{X}$.
  \end{enumerate}
\end{Df}

\begin{rem}
  The condition \ref{item:nbdmaps}{} is often called a \emph{continuity condition} with respect to preneighbourhood systems. 
\end{rem}

\begin{rem}
  The set of preneighbourhood systems on $X$ is \pnhd{X}; likewise \wnhd{X}{,}  \nhd{X}{,} \tnhd{X}{} denote the set of weak neighbourhood systems, neighbourhood systems and topologies on $X$ respectively. Each of \pnhd{X}, \wnhd{X}, \nhd{X} are complete lattices \cite[see][Theorem 3.17 \& Theorem 3.32]{2020}, while \tnhd{X}{} is a complete sublattice on \nhd{X}{} if and only if there exists a largest topology on $X$ \cite[see][Theorem 3.36]{2020}{.} 
\end{rem}

\begin{rem}
  Given the preneighbourhood systems $\mu$ on $X$, $\phi$ on $Y$, a morphism \Arr{f}{X}{Y}{} of \Bb{A} is a preneighbourhood morphism if and only if for any $x \in \Sub{X}{\mathsf{M}}$, $y \in \Sub{Y}{\mathsf{M}}$ any one of the following three conditions is true: $\invfil{f}{\phi(y)} \subseteq \mu(\finv{f}{y})$, $\phi(y) \subseteq \imgfil{f}{\mu(\finv{f}{y})}$, $\invfil{f}{\phi(\img{f}{x})} \subseteq \mu(x)$, \cite[see][Theorem 3.40]{2020}\footnote{The assignment $x \mapsto \invfil{f}{\phi(\img{f}{x})}$, $x \in \Sub{X}{\mathsf{M}}$ (respectively, $y \mapsto \imgfil{f}{\mu(\finv{f}{y})}$, $y \in \Sub{Y}{\mathsf{M}}$) is a preneighbourhood system \cite[see][]{2020}, and is denoted by \invfil{f}{\phi\img{f}{}} (respectively, \imgfil{f}{\mu\finv{f}{}}{}).}. The symbol \Arr{f}{\opair{X}{\mu}}{\opair{Y}{\phi}}{} is used to denote $f$ is a preneighbourhood morphism.
\end{rem}

Contexts abound --- if \Bb{A} is finitely complete, finitely cocomplete
and has {\em all} intersections then there is a
\fact{\Epi{\Bb{A}}}{\ExtMon{\Bb{A}}} system on \Bb{A}{;} in
particular, every small complete, small cocomplete category \Bb{A}, if
well powered have a context
$\CAL{E} = (\Bb{A}, \Epi{\Bb{A}}, \ExtMon{\Bb{A}})$, and if co-well
powered have a context
$\CAL{M} = (\Bb{A}, \ExtEpi{\Bb{A}}, \Mono{\Bb{A}})$. As special cases
are the contexts:
$(\FinSet, \mathsf{Surjections}, \mathsf{Injections})$
\cite[see][Example 3.7]{2020},
$(\Set, \mathsf{Surjections}, \mathsf{Injections})$ \cite[see][Example
3.8]{2020}, $(\Grp, \mathsf{RegEpi}, \mathsf{Mono})$
\cite[see][Example 3.9 \& Proposition 3.10]{2020},
$(\Alg{\opair{\Omega}{\Xi}}, \mathsf{RegEpi}, \mathsf{Mono})$
\cite[see][Example 3.11 \& Proposition 3.12]{2020},
$(\Top, \mathsf{Epi}, \mathsf{ExtMono})$ \cite[see][Example
3.13]{2020}, $(\Loc, \mathsf{Epi}, \mathsf{RegMono})$
\cite[see][Example 3.14]{2020}, every topos with its usual
factorisation structure \cite[see][page 5, (iii)]{2020}, every
lextensive category \cite[see][]{CarboniLackWalters1993} with a proper
factorisation structure \cite[see][page 6, (v)]{2020} (and this
includes \Cat, \opp{\CRing}{}, \Sch{}, \opp{\Bb{A}}{} where \Bb{A}{} is a Zariski category \cite[see][Definition 1.2]{Diers1992}). Also given any context
$\CAL{A} = (\Bb{A}, \mathsf{E}, \mathsf{M})$ and any object $X$ of
\Bb{A},
$\slice{\CAL{A}}{X} = (\slice{\Bb{A}}{X}, \slice{\mathsf{E}}{X},
\slice{\mathsf{M}}{X})$ is the context where
\begin{equation}
  \label{eq:fact-sys-on-bundles}
  \begin{aligned}
    \slice{\mathsf{E}}{X} & = \bigl\{\Arr{e}{(X, x)}{(Y, y)}: e \in
    \mathsf{E}\bigr\} \\
    \slice{\mathsf{M}}{X} & = \bigl\{\Arr{M}{(X, x)}{(Y, y)}: m \in
    \mathsf{M}\bigr\}, \\
  \end{aligned}
\end{equation}
\cite[see][page 5, (iv)]{2020} and \cite[][\S 2.10, for
details]{ClementinoGiuliTholen2004}.

\subsubsection{}
\label{sssec:open}{}
Given a preneighbourhood system $\mu$ on $X$, $x \in \Sub{X}{\mathsf{M}}$:
\begin{equation}
  \label{eq:intr}
  \intr{x}{\mu} = \bigvee
  \bigl\{
  p \in \mo{\mu}: p \leq x
  \bigr\}
\end{equation}
is the \emph{$\mu$-interior} of $x$. Evidently, $\sigma_{X}, \id{X} \in \mo{\mu}$, \mo{\mu}{} is closed under arbitrary joins if and only if each $\intr{x}{\mu} \in \mo{\mu}$ \cite[see][Theorem 3.20]{2020}. Further, if each $\mu$-interior is $\mu$-open then for each $x \in \Sub{X}{\mathsf{M}}$ the following two conditions are equivalent:
\begin{align}
  \label{eq:mu-union-of-up-open}
  \mu(x) & = 
         \bigcup \bigl\{
         \uparrow p: x \leq p \in\mo{\mu}
         \bigr\} \\
  \label{eq:mu-open-generated}
  p \in \mu(x) & \Leftrightarrow x \leq \intr{p}{\mu}.
\end{align}

A preneighbourhood system $\mu$ is said to have \emph{open interiors} if $\intr{x}{\mu}\in\mo{\mu}$ for each $x \in \Sub{X}{\mathsf{M}}$ and is \emph{open generated} if it has open interiors and satisfies \eqref{eq:mu-open-generated}{.} If $\mu$ have open interiors then:
\begin{equation}
  \label{eq:intr-pres-meet}
  \intr{(x \wedge y)}{\mu} = \intr{x}{\mu} \wedge \intr{y}{\mu},
\end{equation}
and \Arr{\intr{}{\mu}}{\Sub{X}{\mathsf{M}}}{\Sub{X}{\mathsf{M}}}{} is a meet preserving \emph{intensional} (i.e., $\intr{x}{\mu} \leq x$ for each $x \in \Sub{X}{\mathsf{M}}$) idempotent endomap with \mo{\mu}{} as its fixed set \cite[see][Corollary 3.26]{2020}{.} Finally, every neighbourhood system $\mu$ is open generated \cite[see][Theorem 3.27]{2020}{.}

\begin{subsubsection}{}
  \label{sssec:refl0}{}
  The condition of a morphism reflecting zero shall be used in the paper. This section state some necessary facts about them.
  
\begin{Thm}
  \label{thm:reflzero-prop}
  Given the morphisms \Arr{f}{X}{\Arr{g}{Y}{Z}} of \Bb{A}, the following statements hold.
  
  \begin{enumerate}[label=(\alph*)]

  \item \label{item:reflzero-alt} \tfae{:}
    \begin{enumerate}[label=(\roman*)]
    \item \label{item:refl0}{} $f$ reflects zero, i.e., $\finv{f}{\sigma_Y} = \sigma_{X}$.

    \item \label{item:refl0-alt}{} For any $x \in \Sub{X}{\mathsf{M}}$:
      \begin{equation}
        \label{eq:reflzero-alt}
        \img{f}{x} = \sigma_Y \Rightarrow x = \sigma_X.
      \end{equation}

    \item \label{item:cont-pren-morph-cond} For all
      $x \in \Sub{X}{\mathsf{M}}$, $y \in \Sub{Y}{\mathsf{M}}$:
      \begin{equation*}
        \label{eq:cont-pren-morph-cond}
        y \wedge \img{f}{x} = \sigma_{Y} \Rightarrow x \wedge \finv{f}{y} = \sigma_{X}.
      \end{equation*}
    \end{enumerate}
 
  \item \label{item:adm-mono-reflzero}{} If $\comp{\finv{f}{}}{\img{f}{}} = \id{\Sub{X}{\mathsf{M}}}$ then $f$ reflects zero. In particular, every admissible morphism reflect zero.
    
  \item \label{item:reflzero-comp-closed} The set of morphisms
    reflecting zero is closed under compositions.
    
  \item \label{item:refl0-left-cancel} If \comp{g}{f}{} reflects zero
    then $f$ reflects zero.
    
  \item \label{item:refl0-her} For any morphism \Arr{f}{X}{Y}
    reflecting zero and $n \in \Sub{Y}{\mathsf{M}}$, the corestriction
    $f_{n}$ on $N$ reflects zero.

  \item \label{item:morph-refl-zero=str-initial} Let \Bb{A}{} be a
    category with pullbacks and initial object \inito{}. If every
    morphism of \Bb{A}{} reflect zero then \inito{} is strict. Further
    if the unique morphism
    \Arr{\inita{\termo}=\terma{\inito}}{\inito}{\termo}{} is an
    admissible monomorphism and \inito{} is strict then every morphism
    reflects zero.
  \end{enumerate}
\end{Thm}

\begin{proof}
  Towards the proof of the equivalence in \ref{item:reflzero-alt}: the
  equivalence of \ref{item:refl0} and \ref{item:refl0-alt}{} follows
  from the adjunction \adjt{\img{f}{}}{\finv{f}{}}; since
  $\img{f}{(x\wedge\finv{f}{y})} \leq y \wedge\img{f}{x}$, \ref{item:refl0-alt}
  implies \ref{item:cont-pren-morph-cond}, while taking
  $y = \sigma_{Y}$ and $x = \id{X}$, 
  \ref{eq:cont-pren-morph-cond}{} implies \ref{item:refl0}{.} The statements in \ref{item:adm-mono-reflzero}-\ref{item:refl0-left-cancel}{} follow from \ref{item:reflzero-alt}{.} Towards a proof of
  \ref{item:refl0-her}{,} in the pullback of
  $m \in \Sub{Y}{\mathsf{M}}$ along $f$, if $f$ reflect zero then
  \begin{equation*}
    \comp{(\finv{f}{m})}{(\finv{f_m}{\sigma_M})} =
    \finv{f}{(\comp{m}{\sigma_M})}
    = \finv{f}{\sigma_Y} = \sigma_X
  \end{equation*}
  implies $\finv{f_m}{\sigma_M} = \sigma_{\finv{f}{M}}$, since
  $\finv{f}{m} \in \mathsf{M}$, proving \ref{item:refl0-her}{.}
  Finally, the proof of
  \ref{item:morph-refl-zero=str-initial} follows from Proposition \ref{prop:smallest-subobjects-connected}{.}
\end{proof}

\begin{rem}
  \label{rem:adm-quasi-pt}{} A finitely complete category with
    an initial object is \emph{quasi-pointed}
    (\cite[see][\S{1}]{Bourn2001}, \cite[][]{GoswamiJanelidze2017}) if
    the unique morphism \Arr{\inita{\termo}}{\inito}{\termo}{} is a
    monomorphism. In many contexts, e.g., in each of
    $(\FinSet, \mathsf{Surjections}, \mathsf{Injections})$,
    $(\Set, \mathsf{Surjections}, \mathsf{Injections})$,
    $(\Top, \mathsf{Epi}, \mathsf{ExtMono})$ or
    $(\Loc, \mathsf{Epi}, \mathsf{RegMon})$ the unique morphism
    \inita{\termo}{} is a regular monomorphism, and hence an
    admissible monomorphism. A context \CAL{A}{} is called
    \emph{admissibly quasi-pointed} if its underlying category
    \Bb{A}{} has the unique morphism \inita{\termo}{} an admissible
    monomorphism. Thus in an
    admissibly quasi-pointed context, the initial object is strict if
    and only if every morphism reflects zero.
  \end{rem}
\begin{rem}
  Using Proposition \ref{prop:smallest-subobjects-connected}{} given the coterminating morphisms $f$ and $g$ consider the diagram in \eqref{eq:refl0-pb}  where the front vertical square is the pullback of $f$ along $g$, the top horizontal square is the pullback of $\sigma_{Z}$ along $f_{g}$; if both $f$ and $g$ reflect zero then the vertical right hand and base horizontal squares are pullback squares, enabling the existence of the unique morphism $w$ to make the whole diagram commute; further all the squares are pullback squares; in particular, $\finv{f_g}{\sigma_Z} = \finv{g_f}{\sigma_{X}}$. Hence $f_{g}$ reflects zero if and only if $g_{f}$ reflects zero.
  \begin{equation*}\tag{$\blacktriangle$}\label{eq:refl0-pb}
    \xymatrixcolsep{3.6em}
    \xymatrix{
      & {\finv{f_g}{\inito_Z}} \ar[rr]^-{(f_{g})_{\sigma_{Z}}}
      \ar@{.>}[dd]|(0.492){\hole}|(0.66){!\,w}
      \ar@{>->}[]!<-1.8ex,-1.92ex>;[dl]_-{\finv{f_{g}}{\sigma_{Z}}} & &
      {\inito_{Z}} \ar@{->>}[dd]^-{\omega_{Z,Y}}
      \ar@{>->}[]!<-1.8ex,-1.92ex>;[dl]|-{\sigma_Z} \\
      {X \times_Y Z} \ar[rr]|(0.66){f_{g}} \ar[dd]_-{g_{f}} & & {Z} \ar[dd]|(0.66){g} \\
      & {\inito_X} \ar@{->>}[rr]|(0.36){\omega_{X,Y}}|(0.54){\hole}
      \ar@{>->}[]!<-1.8ex,-1.92ex>;[dl]|-{\sigma_X} & & {\inito_Y}
      \ar@{>->}[]!<-1.8ex,-1.92ex>;[dl]|-{\sigma_Y} \\
      {X} \ar[rr]_-f & & {Y}
    },
  \end{equation*}
\end{rem}

  \begin{Df}
    \label{df:r0c}{}
      A context $\CAL{A} = (\Bb{A}, \mathsf{E}, \mathsf{M})$ called a \emph{reflecting zero context} if all morphisms reflect
    zero. 
  \end{Df}

\end{subsubsection}

\begin{subsubsection}{}
  \label{sssec:fr+fs}{}
  This section exhibits a connection between formally surjective morphisms and Frobenius morphisms.
\begin{Prop}
  \label{prop:fr-fs-connection}{}
  Every Frobenius \textsf{E}-morphism is formally surjective; if
  \Arr{f}{X}{Y}{} is formally surjective and \img{f}{}{} is a
  homomorphism of meet semilattices then $f$ is Frobenius.
\end{Prop}

\begin{proof}
  If $f$ is a Frobenius morphism then for each
  $y \in \Sub{Y}{\mathsf{M}}$:
  \begin{equation*}
    \img{f}{\finv{f}{y}} = \img{f}{(\id{X}\wedge\finv{f}{y})} = y\wedge \img{f}{\id{X}};
  \end{equation*}
  since $\img{f}{\id{X}} = \id{Y} \Leftrightarrow f \in \mathsf{E}$, every Frobenius
  \textsf{E}-morphism is formally surjective. The second part is
  trivial.
\end{proof}

\begin{rem}
  If $\mathtt{FS}[\Bb{A}]$ (respectively, $\mathtt{FR}[\Bb{A}]$)
  denote the (possibly large) set of formally surjective
  (respectively, Frobenius) morphisms of \Bb{A}{} then the following
  connections:
  \begin{equation}
    \label{eq:fr-fs-connection}
    \xymatrixcolsep{4.8em}
    \xymatrix{
      {\framebox[1.2\width]{\txt{\textsf{E} is pullback stable}}} \ar@2{->}[r] &
      {\framebox[1.2\width]{$\mathsf{E} \subseteq \mathtt{FS}[\Bb{A}]$}} \ar@2{->}[r] &
      {\framebox[1.2\width]{$\Bb{A}_1 \subseteq \mathtt{FR}[\Bb{A}]$}} \ar@2{->}[d] \\
      & &
      {\framebox[1.2\width]{$\mathsf{E} \subseteq \mathtt{FR}[\Bb{A}]$}} \ar@2{->}[ul]
    },
  \end{equation}
  are well known (see \cite[][Theorem 5.13]{2020} or
  \cite[][Proposition 1.3]{ClementinoGiuliTholen1996}).
\end{rem}

\end{subsubsection}

\begin{subsubsection}{}
  \label{sssec:topologicity}{}
  The forgetful functor \Arr{U}{\pNHD{\Bb{A}}}{\Bb{A}}{} is a
topological functor \cite[see][Theorem 4.8(a)]{2020}. Consequently, 
each limit (respectively, colimit) object, unless mentioned to the contrary, is considered as an internal preneighbourhood space with the smallest
(respectively, largest) preneighbourhood system which make each of the
components of the limiting (respectively, colimiting) cone
preneighbourhood morphisms. Thus, for instance:

\begin{enumerate}[label=\roman*., wide, labelindent=0ex]
\item \label{item:terminal-pnbd}{}
  The terminal object \termo{} being the empty product is always equipped with the smallest preneighbourhood system $\nabla_{\termo}$ (see \eqref{eq:pnbd-smallest-largest}{}). Note: \Sub{\termo}{\mathsf{M}}{} is not always trivial. For instance, in $(\opp{\CRing}, \mathsf{Epi}, \mathsf{RegMono})$ the terminal object is the commutative ring \Bb{Z}{} of integers, $\Sub{\Bb{Z}}{\mathsf{RegMono}} = 
\bigl\{
n\Bb{Z}: n \geq 0
\bigr\}$, hence $\nabla_{\Bb{Z}} < \uparrow_{\Bb{Z}}$.

\item \label{item:rest-pnbd}{}
  Given an admissible
monomorphism $\xymatrix{ {M} \ar@{>->}[]!<2.4ex,0ex>;[r]^-m & {X} }$
and a preneighbourhood system $\mu$ on $X$, $M$ is equipped with \rest{\mu}{m}{,} where for any $a \in \Sub{M}{\mathsf{M}}$:
\begin{equation}
  \label{eq:rest-pnbd}
  \begin{aligned}
    \rest{\mu}{m}(a) & = \bigl\{ u \in \Sub{M}{\mathsf{M}}: (\exists v \in
    \mu(\comp{m}{a}))(\finv{m}{v} \leq u)
    \bigr\} \\
    & = \bigl\{ u \in \Sub{X}{\mathsf{M}}: (\exists v \in \mu(\comp{m}{a}))(v \wedge m
    \leq \comp{m}{u}) \bigr\}.
  \end{aligned}
\end{equation}

\item \label{item:product-pnbd}{}
  The binary product  $\xymatrix{ {X} \ar@{<-}[r]^-{p_1} & {X \times Y} \ar[r]^-{p_2} & {Y}  }$ in \Bb{A}{} of the preneighbourhood spaces \opair{X}{\mu}{,} \opair{Y}{\phi} is equipped with $\mu\times\phi$, where for any $\opair{x}{y} \in \Sub{X \times Y}{\mathsf{M}}$:
\begin{equation}
  \label{eq:product-pnbd}
  (\mu\times\phi)\opair{x}{y} = \invfil{p_1}{\mu(x^{\mathsf{M}})} \vee
  \invfil{p_2}{\phi(y^{\mathsf{M}})}.
\end{equation}

\item \label{item:pullback-pnbd}{}
  If $\xymatrix{ {\opair{X}{\mu}} \ar[r]^-f & {\opair{Z}{\psi}} \ar@{<-}[r]^-g & {\opair{Y}{\phi}}  }$ and $\xymatrix{ {X \times_Z Y} \ar[r]^-{f_g} \ar[d]_-{g_f} & {Y} \ar[d]^-g \\ {X} \ar[r]_-f & {Z}  }$ is the pullback of $f$ along $g$ in \Bb{A} then $X \times_Z Y$ is equipped with $\mu \times_{\psi} \phi$, where for any $\opair{x}{y} \in \Sub{X \times_Z Y}{\mathsf{M}}$:
\begin{equation}
  \label{eq:pullback-pnbd}
  (\mu\times_{\psi}\phi)\opair{x}{y} = \invfil{g_f}{\mu(x^{\mathsf{M}})} \vee
  \invfil{f_g}{\phi(y^{\mathsf{M}})}.
\end{equation}

\end{enumerate}

\end{subsubsection}

\subsubsection{}
\label{sssec:examples}
 
In any context $\CAL{A} = (\Bb{A}, \mathsf{E}, \mathsf{M})$ the
(possibly large) set \pnhd{X} of all preneighbourhood systems on $X$
is a complete lattice \cite[see][Theorem 3.17]{2020}. The smallest is
the \emph{indiscrete} neighbourhood system
\Arr{\nabla_X}{\opp{\Sub{X}{\mathsf{M}}}}{\Fil{X}}{} and the largest is the
\emph{discrete} neighbourhood system
\Arr{\uparrow_X}{\opp{\Sub{X}{\mathsf{M}}}}{\Fil{X}}{,} where
\begin{equation}
  \label{eq:pnbd-smallest-largest}
  \begin{aligned}
    \nabla_X(x) =
    \begin{cases}
      \Sub{X}{\mathsf{M}}, & \text{ if }x = \sigma_{X} \\
      \{\id{X}\}, & \text{ if }x \neq \sigma_{X}
    \end{cases}
  \end{aligned}\quad\text{ and }\quad
  \uparrow_X(x) = 
  \bigl\{
  p \in \Sub{X}{\mathsf{M}}: x \leq p
  \bigr\}
\end{equation}
for any $x \in \Sub{X}{\mathsf{M}}$.

\begin{Ex}
  \label{ex:finset}
  In the context
  $(\FinSet, \mathsf{Surjections}, \mathsf{Injections})$ the internal
  preneighbourhood systems are precisely extensional order preserving
  endomaps on the lattice $2^X${} of all subsets of $X$, the
  internal weak neighbourhood systems are the order preserving
  extensional idempotent endomaps on $2^X${} and the internal
  neighbourhood systems are the Kuratowski closure operations on
  $2^X${} \cite[see][Example 3.7]{2020}. Every neighbourhood
  system on a finite set precisely yield topologies,
  \cite[see][Corollary 2.13, Figure 1, for details]{2020}.
\end{Ex}

\begin{Ex}
  \label{ex:set}
  In the context $(\Set, \mathsf{Surjections}, \mathsf{Injections})$
  the internal neighbourhood systems on $X$ are precisely the
  topologies on $X$ \cite[see][Corollary 2.13, Example 3.8, Figure 1
  for details]{2020}{.}
\end{Ex}

\begin{Ex}
  \label{ex:top}{}
  In the context $(\Top, \mathsf{Epi}, \mathsf{ExtMono})$, a
  preneighbourhood system is specified by a preneighbourhood system on
  the underlying set of the topological space; preneighbourhood
  morphisms are continuous functions which are preneighbourhood
  morphisms with respect to the involved preneighbourhood systems.
  In particular, neighbourhood systems on a topological space $X$ is a
  second topology on the underlying set of the space $X$ producing
  bitopological spaces \cite[see][Example 3.13]{2020} and
  preneighbourhood morphisms are continuous functions which are also
  continuous with respect to the second topologies are neighbourhood
  morphisms.
\end{Ex}

\begin{Ex}
  \label{ex:locales}
  In the context $(\Loc, \mathsf{Epi}, \mathsf{RegMono})$
  \cite[see][Example 3.14, for details]{2020}{} a special
  neighbourhood system shall be considered in this paper, namely the
  \emph{$T$-neighbourhood system}. More precisely, given a locale $X$,
  the \emph{$T$-neighbourhood system} on it is given by
  \Arr{\tau_X}{\opp{\Sub{X}{\mathsf{RegMono}}}}{\Fil{X}}{:}
  \begin{equation}
    \label{eq:T-nbd-system}
    {\tau}_X(S) = \bigl\{ T \in \Sub{X}{\mathsf{RegMono}}: (\exists a \in
    X)(S \subseteq \mathfrak{o}(a) \subseteq T ) \bigr\},
  \end{equation}
  where $\mathfrak{o}(a) = \bigl\{ \implyr{a}{x}: x \in L \bigr\}$ is the open
  sublocale for $a \in X$ \cite[see][\S6.1.1]{PicadoPultr2012}, is an
  example of a functorial neighbourhood system on $X$
  \cite[see][Theorem 3.38 \& Definition 4.3]{2020}. $T$-neighbourhood
  systems have been used extensively in
  \cite[][]{DubeIghedo2016b,DubeIghedo2016c}{.} Since for a localic
  map \Arr{f}{X}{Y}, the preimage \finv{f}{}{} does not preserve
  arbitrary joins, with $X$ and $Y$ empowered with $T$-neighbourhood
  systems, $f$ is merely a preneighbourhood morphism and not a
  neighbourhood morphism.  Furthermore, since
  \Sub{X}{\mathsf{RegMono}}{} is a co-frame, and not a frame,
  neighbourhood systems on a locale is not an internal topology,
  internal topologies on a locale is not a reflective subcategory of
  neighbourhood spaces, \cite[see][Theorem 4.8 for details]{2020}{.}
\end{Ex}

\begin{Ex}
  \label{ex:groups}
  In the context $(\Grp, \mathsf{RegEpi}, \mathsf{Mono})$, for a group
  $X$ and a subgroup $A \subseteq X$, let $\mathsf{ncl}_{X}(A)$ denote the
  normal subgroup of $X$ generated by $A$. The order preserving map
  \Arr{\nu_X}{\opp{\Sub{X}{\mathsf{Mono}}}}{\Fil{X}} defined by:
  \begin{align}
    \label{eq:normal-closure-as-pnbd}
    \nu_X(A) & = 
             \bigl\{
             U \in \Sub{X}{\mathsf{Mono}}: \mathsf{ncl}_{X}(A) \subseteq U
             \bigr\} \nonumber \\
           & =
             \bigl\{
             U \in \Sub{X}{\mathsf{Mono}}: (\exists \normalto{N}{X})(A \subseteq N \subseteq U)
             \bigr\}
  \end{align}
  is a preneighbourhood system on $X$.  Since normal subgroups are
  closed under joins and intersections, $\nu_X$ is actually an internal
  neighbourhood system on the group $X$; moreover, every group
  homomorphism \Arr{f}{X}{Y} is a preneighbourhood morphism from the
  internal neighbourhood space \opair{X}{\nu_X}{} to \opair{Y}{\nu_Y},
  making $\nu_{X}$ a functorial preneighbourhood system,
  \cite[see][Definition 4.3]{2020}.
\end{Ex}

\begin{Ex}
  \label{ex:rngs}
  In the context $(\CRng, \mathsf{RegEpi}, \mathsf{Mono})$, since
  \CRng{} have objects commutative rings without identity, every ideal
  is a subring and feature as admissible monomorphisms.  Given any
  ring $X$, a subring $A \subseteq X$, let $\mathsf{idl}_{X}(A)$ denote the
  ideal of $X$ generated by $A$. The order preserving map
  \Arr{\iota_X}{\opp{\Sub{X}{\mathsf{Mono}}}}{\Fil{X}} defined by:
  \begin{align}
    \label{eq:normal-closure-as-pnbd}
    \iota_X(A) & = 
             \bigl\{
             U \in \Sub{X}{\mathsf{Mono}}: \mathsf{idl}_{X}(A) \subseteq U
             \bigr\} \nonumber \\
           & =
             \bigl\{
             U \in \Sub{X}{\mathsf{Mono}}: (\exists I \in \Idl{X})(A \subseteq I \subseteq U)
             \bigr\}
  \end{align}
  is a preneighbourhood system on $X$.  Since ideals are closed under
  joins and intersections, $\nu_X$ is actually an internal neighbourhood
  system on the group $X$; moreover, every ring homomorphism
  \Arr{f}{X}{Y} is a preneighbourhood morphism from the internal
  neighbourhood space \opair{X}{\iota_X}{} to \opair{Y}{\iota_Y}, making
  $\iota_X$ a functorial preneighbourhood system, \cite[see][Definition
  4.3]{2020}.
\end{Ex}


\section{Closure operations}
\label{sec:clos-from-pnbd}

This section introduce a \emph{closure operation} on each preneighbourhood space and investigate its properties.

\subsection{}
\label{ssec:closure-def-ex}

A preneighbourhood system $\mu$ on an object $X$ induces the set:
\begin{equation}
  \label{eq:pnbd-induce-far-subobj}
  \mathtt{Far}_{\mu}p = 
  \bigl\{
  x \in \Sub{X}{\mathsf{M}}: (\exists u \in \mu(x))(u \wedge p = \sigma_{X})
  \bigr\}, \qquad p \in \Sub{X}{\mathsf{M}}
\end{equation}
of admissible subobjects of $X$ which are \emph{far away} from $p$,
with respect to the preneighbourhood system $\mu$.

\begin{Lemma}
  \label{lem:far-prop}{}
  For every object $X$ of \Bb{A}{,}
  \Arr{\mathtt{Far}}{\opp{\Sub{X}{\mathsf{M}}}\times\pnhd{X}}{\mathtt{Dn}(X){}}\footnote{$\mathtt{Dn}(X)$
    is abbreviation for $\mathtt{Dn}(\Sub{X}{\mathsf{M}})$} is an
  order preserving function between complete lattices such that for
  any set $I$, $\mu, \mu_{i} \in \pnhd{X}$ ($i \in I$),
  $p, q \in \Sub{X}{\mathsf{M}}$:
  \begin{align}
    \label{eq:far-pres-bounds-for-each-pnbd}
    \mathtt{Far}_{\mu}\sigma_{X} = \Sub{X}{\mathsf{M}}
    & \text{ and } \mathtt{Far}_{\mu}\id{X} = \{\sigma_X\}, \\[1.2ex]
    \label{eq:far-pres-bounds-for-each-subobj}
    \mathtt{Far}_{\nabla_X}p = \{\sigma_{X}\} \text{ when $p \neq \sigma_{X}$}
    & \text{ and }\mathtt{Far}_{\uparrow_X}p = 
      \bigl\{
      x \in \Sub{X}{\mathsf{M}}: x \wedge p = \sigma_{X}
      \bigr\}, \\[1.2ex]
    \label{eq:far-sets-are-in-incompat}
    x,p \neq \sigma_{X}, x \in \mathtt{Far}_{\mu}p
    & \Rightarrow x \parallel p, \\[1.2ex]
    \label{eq:far-for-meet-of-pnbds}
    \mathtt{Far}_{\bigwedge_{i\in{I}}\mu_i}p
    & = \bigcap_{i\in{I}}\mathtt{Far}_{\mu_i}p,
      \text{ if $p$ is implicative}, \\[1.2ex]
    \label{eq:far-for-joins-of-pnbds}
    \mathtt{Far}_{\bigvee_{i \in I}\mu_i}p
    & = \bigcup_{J \in 2^{I}_{<\aleph_0}}\mathtt{Far}_{\mu_J}p,
      \text{ where }\mu_{J} = \bigvee_{j \in J}\mu_{j}, J \in \subfin{I}, \\[1.2ex]
    \label{eq:far-for-fin-join-of-subobj}
    \mathtt{Far}_{\mu}(p \vee q)
    & = \mathtt{Far}_{\mu}p \cap \mathtt{Far}_{\mu}q,
      \text{ if \Sub{X}{\mathsf{M}}{} is distributive}, \\[1.2ex]
    \label{eq:far-for-nbds-is-complete-ideal}
    (\forall i \in I)(x_{i} \in \mathtt{Far}_{\mu}p)
    & \Rightarrow \bigvee_{i\in{I}}x_{i} \in \mathtt{Far}_{\mu}p,
      \text{ if \Sub{X}{\mathsf{M}}{} is distributive,} \\
    \nonumber
    & \text{ $p$ is implicative and $\mu$ is a neighbourhood system on $X$}.
  \end{align}

\end{Lemma}

\begin{proof}
  The first part of the statement as well as those in
  \eqref{eq:far-pres-bounds-for-each-pnbd},\eqref{eq:far-pres-bounds-for-each-subobj}
  and \eqref{eq:far-sets-are-in-incompat}{} are simple verification.
  If \seq{\mu}{i}{I}{} is a family of preneighbourhood systems on $X$,
  $\mathtt{Far}_{\bigwedge_{i\in{I}}\mu_i}p \subseteq \bigcap_{i\in{I}}\mathtt{Far}_{\mu_i}p$
  ($p \in \Sub{X}{\mathsf{M}}$) since for any fixed $p$,
  $\mathtt{Far}_{-}p$ is monotonic. On the other hand if
  $x \in \bigcap_{i \in I}\mathtt{Far}_{\mu_i}p$ then for each
  $i \in I$ there exists a $u_{i} \in \mu_i(x)$ with
  $u_{i} \wedge p = \sigma_{X}$. If $p$ is implicative then for
  $u = \bigvee_{i \in I}u_{i} \in \bigcap_{i \in I}\mu_i(x)$,
  $u \wedge p = \bigvee_{i \in I}(u_{i} \wedge p) = \sigma_{X}$, proving
  \eqref{eq:far-for-meet-of-pnbds}{.}  Since:
  \begin{align*}
    x \in \mathtt{Far}_{\bigvee_{i\in{I}}\mu_i}p
    & \Leftrightarrow (\exists u \in \bigvee_{i \in I}\mu_i(x))(u \wedge p = \sigma_{X}) \\
    & \Leftrightarrow (\exists n \in \Bb{N})(\exists i_{0}, i_{1}, \dots, i_{n-1} \in I) \\
    & (\exists u_{0} \in \mu_{i_0}(x), u_{1} \in \mu_{i_1}(x), \dots, u_{n-1} \in \mu_{i_{n-1}}(x))(u_{1} \wedge
      u_{2} \wedge \dots \wedge u_{n} \wedge p = \sigma_{X}) \\
    & \Leftrightarrow (\exists J \in \subfin{I})(\exists u \in \mu_{J}(x))(u \wedge p = \sigma_{X}) \\
    & \Leftrightarrow (\exists J \in \subfin{I})(x \in \mathtt{Far}_{\mu_J}p) \\
    & \Leftrightarrow x \in \bigcup_{J \in\subfin{I}}\mathtt{Far}_{\mu_{J}}p,
  \end{align*}
  \eqref{eq:far-for-joins-of-pnbds}{} is proved.  From for any fixed
  $\mu \in \pnhd{X}$, $\mathtt{Far}_{\mu}$ is order reversing,
  $\mathtt{Far}_{\mu}(p\vee{q}) \subseteq \mathtt{Far}_{\mu}p \cap
  \mathtt{Far}_{\mu}q$. On the other hand, if
  $x \in \mathtt{Far}_{\mu}p\cap\mathtt{Far}_{\mu}q$ then there exist
  $u, v \in \mu(x)$ such that $u \wedge p = \sigma_{X} = v \wedge q$; if
  \Sub{X}{\mathsf{M}}{} is distributive then
  $(u \wedge v) \wedge (p \vee q) = (u \wedge v \wedge p) \vee (u \wedge v \wedge q) = \sigma_{X}$, proving
  \eqref{eq:far-for-fin-join-of-subobj}{.}  Finally, if $\mu$ is a
  neighbourhood system on $X$, for a family \seq{x}{i}{I}{} of
  elements from $\mathtt{Far}_{\mu}p$, for each $i \in I$ there exists a
  $u_{i} \in \mu(x_i)$ such that $u_{i} \wedge p = \sigma_{X}$. Since
  $u = \bigvee_{i\in{I}}u_i \in \bigcap_{i \in I}\mu(x_i) = \mu(\bigvee_{i\in{I}}x_i)$,
  $u \wedge p = \bigvee_{i \in I}(u_{i} \wedge p) = \sigma_{X}$, whenever
  $p$ is implicative, \eqref{eq:far-for-nbds-is-complete-ideal} stands
  proved.
\end{proof}

Hence, given any $\mu \in \pnhd{X}$ and $p \in \Sub{X}{\mathsf{M}}$,
\Sub{X}{\mathsf{M}}{} is partitioned into four subsets:
$\mathtt{Far}_{\mu}p$ which is a down-set (and a principal down-set in
the special case when \Sub{X}{\mathsf{M}}{} is a frame and $\mu$ a
neighbourhood system), the second is
$ \bigl\{ x \in \Sub{X}{\mathsf{M}}: x \parallel p \bigr\} \cap
(\mathtt{Far}_{\mu}p)^{c}$, the third is
$ \bigl\{ x \in \Sub{X}{\mathsf{M}}: x > p \bigr\}$ and the fourth is
the principal down-set $\downarrow p$.

\begin{Df}
  \label{df:closure-closed}
  Given any internal preneighbourhood space \opair{X}{\mu}{} define:
  \begin{equation}
    \label{eq:closure}
    \cls{p}{\mu} = \bigvee
    \bigl\{
    x  \in \Sub{X}{\mathsf{M}}: p \not< x \not\in \mathtt{Far}_{\mu}p
    \bigr\}
  \end{equation}
  the \emph{$\mu$-closure of $p$},
  $\mc{\mu} = \Fix{\cls{}{\mu}} = \bigl\{ p \in \Sub{X}{\mathsf{M}}: p =
  \cls{p}{\mu} \bigr\}$ is the (possibly large) set of
  \emph{$\mu$-closed admissible subobjects} of $X$.
\end{Df}

\begin{rem}
  Evidently
  $\cls{p}{\mu} = p \vee \bigvee \bigl\{ x \in \Sub{X}{\mathsf{M}}: x
  \parallel p \text{ and }x \notin \mathtt{Far}_{\mu}p \bigr\}$;
  moreover:
  \begin{equation*}
    \sigma_{X} \neq x < \cls{p}{\mu} \Leftrightarrow p \not< x \not\in \mathtt{Far}_{\mu}p;
  \end{equation*}
  noting the \emph{strict inequality} on the left hand side
  above. Thus for any $p \in \Sub{X}{\mathsf{M}}$ the statements:
  \begin{enumerate}[label=(\roman*)]
  \item the set
    $\mathsf{N}_{\mu,p} = \bigl\{ x \in \Sub{X}{\mathsf{M}}: p \not< x
    \not\in \mathtt{Far}_{\mu}p \bigr\}$ has a largest element;

  \item $x \parallel p \Rightarrow x \in \mathtt{Far}_{\mu}p$ for any
    $x \in \Sub{X}{\mathsf{M}}$;

  \item the sets:
    \begin{equation*}
      \bigl\{
      x \in \Sub{X}{\mathsf{M}}: x \leq p
      \bigr\},\quad
      \bigl\{
      x \in \Sub{X}{\mathsf{M}}: x > p
      \bigr\}\quad\text{ and }\quad \mathtt{Far}_{\mu}{p}    
    \end{equation*}
    make a partition of \Sub{X}{\mathsf{M}}{;}

  \item $\cls{p}{\mu} = p$;
  \end{enumerate}
  are equivalent.

  In particular:
  \begin{equation}
    \label{eq:near-elements-have-no-largest-iff-not-closed}
    p < \cls{p}{\mu} \Leftrightarrow 
    \Bigl(
    p \not< x \notin \mathtt{Far}_{\mu}p \Rightarrow (\exists y)(p \not< y \notin \mathtt{Far}_{\mu}p \text{ and }x < y)
    \Bigr).
  \end{equation}
\end{rem}

Evidently \Arr{\cls{}{\mu}}{\Sub{X}{\mathsf{M}}}{\Sub{X}{\mathsf{M}}}{}
is a closure operation on \Sub{X}{\mathsf{M}}, the idempotent hull of
\cls{}{\mu}{} is
\Arr{\clsm{}{\mu}}{\Sub{X}{\mathsf{M}}}{\Sub{X}{\mathsf{M}}}{,} where:
\begin{equation*}
  \clsm{p}{\mu} = \bigwedge
  \bigl\{
  t \in \mc{\mu}: p \leq t
  \bigr\},
\end{equation*}
both \cls{}{\mu}{,} \clsm{}{\mu}{} have the same fixed set \mc{\mu} and for
any $p \in \Sub{X}{\mathsf{M}}$, \clsm{p}{\mu}{} is the smallest
$\mu$-closed admissible subobject of $X$ larger than $p$. In view of
Remark \ref{rem:idemp-hull-as-small-join-of-succ-clos-operations}{,}
if \Sub{X}{\mathsf{M}} is a small set then
$\clsm{}{\mu} = \bigvee_{\beta\leq\alpha}\cls{}{\mu}^{\beta}$ for some ordinal $\alpha$.

\begin{rem}\label{rem:obstruction-familiar-prop}{}
  It shall turn out the condition $p \not< x$ in Definition
  \ref{df:closure-closed}{} is an obstruction to many familiar
  properties of closure (see Remark
  \ref{rem:explains-top-or-locale-or-many} and Remark
  \ref{rem:obstruction-again}{}). One antidote, as shall be exhibited,
  is \Sub{X}{\mathsf{M}} being \emph{atom generated}.
\end{rem}

\begin{rem}\label{rem:closure-when-atom-generated}{}
  Obviously,
  $\cls{p}{\mu} \geq \bigvee \bigl\{ a\in\atom{X}:a\not\in\mathtt{Far}_{\mu}p \bigr\}$;
  if \Sub{X}{\mathsf{M}}{} is atomic, then for each
  $\sigma_{X}\neq x < \cls{p}{\mu}$ there exists an atom $a \leq x$ and hence
  $a < \cls{p}{\mu}$ implying $a \not\in\mathtt{Far}_{\mu}p$. If further
  \Sub{X}{\mathsf{M}}{} is atom generated,
  $x = \bigvee \bigl\{ a \in \atom{X}:a\leq x \bigr\}$ and hence
  $\cls{p}{\mu} \geq \bigvee \bigl\{ a\in\atom{X}:a\not\in\mathtt{Far}_{\mu}p \bigr\} \geq
  x$, proving
  $\cls{p}{\mu} = \bigvee \bigl\{ a \in \atom{X}: a \not\in\mathtt{Far}_{\mu}p
  \bigr\}$.

  Thus, in \emph{atom generated} \Sub{X}{\mathsf{M}}{,} a convenient
  formula for computing the closure of a subobject is available.
\end{rem}

\begin{Thm}
  \label{thm:clos-in-pscompl}
  The following statements are true.
  \begin{enumerate}[label=(\alph*)]
  \item \label{item:closure-at-implicative-pres-join}{} If
    $p \in \Sub{X}{\mathsf{M}}$ is implicative then for any family
    \seq{\mu}{i}{I}{} of preneighbourhood systems on $X$:
    \begin{equation}
      \label{eq:closure-at-implicative-pres-join}
      \cls{p}{\bigwedge_{i \in I}\mu_{i}} = \bigvee_{i \in I}\cls{p}{\mu_{i}}.
    \end{equation}

  \item \label{item:closure-at-pt-of-join-is-meet-of-finites}{} For
    any family \seq{\mu}{i}{I}{} of preneighbourhood systems on $X$ and
    $p \in \Sub{X}{\mathsf{M}}$:
    \begin{equation}
      \label{eq:closure-at-pt-of-join-is-meet-of-finites}
      \cls{p}{\bigvee_{i \in I}\mu_i} = \bigwedge_{J \in \subfin{I}}\cls{p}{\mu_J}, \quad\text{ where }\mu_J = \bigvee_{j \in J}\mu_j.
    \end{equation}
    
  \item \label{item:pscompl-case} If \Sub{X}{\mathsf{M}}{} is
    pseudocomplemented then for any preneighbourhood system $\mu$ on
    $X$:
    \begin{equation}
      \label{eq:cond-for-bdd-by-cl}
      \sigma_{X} \neq x < \cls{p}{\mu} \Leftrightarrow p \not< x\text{ and }p^{*} \not\in \mu(x).
    \end{equation}

    Hence: $p \in \mc{\mu}$ (respectively, $p \in \mo{\mu}$) if and
    only if $p^{*} \in \mo{\mu}$ (respectively, $p^{*} \in \mc{\mu}$).

  \item \label{item:dist-atomic-give-additive-cls}{} If
    \Sub{X}{\mathsf{M}}{} is distributive and atom generated then for
    any preneighbourhood space \opair{X}{\mu}{,}
    $p, q \in \Sub{X}{\mathsf{M}}$:
    \begin{equation}
      \label{eq:dist-atomic-give-additive-cls}
      \cls{(p\vee{q})}{\mu} = \cls{p}{\mu} \vee \cls{q}{\mu}.
    \end{equation}
  
  \item \label{item:pnbd-open-gen-atom-gen-give-idemp-cls}{} If $\mu$
    is open generated and \Sub{X}{\mathsf{M}}{} is atom generated then
    \cls{}{\mu}{} is idempotent.

  \end{enumerate}
\end{Thm}

\begin{proof}
  Evidently from Lemma \ref{lem:far-prop}{} and Definition
  \ref{df:closure-closed}{,} for $\mu, \psi \in \pnhd{X}$, if
  $\mu \leq \psi$ then $\cls{}{\psi} \leq \cls{}{\mu}$. Hence for any family
  \seq{\mu}{i}{I}{} of preneighbourhood systems on $X$,
  $\cls{}{\bigvee_{i\in I}\mu_i}\leq\cls{}{\mu_i}\leq\cls{}{\bigwedge_{i\in I}\mu_i}$, entailing
  $\cls{p}{\bigvee_{i\in{I}}\mu_i} \leq \bigwedge_{i \in I}\cls{p}{\mu_i} \leq \bigvee_{i \in
    I}\cls{p}{\mu_i} \leq \cls{p}{\bigwedge_{i\in{I}}\mu_i}$, for any
  $p \in \Sub{X}{\mathsf{M}}$. Since:
  \begin{align*}
    \sigma_{X} \neq x < \cls{p}{\bigwedge_{i\in{I}}\mu_i}
    & \Leftrightarrow p \not< x \notin \mathtt{Far}_{\bigwedge_{i\in{I}}\mu_i}p \\
    & \Leftrightarrow p \not< x \notin \bigcap_{i \in I}\mathtt{Far}_{\mu_i}p
    & \text{ (if $p$ is implicative, \eqref{eq:far-for-meet-of-pnbds})} \\
    & \Leftrightarrow (\exists i \in I)(p \not< x \notin \mathtt{Far}_{\mu_i}p) \\
    & \Leftrightarrow (\exists i \in I)(\sigma_{X} \neq x < \cls{p}{\mu_i}) \\
    & \Rightarrow \sigma_{X} \neq x < \bigvee_{i \in I}\cls{p}{\mu_i},
  \end{align*}
  $\cls{p}{\bigwedge_{i\in{I}}\mu_i} \leq \bigvee_{i\in{I}}\cls{p}{\mu_i}$, proving
  \eqref{eq:closure-at-implicative-pres-join}{.} Since
  $\bigvee_{i \in I}\mu_i = \bigvee_{J \in\subfin{I}}\mu_J$ with
  $\mu_J = \bigvee_{j \in J}\mu_j$ for $J \in \subfin{I}$,
  $\cls{p}{\bigvee_{i\in{I}}\mu_{i}} = \cls{p}{\bigvee_{J \in\subfin{I}}\mu_J} \leq \bigwedge_{J
    \in\subfin{I}}\cls{p}{\mu_J}$, for any
  $p \in \Sub{X}{\mathsf{M}}$; further:
  \begin{align*}
    \sigma_{X} \neq x < \bigwedge_{J\in\subfin{I}}\cls{p}{\mu_J}
    & \Leftrightarrow (\forall J \in \subfin{I})(\sigma_{X} \neq x < \cls{p}{\mu_J}) \\
    & \Leftrightarrow (\forall J \in \subfin{I})(p \not< x \notin \mathtt{Far}_{\mu_J}p) \\
    & \Leftrightarrow p \not< x \notin \mathtt{Far}_{\bigvee_{i\in{I}}\mu_i}p
    & \text{(using \eqref{eq:far-for-joins-of-pnbds}{})} \\
    & \Leftrightarrow \sigma_{X} \neq x < \cls{p}{\bigvee_{i\in{I}}\mu_i},
  \end{align*}
  shows
  $\bigwedge_{J\in\subfin{I}}\cls{p}{\mu_J} \leq \cls{p}{\bigvee_{i\in{I}}\mu_i}$, proving
  \eqref{eq:closure-at-pt-of-join-is-meet-of-finites}{.}  In case when
  $X$ is pseudocomplemented,
  $x \in \mathtt{Far}_{\mu}p \Leftrightarrow (\exists u \in \mu(x))(u \wedge p = \sigma_{X}) \Leftrightarrow (\exists u \in
  \mu(x))(u \leq p^{*}) \Leftrightarrow p^{*} \in \mu(x)$, proving
  \eqref{eq:cond-for-bdd-by-cl}{.} Clearly,
  $p \parallel p^{*} \Leftrightarrow p, p^{*} \neq \sigma_{X}$; if further:
  \begin{enumerate}[label=\roman*.]
  \item $p \in \mc{\mu}$, then
    $\sigma_{X} \neq p^{*} \not< p = \cls{p}{\mu} \Leftrightarrow p < p^{*} \text{ or }p^{*} \in
    \mu(p^*) \Rightarrow p^{*} \in \mu(p^*)$ from \eqref{eq:cond-for-bdd-by-cl} and
    assumption that $p \parallel p^{*}$. Hence $p^{*} \in \mo{\mu}$.

  \item if $p \in \mo{\mu}$ then
    $x \leq p \Leftrightarrow p \in \mu(x) \Rightarrow p^{**} \in \mu(x)$, so that
    \eqref{eq:cond-for-bdd-by-cl}{} implies:
    \begin{equation*}
      \sigma_{X} \neq x < \cls{p^*}{\mu} \Rightarrow p^{*} \not< x\text{ and }x\nleq p.
    \end{equation*}
    Since $p \parallel p^{*}$,
    $p \wedge \cls{p^*}{\mu} = \sigma_{X} \Leftrightarrow \cls{p^*}{\mu} \leq p^{*} \Leftrightarrow p^{*} \in
    \mc{\mu}$.
  \end{enumerate}
  Since
  $\bigl\{\sigma_X,\id{X}\bigr\} \subseteq \mc{\mu} \cap \mo{\mu}$, the statements of the
  second part are trivially true if $p = \sigma_{X}$ or
  $p^{*} = \sigma_{X}$. This completes the proof of
  \ref{item:pscompl-case}{.}  In case \Sub{X}{\mathsf{M}}{} is
  distributive and atom generated then:
  \begin{align*}
    \cls{(p\vee{q})}{\mu} & = \bigvee\bigl(\atom{X}\cap(\mathtt{Far}_{\mu}(p\vee{q}))^c\bigr)
    & \text{(Remark \ref{rem:closure-when-atom-generated})} \\
                     & = \bigvee\bigl(\atom{X}\cap(\mathtt{Far}_{\mu}p\cap\mathtt{Far}_{\mu}q)^c\bigr)
    & \text{(using \eqref{eq:far-for-fin-join-of-subobj})}\\
                     & = \bigvee\biggl(\bigl(\atom{X}\cap(\mathtt{Far}_{\mu}p)^c\bigr)
                       \cup \bigl(\atom{X}\cap(\mathtt{Far}_{\mu}q)^c\bigr)\biggr) \\
                     & = \bigvee\bigl(\atom{X}\cap(\mathtt{Far}_{\mu}p)^c\bigr) \vee
                       \bigvee\bigl(\atom{X}\cap(\mathtt{Far}_{\mu}q)^c\bigr) \\
                     & = \cls{p}{\mu} \vee \cls{q}{\mu}
    & \text{(Remark \ref{rem:closure-when-atom-generated})},
  \end{align*}
  proving \ref{item:dist-atomic-give-additive-cls}{.}  Towards a proof
  of \ref{item:pnbd-open-gen-atom-gen-give-idemp-cls}{,} if
  $\cls{p}{\mu}$ is an atom then $p \in \mc{\mu}$ or $p = \sigma_{X}$, and in
  either case $p = \cls{p}{\mu} = \cls{\cls{p}{\mu}}{\mu}$. If
  $a \in\atom{X}$ and $a < \cls{\cls{p}{\mu}}{\mu}$, then since $\mu$ is open
  generated, $a \not\in \mathtt{Far}_{\mu}{\cls{p}{\mu}}$ implies
  $(\intr{u}{\mu}) \wedge \cls{p}{\mu} \neq \sigma_{X}$ for every
  $u \in \mu(a)$. Choose and fix any $u \in \mu(a)$. Since
  \Sub{X}{\mathsf{M}}{} is atomic, for every atom
  $b \leq (\intr{u}{\mu}) \wedge \cls{p}{\mu}$, $u \in \mu(b)$ and
  $(\intr{u}{\mu}) \wedge p \neq \sigma_{X}$ ($\because$,
  $b < \cls{p}{\mu}$). Since this happens for each $u \in \mu(a)$,
  $a \not\in \mathtt{Far}_{\mu}(p)$, and hence $a < \cls{p}{\mu}$. Since
  \Sub{X}{\mathsf{M}}{} is atom generated,
  \ref{item:pnbd-open-gen-atom-gen-give-idemp-cls}{} stands proved.
\end{proof}

\begin{rem}
  A closure operation which preserves finite joins is called
  \emph{additive} \cite[see conditions (\texttt{AD}) and (\textsf{GR})
  of][\S 2.6]{DikranjanTholen1995}). The condition (\textsf{GR}) of
  \cite{DikranjanTholen1995}{,} however, is already embedded in the
  definition of a closure operation in this paper.
\end{rem}

\begin{rem}
  \label{rem:explains-top-or-locale-or-many}{}
  Theorem \ref{thm:clos-in-pscompl}{} shows if \Sub{X}{\mathsf{M}}{}
  is atom generated then for every neighbourhood system $\mu$ on $X$,
  \cls{}{\mu}{} is idempotent; furthermore, if \Sub{X}{\mathsf{M}}{} is
  distributive then for any neighbourhood system $\mu$, \cls{}{\mu}{} is
  additive (and hence a Kuratowski closure operation).
\end{rem}

\begin{Prop}
  \label{prop:intr-clos-reciprocal}{}
  If \Sub{X}{\mathsf{M}}{} is pseudocomplemented and the
  preneighbourhood system $\mu$ have open interiors then for any
  $p \in \Sub{X}{\mathsf{M}}$:
  \begin{equation}
    \label{eq:intr-clos-reciprocal}
    (\clsm{p}{\mu})^{*} = \intr{p^*}{\mu}.
  \end{equation}
\end{Prop}

\begin{proof}
  From Theorem \ref{thm:clos-in-pscompl}{\ref{item:pscompl-case}{}}
  for any $p \in \Sub{X}{\mathsf{M}}$,
  $p \leq \clsm{p}{\mu} \Rightarrow (\clsm{p}{\mu})^{*} \leq p^{*} \Leftrightarrow (\clsm{p}{\mu})^{*} \leq
  \intr{p^{*}}{\mu}$, since $\mu$ has interiors open. On the other hand,
  if
  $\cls{p}{\mu} = \cls{p}{\mu}\wedge\intr{p^*}{\mu} \Leftrightarrow p \leq \cls{p}{\mu} \leq
  \intr{p^*}{\mu} \leq p^{*} \Leftrightarrow p = \sigma_{X}$. Hence for each
  $p \neq \sigma_{X}$,
  $\cls{p}{\mu} \wedge \intr{p^*}{\mu} < \cls{p}{\mu}$ and hence
  $p \not< \cls{p}{\mu} \wedge \intr{p^*}{\mu} \text{ and } p^{*} \not\in
  \mu(\cls{p}{\mu}\wedge\intr{p^*}{\mu})$, provided
  $\cls{p}{\mu}\wedge\intr{p^*}{\mu}\neq\sigma_X$ (using
  \eqref{eq:cond-for-bdd-by-cl}{}). Since \intr{p^*}{\mu} is
  $\mu$-open,
  $\intr{p^*}{\mu} \in \mu(\intr{p^*}{\mu}) \subseteq \mu(\cls{p}{\mu}\wedge\intr{p^*}{\mu}) \Rightarrow
  p^{*} \in \mu(\cls{p}{\mu}\wedge\intr{p^*}{\mu})$, forcing
  $\cls{p}{\mu}\wedge\intr{p^*}{\mu} = \sigma_{X}$. Since $\mu$ has open interiors an
  use of Theorem \ref{thm:clos-in-pscompl}\ref{item:pscompl-case}
  yields:
  $\cls{p}{\mu} \leq (\intr{p^*}{\mu})^{*} \Leftrightarrow \clsm{p}{\mu} \leq
  (\intr{p^*}{\mu})^{*} \Rightarrow \intr{p^*}{\mu} \leq (\intr{p^*}{\mu})^{**} \leq
  (\clsm{p}{\mu})^{*}$. This completes the proof on observing trivial
  satisfaction of equality for $p = \sigma_{X}$.
\end{proof}

\begin{Cor}
  \label{cor:pscompl-atgen-cl-int-connection}{}
  If \Sub{X}{\mathsf{M}}{} is pseudocomplemented and atom generated
  then for any open generated preneighbourhood system $\mu$ on $X$,
  $(\cls{p}{\mu})^{*} = \intr{p^*}{\mu}$.
\end{Cor}

\begin{rem}\label{rem:reg-clos-reg-open-dual-equiv}{}
  Theorem \ref{thm:clos-in-pscompl}{\ref{item:pscompl-case}{}} yields
  the adjunction on the top row of the diagram:

  \begin{equation}
    \label{eq:reg-open-iso-reg-closed}
    \xymatrix{
      {\mc{\mu}} \ar@<1.2ex>[r]^-{\star} \ar@{<-}@<-1.2ex>[r]_-{\star}
      \ar@{}[r]|-{\bot} & {\opp{\mo{\mu}}} \\
      {\mc{\mu}^{*}} \ar@<1.2ex>[r]^-{\star} \ar@{<-}@<-1.2ex>[r]_-{\star}
      \ar@{>->}[]!<0ex,3.6ex>;[u] \ar@{}[r]|-{\txt{\rotatebox{90}{$\simeq$}}} &
      {\opp{(\mo{\mu}^{*})}} \ar@{>->}[]!<0ex,3.6ex>;[u]      
    }.
  \end{equation}
  The adjunction restricts to an equivalence to the (possibly large)
  sets $\mc{\mu}^{\star} = \bigl\{ p \in \mc{\mu}: p = p^{**} \bigr\}$ of
  \emph{regular closed} subobjects and
  $\mo{\mu}^{\star} = \bigl\{ u \in \mo{\mu}: u = u^{**} \bigr\}$ of
  \emph{regular open} subobjects. Evidently, regular closed subobjects
  are closed under arbitrary meets and finite joins and hence regular
  open subobjects are open under arbitrary joins and finite meets,
  \cite[compare with][Theorem 3.20]{2020}.
\end{rem}


\subsection{}
\label{ssec:closure-vs-pnbd}

Given any object $X$, the symbols $\mathtt{EGM}(X)$,
$\mathtt{CBSMSL}(X)$ abbreviate $\mathtt{EGM}(\Sub{X}{\mathsf{M}})$,
$\mathtt{CBSMSL}(\Sub{X}{\mathsf{M}})$ respectively.

\begin{Thm}
  \label{thm:closures-in-pnbd}{}
  There is an adjunction
  $\xymatrix{ {\mathtt{EGM}(X)} \ar@<1ex>[r]^-{\Phi} \ar@{}[r]|-{\bot}
    \ar@{<-}@<-1ex>[r]_-{\Psi} & {\opp{\pnhd{X}}} }$ with
  $\comp{\Psi}{\Phi} = \id{\mathtt{EGM}(X)}$, where
  $\Phi(c) = \comp{\uparrow_X}{\opp{c}}$ ($c \in \mathtt{EGM}(X)$) and
  $\Psi(\mu) = \bigwedge\mu$ defined by
  $(\bigwedge\mu)(x) = \bigwedge\mu(x)$ ($x \in \Sub{X}{\mathsf{M}}$).
\end{Thm}

\begin{proof}
  Evidently, both $\Phi, \Psi$ are order preserving and
  $\comp{\Psi}{\Phi} = \id{\mathtt{EGM}(X)}$; furthermore, for any
  $c \in \mathtt{EGM}(X)$, $\mu\in\pnhd{X}$:
  \begin{align*}
    \Phi(c) \opp{\leq}\mu
    & \Leftrightarrow \mu\leq\Phi(c) \\
    & \Leftrightarrow \bigl(x\in\Sub{X}{\mathsf{M}} \Rightarrow (\mu(x) \subseteq \Phi(c)(x))\bigr) \\
    & \Leftrightarrow \bigl(x\in\Sub{X}{\mathsf{M}} \Rightarrow (p \in\mu(x) \Rightarrow c(x) \leq p)\bigr) \\
    & \Leftrightarrow \bigl(x\in\Sub{X}{\mathsf{M}} \Rightarrow (c(x) \leq \bigwedge\mu(x))\bigr) \\
    & \Leftrightarrow c \leq \Psi(\mu),
  \end{align*}
  completing the proof.
\end{proof}

\begin{Prop}
  \label{prop:sp-closures-as-pnbd}{}
  Using notation of Theorem \ref{thm:closures-in-pnbd}{,}
  $c \in \mathtt{EGM}(X)$ is idempotent (respectively, idempotent and
  join preserving) if and only if $\Phi(c)$ is a weak neighbourhood
  system (respectively, neighbourhood system) on $X$.
\end{Prop}

\begin{proof}
  $\Phi(c)$ is a weak neighbourhood system if and only if for each
  $x \in \Sub{X}{\mathsf{M}}$:
  \begin{equation*}
    \begin{array}{lrl}
      & u \in \Phi(c)(x) & \Leftrightarrow (\exists v \in \Phi(c)(x))(u \in \Phi(c)(v)) \\
      \Leftrightarrow & u \geq c(x) & \Leftrightarrow (\exists v \geq c(x))(u \geq c(v)) \\
      \Leftrightarrow & c^{2}(x) = c(x),
    \end{array}    
  \end{equation*}
  proving the idempotence of $c$. Further, $\Phi(c)$ is a neighbourhood
  system if and only if for every $A \subseteq \Sub{X}{\mathsf{M}}$:
  \begin{equation*}
    \begin{array}{lrl}
      & \Phi(c)(\bigvee A) & \supseteq \bigcap_{x\in A}\Phi(c)(x) \\
      \Leftrightarrow &  
          \bigl(
          u \geq \bigvee_{x\in A}c(x) & \Rightarrow u \geq c(\bigvee A) 
                             \bigr) \\
      \Leftrightarrow & \bigvee_{x \in A}c(x) & \geq c(\bigvee A)
    \end{array}
  \end{equation*}
  completing the proof.
\end{proof}

\begin{rem}
  A closure operation is said to be \emph{fully additive} if and only
  if it preserves arbitrary non-empty joins \cite[see][\S{2.6}{
    condition (\texttt{FA})}]{DikranjanTholen1995}. Proposition
  \ref{prop:sp-closures-as-pnbd}{} shows fully additive idempotent
  closure operators are special neighbourhood systems on $X$;
  evidently, every neighbourhood system is not of the form $\Phi(c)$ ---
  for instance in $(\Set, \mathsf{Surjections}, \mathsf{Injections})$,
  the usual topology on the real line \Bb{R}{} is not of the form
  $\Phi(c)$, for any $c \in \mathtt{EGM}(\Bb{R})$.
\end{rem}

\begin{rem}
  \label{rem:idemp-cls-cls-pnbd-summary}{}
  Alongside Proposition \ref{prop:cbsmsl-dually-refl-in-egm}, Theorem
  \ref{thm:closures-in-pnbd}{} asserts: the idempotent closure
  operations on \Sub{X}{\mathsf{M}}{,} identified as complete bounded
  sub-$\wedge$-semilattices of the lattice \Sub{X}{\mathsf{M}}{} are
  embedded reflectively inside the complete lattice of grounded
  closure operations on \Sub{X}{\mathsf{M}}, which in turn are
  embedded coreflectively and dually inside the complete lattice of
  preneighbourhood systems on $X$ as some special weak neighbourhood
  systems.  Further, preneighbourhood systems induce closure
  operations yielding the order preserving maps
  \arrp{\cls{}{}}{\clsm{}{}{}}{\opp{\pnhd{X}}}{\mathtt{EGM}(X)}{} with
  $\cls{}{}\leq\clsm{}{}$, $\Phi(\clsm{}{\mu}) \leq \Phi(\cls{}{\mu})$
  and $\Phi(\clsm{}{\mu})$ a weak neighbourhood system on $X$. The
  diagram below summarise this.

  \begin{equation}
    \label{eq:summary-for-pnbd}
    \xymatrixcolsep{6em}
    \xymatrix{
      {\opp{\mathtt{CBSMSL}(X)}} \ar@{>->}[]!<8.4ex,0ex>;[r]_-{\nu}
      \ar@<1.8ex>@{<-}[r]^-{\mathtt{Fix}} \ar@<0.9ex>@{}[r]|-{\bot} &
      {\mathtt{EGM}(X)} \ar@{>->}[]!<5.4ex,0ex>;[r]^-{\Phi}
      \ar@<-1.8ex>@{<-}[r]_-{\Psi} \ar@<-0.9ex>@{}[r]|-{\bot}
      \ar@{<-}@/^{6ex}/[r]_-{\cls{}{}}="l"
      \ar@{<-}@<2.4ex>@/^{6ex}/[r]^-{\clsm{}{}}="u"
      & {\opp{\pnhd{X}}}
      \ar@{} "l";"u"|-{\txt{\rotatebox{90}{$\leq$}}}
    }
  \end{equation}

  In the context
  $(\FinSet, \mathsf{Surjections}, \mathsf{Injections})$ $\Phi$ is an
  isomorphism; the presence of non-discrete Hausdorff topological
  spaces ensure in the context
  $(\Set, \mathsf{Surjections}, \mathsf{Injections})$, $\Phi$ is not
  an isomorphism.
\end{rem}


\subsection{}
\label{ssec:closure-hereditary}
Given any preneighbourhood space \opair{X}{\mu}{,} there are two closure
operations, \cls{}{\mu}{} and \clsm{}{\mu}. The latter is idempotent, and
both of them describe the same closed subobjects. The notion of
\emph{continuity} of morphisms with respect to \cls{}{\mu}{} and
\clsm{}{\mu}{} needs mention.

\begin{Prop}
  \label{prop:cont-wrt-cl-idem-hull}
  Given the internal preneighbourhood spaces \opair{X}{\mu},
  \opair{Y}{\phi} and a morphism \Arr{f}{X}{Y}, consider the statements:
  \begin{enumerate}[label=(\alph*)]
  \item \label{item:cont-wrt-cl} For every
    $p \in \Sub{X}{\mathsf{M}}$,
    $\img{f}{\cls{p}{\mu}} \leq \cls{\img{f}{p}}{\phi}$.

  \item \label{item:cont-wrt-idem-hull} For every
    $p \in \Sub{X}{\mathsf{M}}$,
    $\img{f}{\clsm{p}{\mu}} \leq \clsm{\img{f}{p}}{\phi}$.

  \item \label{item:preimage-pres-closed}{} For every
    $t \in \mc{\phi}$, $\finv{f}{t} \in \mc{\mu}$.
    
  \item \label{item:sm-closed-cont-image-is-sm-cl-cont-image-of-cls}
    For every $p \in \Sub{X}{\mathsf{M}}$,
    $\clsm{\img{f}{p}}{\phi} = \clsm{\img{f}{\cls{p}{\mu}}}{\phi}$.
  \end{enumerate}
  Then: \ref{item:cont-wrt-cl}{} implies
  \ref{item:cont-wrt-idem-hull}; the statements
  \ref{item:cont-wrt-idem-hull},
  \ref{item:sm-closed-cont-image-is-sm-cl-cont-image-of-cls} and
  \ref{item:preimage-pres-closed}{} are equivalent.
\end{Prop}

  \begin{proof}
    Assuming \ref{item:cont-wrt-idem-hull}, for any $t \in \mc{\phi}$ since
    $\img{f}{\clsm{\finv{f}{t}}{\mu}} \leq \clsm{\img{f}{\finv{f}{t}}}{\phi} \leq
    \clsm{t}{\phi} = t$, $\finv{f}{t} \in \mc{\mu}$, proving
    \ref{item:preimage-pres-closed}{.} Conversely, assuming
    \ref{item:preimage-pres-closed}:
    \begin{equation*}
      \begin{aligned}
        \finv{f}{\clsm{\img{f}{p}}{\phi}} & = \finv{f}{\bigwedge \bigl\{ t \in
                                         \mc{\phi}: \img{f}{p} \leq t
                                         \bigr\}} \\
                                       & = \bigwedge \bigl\{ \finv{f}{t}: p \leq \finv{f}{t}, t \in \mc{\phi}
                                         \bigr\}
                                       & \text{ ($\because$, \adjt{\img{f}{}}{\finv{f}{}}{})}\\
                                       & \geq \bigwedge \bigl\{ s \in \mc{\mu}: p \leq s \bigr\} = \clsm{p}{\mu},
      \end{aligned}
    \end{equation*}
    proving \ref{item:cont-wrt-idem-hull}. Further, assuming the
    statement in \ref{item:cont-wrt-idem-hull}{}:
    \begin{equation*}
      \begin{array}{lrll}
        & \img{f}{\clsm{p}{\mu}} & \leq \clsm{\img{f}{p}}{\phi}
        & \text{ (statement in \ref{item:cont-wrt-idem-hull}) } \\[2ex]
        \Leftrightarrow & \clsm{p}{\mu} & \leq \finv{f}{\clsm{\img{f}{p}}{\phi}} \\[2ex]
        \Leftrightarrow & \cls{p}{\mu} & \leq \finv{f}{\clsm{\img{f}{p}}{\phi}}
        & \text{ (since \ref{item:cont-wrt-idem-hull}{} implies
          \ref{item:preimage-pres-closed}{})} \\[2ex]
        \Leftrightarrow & \img{f}{\cls{p}{\mu}} & \leq \clsm{\img{f}{p}}{\phi} \\[2ex]
        \Leftrightarrow & \bigl( t\in\mc{\phi} \Rightarrow (\img{f}{p} \leq t  & \Leftrightarrow \img{f}{\cls{p}{\mu}} \leq t )\bigr)
        \\[2ex] \\
        \Leftrightarrow &  \clsm{\img{f}{\cls{p}{\mu}}}{\phi} & \leq \clsm{\img{f}{p}}{\phi} \\[2ex]
        \Leftrightarrow & \clsm{\img{f}{p}}{\phi} & = \clsm{\img{f}{\cls{p}{\mu}}}{\phi}
        & \text{ ($\because$, $p \leq \cls{p}{\mu}$)} \\[2ex] 
      \end{array}    
    \end{equation*}
    proves the equivalence of the statements
    \ref{item:cont-wrt-idem-hull}{,}
    \ref{item:sm-closed-cont-image-is-sm-cl-cont-image-of-cls} and
    completing the proof of the equivalence of
    \ref{item:cont-wrt-idem-hull}-\ref{item:preimage-pres-closed}{.}
    Finally, assuming \ref{item:cont-wrt-cl} for every
    $t \in \mc{\phi}$ since
    $\img{f}{\cls{\finv{f}{t}}{\mu}} \leq
    \cls{\img{f}{\finv{f}{t}}}{\phi} \leq \cls{t}{\phi} = t$,
    $\finv{f}{t} \in \mc{\mu}$, proving
    \ref{item:preimage-pres-closed}.
  \end{proof}

  \begin{Df}
    \label{df:continuity}{}
    Given the internal preneighbourhood spaces \opair{X}{\mu} and \opair{Y}{\phi}, a morphism \Arr{f}{X}{Y}{} is called \emph{$\mu$-$\phi$ continuous} or simply
  \emph{continuous} if \finv{f}{}{} preserves closed subobjects, i.e.,
  the statement Proposition \ref{prop:cont-wrt-cl-idem-hull}\ref{item:preimage-pres-closed} holds
  good; if $f$ satisfies the statement Proposition \ref{prop:cont-wrt-cl-idem-hull}\ref{item:cont-wrt-cl},
  then $f$ is called \emph{$\mu$-$\phi$ continuous with respect to
    closures} or simply \emph{continuous with respect to
    closures}. 
\end{Df}

\begin{rem}
  If \Arr{c_X}{\Sub{X}{\mathsf{M}}}{\Sub{X}{\mathsf{M}}} and
  \Arr{c_Y}{\Sub{Y}{\mathsf{M}}}{\Sub{Y}{\mathsf{M}}} are both
  monotonic then the adjunction \adjt{\img{f}{}}{\finv{f}{}}{} yields:
  \begin{equation*}
    (\forall x \in \Sub{X}{\mathsf{M}})
    \left(
      \img{f}{c_X(x)} \leq c_{Y}(\img{f}{x})
    \right) \Leftrightarrow (\forall y \in \Sub{Y}{\mathsf{M}})
    \left(
      c_{X}(\finv{f}{y}) \leq \finv{f}{c_Y(y)}
    \right).
  \end{equation*}
  This provides alternative formulations of
  \ref{item:cont-wrt-cl}{} and \ref{item:cont-wrt-idem-hull}
  respectively:
  \begin{align*}
    \cls{\finv{f}{q}}{\mu} \leq \finv{f}{\cls{q}{\phi}}, \qquad\text{ for all }q \in \Sub{X}{\mathsf{M}}, \\ \intertext{and}
    \clsm{\finv{f}{q}}{\mu} \leq \finv{f}{\clsm{q}{\phi}}, \qquad\text{ for all }q \in \Sub{X}{\mathsf{M}}.    
  \end{align*}
\end{rem}

\begin{Thm}
  \label{thm:closure-continuous}
  Given an internal preneighbourhood space \opair{X}{\mu}{} and
  admissible monomorphisms
  $\xymatrix{ {A} \ar@{>->}[]!<2.4ex,0ex>;[r]^-a & {M}
    \ar@{>->}[]!<2.4ex,0ex>;[r]^-m & {X} }$:
  \begin{equation}
    \label{eq:adm-sub-cont}
    \cls{a}{\rest{\mu}{m}} \leq \finv{m}{\cls{(\comp{m}{a})}{\mu}}.
  \end{equation}

  Furthermore:
  \begin{align}
    \label{eq:closure-CC}
    m \in \mc{\mu}, a \in \mc{\rest{\mu}{m}} \Rightarrow \comp{m}{a} \in \mc{\mu} \\ \intertext{and}
    \label{eq:closure-idemp-hereditary}
    m \in \mc{\mu} \Rightarrow 
    \bigl(
    \clsm{a}{\rest{\mu}{m}} = \finv{m}{\clsm{(\comp{m}{a})}{\mu}}
    \bigr).
  \end{align}

  Finally, if \Sub{X}{\mathsf{M}}{} is atom generated or
  \Arr{\finv{m}{}}{\Sub{X}{\mathsf{M}}}{\Sub{M}{\mathsf{M}}}{}
  preserve joins then:
  \begin{equation}
    \label{eq:closure-hereditary}
    \cls{a}{\rest{\mu}{m}} = \finv{m}{\cls{(\comp{m}{a})}{\mu}}.
  \end{equation}
\end{Thm}

\begin{proof}
  Using the adjunction \adjt{\img{m}{}}{\finv{m}{}}{,} for any
  $x \in \Sub{X}{\mathsf{M}}$,
  $x \geq \comp{m}{a} \Leftrightarrow \finv{m}{x} \geq a$, i.e.,
  $x \not> \comp{m}{a} \Leftrightarrow \finv{m}{x} \not> a$ ($\because, m \in \mathsf{M}$). Also, since
  $\finv{m}{(u\wedge\comp{m}{a})} = \finv{m}{u} \wedge a$, implying
  $\comp{m}{(a\wedge\finv{m}{u})} =
  \comp{m}{\left(\finv{m}{(u\wedge\comp{m}{a})}\right)} = m \wedge u \wedge
  \comp{m}{a} = u \wedge \comp{m}{a}$, the computations:
  \begin{align*}
    x \in \mathtt{Far}_{\mu}(\comp{m}{a})
    & \Leftrightarrow (\exists u \in \mu(x))(u \wedge \comp{m}{a} = \sigma_{X}) \\
    & \Leftrightarrow (\exists u \in \mu(x))(\comp{m}{(a\wedge\finv{m}{u})} = \sigma_{X}) \\
    & \Leftrightarrow (\exists u \in \mu(x)(a \wedge\finv{m}{u} = \sigma_{M}))
    & \text{ ($\because, m \in \mathsf{M}$)}\\
    & \Leftrightarrow (\exists v \in \rest{\mu}{m}(\finv{m}{x}))(a \wedge v = \sigma_{M})
    & \text{ (using \eqref{eq:rest-pnbd})}\\
    & \Leftrightarrow \finv{m}{x} \in \mathtt{Far}_{\rest{\mu}{m}}(a),
  \end{align*}
  show:
  $x \in \mathtt{Far}_{\mu}(\comp{m}{a}) \Leftrightarrow \finv{m}{x} \in
  \mathtt{Far}_{\rest{\mu}{m}}(a)$.

  Therefore:
  \begin{align*}
    \img{m}{\cls{a}{\rest{\mu}{m}}}
    & = \img{m}{\bigvee
      \bigl\{
      x \in \Sub{M}{\mathsf{M}}: a \not< x \not\in \mathtt{Far}_{\rest{\mu}{m}}(a)
      \bigr\}
      } \\
    & = \bigvee
      \bigl\{
      \img{m}{x}: x \in \Sub{M}{\mathsf{M}},
      a \not< \finv{m}{(\comp{m}{x})} \not\in \mathtt{Far}_{\rest{\mu}{m}}(a)
      \bigr\} \\
    & = \bigvee
      \bigl\{
      \comp{m}{x}: x \in \Sub{M}{\mathsf{M}}, \comp{m}{a} \not< \comp{m}{x}
      \not\in \mathtt{Far}_{\mu}{(\comp{m}{a})}
      \bigr\} \\
    & \leq \bigvee
      \bigl\{
      y \in \Sub{X}{\mathsf{M}}: \comp{m}{a}\not< y \not\in
      \mathtt{Far}_{\mu}(\comp{m}{a})
      \bigr\} \\
    & = \cls{(\comp{m}{a})}{\mu},
  \end{align*}
  proving
  $\comp{m}{\cls{a}{\rest{\mu}{m}}} \leq \cls{(\comp{m}{a})}{\mu} \Leftrightarrow
  \cls{a}{\rest{\mu}{m}} \leq \finv{m}{\cls{(\comp{m}{a})}{\mu}}$.

  The equations \eqref{eq:closure-CC} \&
  \eqref{eq:closure-idemp-hereditary} are trivially true when $a$ is
  an isomorphism; hence it is enough to prove them for $a$ not an
  isomorphism.

  For $m \in \mc{\mu}$ and $\id{M} \neq a \in \mc{\rest{\mu}{m}}$:
  \begin{align*}
    \cls{(\comp{m}{a})}{\mu}
    & = \bigvee
      \bigl\{
      x \in \Sub{X}{\mathsf{M}}: \comp{m}{a}\not< x \not\in
      \mathtt{Far}_{\mu}(\comp{m}{a})
      \bigr\} \\
    & = \bigvee
      \bigl\{
      x \in \Sub{X}{\mathsf{M}}: a \not<
      \finv{m}{x} \not\in \mathtt{Far}_{\rest{\mu}{m}}(a)
      \bigr\} \\
    & = \bigvee
      \bigl\{
      x \in \Sub{X}{\mathsf{M}}: \finv{m}{x} \leq a
      \bigr\} \\
    & \text{ (since $a \in \mc{\rest{\mu}{m}}$)} \\
    & = \bigvee
      \bigl\{
      x \in \Sub{X}{\mathsf{M}}: \sigma_{M} \neq \finv{m}{x} \leq a
      \bigr\} \vee \bigvee
      \bigl\{
      x \in \Sub{X}{\mathsf{M}}: \finv{m}{x}=\sigma_M
      \bigr\} \\
    & = p \vee q,
  \end{align*}
  where
  $p = \bigvee \bigl\{ x \in \Sub{X}{\mathsf{M}}: \sigma_{M} \neq \finv{m}{x} \leq a
  \bigr\}$ and
  $q = \bigvee \bigl\{ x \in \Sub{X}{\mathsf{M}}: \finv{m}{x}=\sigma_M
  \bigr\}$. For $x \in \Sub{X}{\mathsf{M}}$ with
  $x < \cls{(\comp{m}{a})}{\mu}$ such that
  $\finv{m}{x}=\sigma_M \Leftrightarrow x \wedge m = \sigma_{X}$, since
  $x \not\in \mathtt{Far}_{\mu}(\comp{m}{a})$,
  $x \not\in \mathtt{Far}_{\mu}(m)$ ($\because, \comp{m}{a} < m$ and
  Lemma \ref{lem:far-prop}), and hence $x \geq m$ or $x \leq m$
  ($\because, m \in \mc{\mu}$). Since $x \not>\comp{m}{a}$,
  $x \not\geq m$ forces $x \leq m$. Hence there exists only one $x$
  contributing to the join for $q$, namely $x = \sigma_{X}$, i.e.,
  $q = \sigma_{X}$. On the other hand, using Lemma
  \ref{lem:non-zero-meets-for-obj-not-far} and $m \in \mc{\mu}${,} for
  each $x$ contributing to the join for $p$, $x \leq \comp{m}{a}$. Thus:
  \begin{equation*}
    \cls{(\comp{m}{a})}{\mu} = p = \bigvee
    \bigl\{
    x \in \Sub{X}{\mathsf{M}}: x \leq \comp{m}{a} 
    \bigr\} = \comp{m}{a}
  \end{equation*}
  
  yielding $\comp{m}{a} \in \mc{\mu}$; also:
  \begin{align*}
    \finv{m}{\clsm{(\comp{m}{a})}{\mu}}
    & = \finv{m}{\bigwedge
      \bigl\{
      s \in \mc{\mu}: \comp{m}{a} \leq s
      \bigr\}} \\
    & = \bigwedge
      \bigl\{
      \finv{m}{s}: a \leq \finv{m}{s}, s \in \mc{\mu}
      \bigr\} \\
    & = \bigwedge
      \bigl\{
      t \in \mc{\rest{\mu}{m}}: a \leq t
      \bigr\}
    & \text{ (using \eqref{eq:closure-CC}{}) } \\
    & = \clsm{a}{\rest{\mu}{m}},
  \end{align*}
  yielding \eqref{eq:closure-idemp-hereditary}{.}

  If \Arr{\finv{m}{}}{\Sub{X}{\mathsf{M}}}{\Sub{M}{\mathsf{M}}}{}
  preserve arbitrary joins then:
  \begin{align*}
    \finv{m}{\cls{(\comp{m}{a})}{\mu}}
    & = \finv{m}{\bigvee
      \bigl\{
      x \in \Sub{X}{\mathsf{M}}: \comp{m}{a} \not< x \not\in
      \mathtt{Far}_{\mu}(\comp{m}{a})
      \bigr\}} \\
    & = \bigvee
      \bigl\{
      \finv{m}{x}: x \in \Sub{X}{\mathsf{M}}, a \not< \finv{m}{x} \not\in
      \mathtt{Far}_{\rest{\mu}{m}}(a)
      \bigr\} \\
    & \leq \bigvee
      \bigl\{
      y \in \Sub{M}{\mathsf{M}}: a \not< y \not\in \mathtt{Far}_{\rest{\mu}{m}}(a)
      \bigr\} = \cls{a}{\rest{\mu}{m}},
  \end{align*}
  proving \eqref{eq:closure-hereditary}{} in this case. On the other
  hand if \Sub{X}{\mathsf{M}}{} is atom generated then for
  an atom $b$, $b \leq m \wedge \cls{(\comp{m}{a})}{\mu}$ if and only if
  $b \leq m$ and
  $b \not\in \mathtt{Far}_{\mu}(\comp{m}{a}) \Leftrightarrow \finv{m}{b} \not\in
  \mathtt{Far}_{\rest{\mu}{m}}(a)$; since $\finv{m}{b}$ is also an atom,
  it implies
  $\finv{m}{b} \leq \cls{a}{\rest{\mu}{m}} \Leftrightarrow m \wedge b = b \leq
  \comp{m}{\cls{a}{\rest{\mu}{m}}}$. Hence
  $m \wedge \cls{(\comp{m}{a})}{\mu} \leq \comp{m}{\cls{a}{\rest{\mu}{m}}} \Leftrightarrow
  \finv{m}{\cls{(\comp{m}{a})}{\mu}} \leq \cls{a}{\rest{\mu}{m}}$, proving
  \eqref{eq:closure-hereditary}{} in this case also. This completes
  the proof.
\end{proof}

\begin{Lemma}
  \label{lem:non-zero-meets-for-obj-not-far}{}
  Given the admissible monomorphisms
  $\xymatrix{ {A} \ar@{>->}[]!<2.4ex,0ex>;[r]^-a & {M}
    \ar@{>->}[]!<2.4ex,0ex>;[r]^-m & {X} }$ with $a \neq \id{M}$, if
  $\sigma_{M} \neq (x \wedge m) \leq \comp{m}{a}$ then $x < \cls{m}{\mu}$.
\end{Lemma}

\begin{proof}
  If $x > m$ then $m = x \wedge m \leq \comp{m}{a} \leq m$ implies $a$ is an
  isomorphism; hence if $a \neq \id{M}$ then $x \not> m$. Also,
  $x \wedge m \leq \comp{m}{a} \Leftrightarrow \finv{m}{x} \leq a \Rightarrow \finv{m}{x} \not\in
  \mathtt{Far}_{\rest{\mu}{m}}(a) \Leftrightarrow x \not\in\mathtt{Far}_{\mu}(\comp{m}{a})
  \Rightarrow x \not\in \mathtt{Far}_{\mu}(m)$. Thus, the conditions imply:
  $m \not< x \not\in \mathtt{Far}_{\mu}(m) \Leftrightarrow x < \cls{m}{\mu}$.
\end{proof}

\begin{Df}[{\cite[see][condition
  \texttt{(HE)}, \S{2.5}]{DikranjanTholen1995}}]
\label{df:hereditary}{}
Given an internal preneighbourhood space \opair{X}{\mu} and a $m \in \Sub{X}{\mathsf{M}}$, \cls{}{\mu}{} is \emph{hereditary for $m$} if for all $a \in \Sub{M}{\mathsf{M}}$:
\begin{equation}
  \label{eq:closure-hereditary-df}
  \cls{a}{\rest{\mu}{m}} = \finv{m}{\cls{(\comp{m}{a})}{\mu}}.
\end{equation}

Similarly, \clsm{}{\mu}{} is \emph{hereditary for $m$} if for all $a \in \Sub{M}{\mathsf{M}}$:
\begin{equation}
  \label{eq:idemp-hull-hereditary-df}
  \clsm{a}{\rest{\mu}{m}} = \finv{m}{\clsm{(\comp{m}{a})}{\mu}}.
\end{equation}

For a subset $\mathfrak{a}\subseteq\Sub{X}{\mathsf{M}}$, \cls{}{\mu}{} (respectively, \clsm{}{\mu}{}) is \emph{\eulf{a}-hereditary{}} if \cls{}{\mu}{} (respectively, \clsm{}{\mu}{}) is hereditary for each $m \in \eulf{a}$; \textsf{M}-hereditary is shortened to \emph{hereditary}.
\end{Df}

\begin{rem}
  In terms of Definition \ref{df:hereditary}{,} equation \eqref{eq:closure-hereditary}{} suggests in the special case when \Sub{X}{\mathsf{M}}{} is atom generated, for every preneighbourhood system $\mu$ on $X$, \cls{}{\mu}{} is hereditary; in general, if for an admissible subobject $m$, \finv{m}{}{} preserve joins then \cls{}{\mu}{} is hereditary for $m$. For every admissible monomorphism $m \in \Sub{X}{\mathsf{M}}$, \finv{m}{}{} preserve joins if and only if \Sub{X}{\mathsf{M}}{} is a frame (\cite[see][]{FreydScedrov1990} or \cite[see][Theorem 2.11(a)]{2020}{}). Thus, \cls{}{\mu}{} is hereditary for any preneighbourhood system $\mu$ on $X$ if \Sub{X}{\mathsf{M}}{} is a frame.

  However, in general from \eqref{eq:closure-idemp-hereditary}{,} for every preneighbourhood system $\mu$ on $X$, \clsm{}{\mu}{} is \mc{\mu}{-}hereditary. 
\end{rem}

\begin{rem}\label{rem:closed-embed}{}
  Equations \eqref{eq:adm-sub-cont} and \eqref{eq:closure-CC} together
  suggest:
  \begin{equation*}
    m \in \mc{\mu} \Rightarrow \bigl( \comp{m}{a} \in \mc{\mu} \Leftrightarrow a \in \mc{\rest{\mu}{m}}
  \bigr),
\end{equation*}
i.e., the closed subobjects of a closed subobject $M$ are precisely those which are closed in $X$.  
\end{rem}

\begin{rem}
  Theorem \ref{thm:closure-continuous}{} ensures each admissible
  monomorphism is continuous with respect to closures and satisfies
  the condition \texttt{(CC)} \cite[see][condition \texttt{(CC)},
  \S{2.4}]{DikranjanTholen1995}. 
\end{rem}


\subsection{}
\label{ssec:cont-pren-morph}

Theorem \ref{thm:closure-continuous}{} ensure every admissible
monomorphism is continuous with respect to closures and hence
continuous (Prop \ref{prop:cont-wrt-cl-idem-hull}{}).

\begin{Prop}
  \label{prop:cont-suffices-weak}{}
  Given any preneighbourhood morphism
  \Arr{f}{\opair{X}{\mu}}{\opair{Y}{\phi}},
  if $f = \comp{m}{h}$ where \clsm{}{\mu}{} is hereditary for $m$ and
  $h$ is a preneighbourhood morphism then $f$ is continuous if and
  only if $h$ is continuous.
\end{Prop}

\begin{proof}
  If
  $\xymatrix{ {\opair{X}{\mu}} \ar@/_{3ex}/[rr]_-f \ar[r]^-h &
    {\opair{I}{\rest{\phi}{m}}} \ar@{>->}[]!<6.6ex,0ex>;[r]^-m &
    {\opair{Y}{\phi}} }$ with \clsm{}{\mu}{} hereditary for $m$ then:
  \begin{equation*}
    \begin{array}{lrlr}
      & \img{f}{\clsm{p}{\mu}} & \leq \clsm{\img{f}{p}}{\phi} & \\[1.2ex]
      \Leftrightarrow & \comp{m}{\img{h}{\clsm{p}{\mu}}} & \leq \clsm{(\comp{m}{\img{h}{p}})}{\phi} & \\[1.2ex]
      \Leftrightarrow & \img{h}{\clsm{p}{\mu}} & \leq
                                 \finv{m}{\clsm{(\comp{m}{\img{h}{p}})}{\phi}} & \\[1.2ex]
      \Leftrightarrow & \img{h}{\clsm{p}{\mu}} & \leq \clsm{\img{h}{p}}{\rest{\phi}{m}}
                                                      & \text{ (since \clsm{}{\mu}{} is $m$-hereditary)},
    \end{array}
  \end{equation*}
  completing the proof.
\end{proof}

Using Definition \ref{df:dense-mor} \& Remark
\ref{rem:dense-is-closure-dense}:
\begin{Cor}
  \label{cor:cont-iff-dense-are-cont}
  In the context of Proposition \ref{prop:cont-suffices-weak}{,} every
  preneighbourhood morphism is continuous if and only if every dense
  morphism is continuous.
\end{Cor}

\begin{Lemma}
  \label{lem:far-of-dom-cod}{}
  Consider internal preneighbourhood spaces \opair{X}{\mu}{,}
  \opair{Y}{\phi}, a morphism \Arr{f}{X}{Y}, admissible subobjects
  $x, p \in \Sub{X}{\mathsf{M}}$.

  \begin{enumerate}[label=(\alph*)]
  \item \label{item:far-cod-to-dom}{} If $f$ reflects zero and
    $\mu \supseteq \invfil{f}{\phi\img{f}{}}$ then:
    \begin{equation}
      \label{eq:far-cod-to-dom}
      \img{f}{x} \in \mathtt{Far}_{\phi}(\img{f}{p}) \Rightarrow x \in \mathtt{Far}_{\mu}p.
    \end{equation}

  \item \label{item:far-dom-to-cod}{} If $f$ is a Frobenius morphism
    and $\mu \subseteq \invfil{f}{\phi\img{f}{}}$ then:
    \begin{equation}
      \label{eq:far-dom-to-cod}
      x \in \mathtt{Far}_{\mu}p \Rightarrow \img{f}{x} \in \mathtt{Far}_{\phi}(\img{f}{p}).
    \end{equation}
  \end{enumerate}
  
\end{Lemma}

\begin{proof}
  If $\img{f}{x} \in \mathtt{Far}_{\mu}(\img{f}{p})$ then there exists a
  $v \in \phi(\img{f}{x})$ with $v \wedge \img{f}{p} = \sigma_{Y}$; if
  $f$ reflects zero then
  $\sigma_{X} = \finv{f}{(v\wedge\img{f}{p})} \geq p \wedge \finv{f}{v} \Rightarrow p \wedge
  \finv{f}{v} = \sigma_{X}$. Hence if
  $\mu \supseteq \invfil{f}{\phi\img{f}{}}$ then
  $\finv{f}{v} \in \mu(x)$ proving \ref{item:far-cod-to-dom}{.}  On the
  other hand, if $x \in \mathtt{Far}_{\mu}p$ then there exists a
  $u \in \mu(x)$ such that $u \wedge p = \sigma_{X}$. Since
  $\mu \subseteq \invfil{f}{\phi\img{f}{}}$, there exists a
  $v \in \phi(\img{f}{x})$ such that $\finv{f}{v} \leq u$ entailing
  $p \wedge \finv{f}{v} = \sigma_{X}$ for some
  $v \in \phi(\img{f}{x})$. Since $f$ is a Frobenius morphism,
  $\sigma_{Y} = \img{f}{(p\wedge\finv{f}{v})} = v \wedge \img{f}{p}$, implying
  $\img{f}{x} \in \mathtt{Far}_{\phi}(\img{f}{p})$. This completes the
  proof.
\end{proof}

\begin{rem}
  \label{rem:obstruction-again}
  Lemma \ref{lem:far-of-dom-cod} illustrates the obstruction in
  establishing continuity or continuity with respect to
  closures. Thus, in context of Lemma \ref{lem:far-of-dom-cod}{,} if
  $t \in \mc{\phi}$ and $f$ is a zero reflecting preneighbourhood morphism
  then
  $x \notin \mathtt{Far}_{\mu}(\finv{f}{t}) \Rightarrow \img{f}{x} \notin
  \mathtt{Far}_{\phi}(\img{f}{\finv{f}{t}}) $ and hence
  $\img{f}{x} \notin \mathtt{Far}_{\phi}{t}$. Since $t \in \mc{\phi}$,
  $\img{f}{x} \leq t \Leftrightarrow x \leq \finv{f}{t}$ or else
  $\img{f}{x} \geq t$. Consequently, \finv{f}{t}{} would be closed if
  $\img{f}{x} \geq t \Rightarrow x \geq \finv{f}{t}$. With regards to continuity with
  respect to closures:
  \begin{align*}
    x < \cls{p}{\mu} & \Leftrightarrow p \not< x \notin \mathtt{Far}_{\mu}p \\
                   & \Rightarrow p \not< x \text{ and }\img{f}{x} \notin \mathtt{Far}_{\phi}(\img{f}{p})
                   & \text{ (if $\mu \supseteq \invfil{f}{\phi\img{f}{}}$ and $f$ reflects zero)};
  \end{align*}
  hence if $\img{f}{p} < \img{f}{x} \Rightarrow p < x$ then
  $x < \cls{p}{\mu} \Rightarrow \img{f}{x} < \cls{\img{f}{p}}{\phi}$
  entailing continuity with respect to closures.
\end{rem}

\begin{Cor}
  \label{cor:cont-examples}{}
  Given \Arr{f}{\opair{X}{\mu}}{\opair{Y}{\phi}}, the following
  statements hold.
  \begin{enumerate}[label=(\alph*)]
  \item \label{item:refl-0-adm-mono-like-is-cont}{} If
    $\comp{\finv{f}{}}{\img{f}{}} = \id{\Sub{X}{\mathsf{M}}}$ then $f$
    is continuous with respect to closures.

  \item \label{item:refl-0-atomic-atom-gen-pres-atom-is-cont}{} If
    \Sub{X}{\mathsf{M}}{} is atom generated, $f$ reflects zero and
    preserve atoms then $f$ is continuous with respect to closures.
  \end{enumerate}
\end{Cor}

\begin{rem}
  Every admissible monomorphism satisfy condition of Corollary
  \ref{cor:cont-examples}{\ref{item:refl-0-adm-mono-like-is-cont}{}};
  however the property do not characterise admissible monomorphisms
  --- in the context $(\Top, \mathsf{Epi}, \mathsf{ExtMon})$ injective
  continuous maps may not be extremal monomorphisms and yet satisfy
  the condition.
\end{rem}

\begin{rem}
  With regards to preservation of atoms, every Frobenius morphism
  preserve atoms: if \Arr{f}{X}{Y}{} is a Frobenius morphism and
  $a \in \atom{X}$, then for each $y \in \Sub{Y}{\mathsf{M}}$ either
  $a \leq \finv{f}{y} \Leftrightarrow \img{f}{a} \leq y$ or else
  $a \wedge \finv{f}{y} = \sigma_{X} \Rightarrow \sigma_{Y} =
  \img{f}{(a\wedge\finv{f}{y})} = y \wedge \img{f}{a}$, completing the
  proof. Hence, from Corollary
  \ref{cor:cont-examples}\ref{item:refl-0-atomic-atom-gen-pres-atom-is-cont}{,}
  if \Sub{X}{\mathsf{M}}{} is atom generated then every reflecting
  zero Frobenius preneighbourhood morphism with \opair{X}{\mu}{} as
  domain is continuous with respect to closures.
\end{rem}


\subsection{}
\label{ssec:closure-ex}

As observed \emph{continuity} for morphisms with respect to induced
closure operations is not automatic, even for preneighbourhood
morphisms. This section illustrate its presence in many familiar
contexts, as well as exhibit instances of several properties of
closure operations from a preneighbourhood system discussed so
far. However, before embarking on the examples, let for
$\mu \in \pnhd{X}$ and $\phi \in \pnhd{Y}$:
\begin{align}
  \label{eq:sets-of-clos-cont-pnhd-maps}
  \mathtt{CC}\opair{\mu}{\phi}
  & = 
    \bigl\{
    \Arr{f}{X}{Y}: \mu \supseteq \invfil{f}{\phi\img{f}{}}\text{ and }(\forall p \in \Sub{X}{\mathsf{M}})(\img{f}{\cls{p}{\mu}} \leq \cls{\img{f}{p}}{\phi})
    \bigr\}, \\ \intertext{and}
  \label{eq:sets-of-cont-pnhd-maps}
  \mathtt{C}\opair{\mu}{\phi}
  & = 
    \bigl\{
    \Arr{f}{X}{Y}: \mu \supseteq \invfil{f}{\phi\img{f}{}}\text{ and }(\forall p \in \Sub{X}{\mathsf{M}})(\img{f}{\clsm{p}{\mu}} \leq \clsm{\img{f}{p}}{\phi})
    \bigr\}.
\end{align}
Evidently,
\arrp{\mathtt{C}}{\mathtt{CC}}{\pnhd{X}\times\opp{\pnhd{Y}}}{2^{\Homo{\Bb{A}}{X}{Y}}}{}
are both order preserving maps with
$\mathtt{CC}\opair{\mu}{\phi} \subseteq \mathtt{C}\opair{\mu}{\phi}$.

\begin{subsubsection}{}
  \label{sssec:sets-with-structure}
  
  Let \CAL{A}{} be any of the contexts:
  $(\FinSet, \mathsf{Surjections}, \mathsf{Injections})$,
  $(\Set, \mathsf{Surjections}, \newline \mathsf{Injections})$,
  $(\Top, \mathsf{Epi}, \mathsf{ExtMon})$. In each of them for each
  object $X$, \Sub{X}{\mathsf{M}} is distributive, complemented, atom
  generated, morphisms preserve atoms and the contexts are admissibly
  quasi-pointed with strict initial object. Hence for every
  preneighbourhood system $\mu$ on $X$, \cls{}{\mu}{} is additive (Theorem
  \ref{thm:clos-in-pscompl}\ref{item:dist-atomic-give-additive-cls}{}),
  hereditary (Theorem \ref{thm:closure-continuous}{,} equation
  \eqref{eq:closure-hereditary}{}),
  $p \in \mc{\mu} \Rightarrow p^{*} \in \mo{\mu}$ (respectively,
  $p \in \mo{\mu} \Rightarrow p^{*} \in \mc{\mu}$) (Theorem
  \ref{thm:clos-in-pscompl}{\ref{item:pscompl-case}{}}) and each
  preneighbourhood morphism is continuous (Corollary
  \ref{cor:cont-examples}\ref{item:refl-0-atomic-atom-gen-pres-atom-is-cont}{,}
  Theorem
  \ref{thm:reflzero-prop}\ref{item:morph-refl-zero=str-initial}{,}
  Remark \ref{rem:adm-quasi-pt}{}); if further $\mu$ has open interiors
  then $(\clsm{p}{\mu})^{*} = \intr{p^*}{\mu}$ (Proposition
  \ref{prop:intr-clos-reciprocal}{}). Moreover, if $\mu$ is open
  generated then \cls{}{\mu}{} is idempotent (Theorem
  \ref{thm:clos-in-pscompl}{\ref{item:pnbd-open-gen-atom-gen-give-idemp-cls}{}}). Thus
  for each neighbourhood system $\mu$ on $X$, \cls{}{\mu}{} is a
  hereditary Kuratowski closure operation yielding internal
  topologies.

\end{subsubsection}{}

\begin{subsubsection}{}
  \label{sssec:locales}{}
  The context $(\Loc, \mathsf{Epi}, \mathsf{RegMono})$ has each of its
  \Sub{X}{\mathsf{M}}{} a coframe and the context is admissibly
  quasi-pointed with the initial object strict; however, unlike
  contexts in \S\ref{sssec:sets-with-structure}{,} preimages do not
  preserve arbitrary joins and images do not preserve
  atoms. Furthermore, each subobject (also called \emph{sublocale})
  $S \subseteq X$ is generated by \emph{principal sublocales}:
  $S = \bigvee_{a \in S}[a]$, where
  $[a] = \bigl\{ \implyr{t}{a}: t \in X \bigr\}$ is the smallest
  sublocale containing $a$ \cite[see][\S
  III.10.2{}]{PicadoPultr2012}. Hence, for every preneighbourhood
  system $\mu$ on a locale $X$, \cls{}{\mu}{} is hereditary (Theorem
  \ref{thm:closure-continuous}{,} equation
  \eqref{eq:closure-hereditary}{}).

  Consider the $T$-neighbourhood system for a locale $X$ (see Example
  \ref{ex:locales}{}).  The closed sublocale for $a \in X$ is
  $\atleast{a} = \bigl\{ x \in X: x \geq a \bigr\}$ is the complement of
  $\open{a} = \bigl\{ \implyr{a}{x}: x \in X \bigr\} = \bigl\{ x \in X: x
  = \implyr{a}{x} \bigr\}$ in the lattice \Sub{X}{\mathsf{M}} of all
  sublocales of $X$ \cite[see][Proposition
  6.1.3]{PicadoPultr2012}. Since for any $x \in X$:
  $1 \neq x \in \open{a} \Leftrightarrow x = \implyr{a}{x} \neq 1 \Leftrightarrow a \nleq x$, for any
  sublocale $S$:
  \begin{align*}
    x \in \cls{S}{\tau_X}
    & \Leftrightarrow \gen{x} \subseteq \cls{S}{\tau_X} \\
    & \Leftrightarrow 
      \left(
      T \in \tau_X(\gen{x}) \Rightarrow T \cap S \neq \{1\}
      \right) \\
    & \Leftrightarrow 
      \left(
      \gen{x} \subseteq \open{a} \Rightarrow \open{a} \cap S \neq \{1\}
      \right) \\
    & \Leftrightarrow
      \left(
      \open{a} \cap S = \{1\} \Rightarrow \gen{x}\nsubseteq \open{a}
      \right) \\
    & \Leftrightarrow
      \left(
      S \subseteq \atleast{a} \Rightarrow x \geq a
      \right) \\
    & \Leftrightarrow x \geq \bigwedge S.
  \end{align*}
  This proves $\cls{S}{\tau_X} = \atleast{(\bigwedge S)}$, the usual localic
  closure of $S$ \cite[see][\S III.8{}]{PicadoPultr2012}; In
  particular, \cls{}{\tau_X} is additive and idempotent, i.e., a
  hereditary Kuratowski closure operation; since \Sub{X}{\mathsf{M}}{}
  is a coframe, this is not an internal topology in general. Further,
  for every localic map \Arr{f}{X}{Y}{,} since
  $\img{f}{\cls{S}{\tau_X}} \subseteq \cls{\img{f}{S}}{\tau_Y}$
  ($S \in \Sub{X}{\mathsf{RegMono}}$)
  \cite[see][\S{III.8.4}{}]{PicadoPultr2012}, the preneighbourhood
  morphism \Arr{f}{\opair{X}{\tau_X}}{\opair{Y}{\tau_Y}}{} is continuous,
  and for preneighbourhoods $\mu \supseteq \tau_X$,
  $\tau_Y \supseteq \phi$, \Arr{f}{\opair{X}{\mu}}{\opair{Y}{\phi}}{} is continuous.
\end{subsubsection}{}

\begin{subsubsection}{}
  \label{sssec:groups}{}
  In the context $(\Grp, \mathsf{RegEpi}, \mathsf{Mono})$, for a group
  $X$ consider the neighbourhood system $\nu_X$ of Example
  \ref{ex:groups}; in fact $\nu_X = \Phi(\mathtt{ncl}_X)$, where
  $\mathtt{ncl}_X$ is the normal closure of a subgroup $G$ of
  $X$. Since $\mathtt{ncl}_X$ is an idempotent and join preserving
  closure operation, $\Phi(\mathtt{ncl}_X)$ is a neighbourhood system on
  $X$, \cls{}{\Phi(\mathtt{ncl}_X)}{} is idempotent,
  \mo{\Phi(\mathtt{ncl}_X)} is precisely the set of all normal subgroups
  of $X$, for any subgroup $A$:
  \begin{align*}
    x \in \cls{A}{\Phi(\mathtt{ncl}_X)}
    & \Leftrightarrow \gen{x} \in \cls{A}{\Phi(\mathtt{ncl}_X)} \\
    & \Leftrightarrow \mathtt{ncl}_X(x) \cap A \neq \{0\}, \\
  \end{align*}
  where \gen{x}{} is the cyclic group generated by $x$,
  $A \in \mc{\Phi(\mathtt{ncl}_X)}$ if and only if for any $x \in X$:
  \begin{equation*}
    \mathtt{ncl}_X(x) \cap A \neq \{0\} \Rightarrow x \in A,
  \end{equation*}
  $N \trianglelefteq X \Rightarrow \cls{N}{\Phi(\mathtt{ncl}_X)} \trianglelefteq X$,
  $\intr{A}{\Phi(\mathtt{ncl}_X)}$ is the normal core of $A$ (i.e., the
  largest normal subgroup contained in a subgroup), for each group
  homomorphism \Arr{f}{X}{Y} the preneighbourhood morphism
  \Arr{f}{\opair{X}{\Phi(\mathtt{ncl}_X)}}{\opair{Y}{\Phi(\mathtt{ncl}_Y)}}{}
  is continuous and every \Arr{f}{\opair{X}{\mu}}{\opair{Y}{\phi}}{} where
  $\mu \supseteq \mathtt{ncl}_{X}$ and $\phi \subseteq \mathtt{ncl}_{Y}$ is continuous.
\end{subsubsection}{}

\begin{subsubsection}{}
  \label{sssec:rngs}
  In the context $(\CRng, \mathsf{RegEpi}, \mathsf{Mono})$ for a ring
  $X$ consider the neighbourhood system $\iota_X$ of Example
  \ref{ex:rngs}; in fact $\iota_X = \Phi(\mathtt{idl}_X)$, where
  $\mathtt{idl}_X$ is the ideal closure of a subring of $X$. Since
  $\mathtt{idl}_X$ is an idempotent and join preserving closure
  operation, $\Phi(\mathtt{idl}_X)$ is a neighbourhood system,
  \cls{}{\Phi(\mathtt{idl}_X)}{} is idempotent,
  $\mo{\Phi(\mathtt{idl}_X)} = \Idl{X}$ the set of all ideals of the ring
  $X$, for any subring $A$:
  \begin{align*}
    x \in \cls{A}{\Phi(\mathtt{idl}_X)}
    & \Leftrightarrow \gen{x} \in \cls{A}{\Phi(\mathtt{idl}_X)} \\
    & \Leftrightarrow \mathtt{idl}_X(x) \cap A \neq \{0\} \\
    & \Leftrightarrow (\exists r \in X)(rx \in A),
  \end{align*}
  where \gen{x}{} is the subring generated by $x$,
  $A \in \mc{\Phi(\mathtt{idl}_X)}$ if and only if for any $x \in X$:
  \begin{equation*}
    (\exists r \in X)(rx \in A) \Rightarrow x \in A,
  \end{equation*}
  $I \in \Idl{X} \Rightarrow \cls{I}{\Phi(\mathtt{idl}_X)} \in \Idl{X}$,
  $\intr{A}{\Phi(\mathtt{idl}_X)}$ is the \emph{ideal core} of $A$ (i.e.,
  the largest ideal contained in a subring), for each ring
  homomorphism \Arr{f}{X}{Y} the preneighbourhood morphism
  \Arr{f}{\opair{X}{\Phi(\mathtt{idl}_X)}}{\opair{Y}{\Phi(\mathtt{idl}_Y)}}{}
  is continuous and every \Arr{f}{\opair{X}{\mu}}{\opair{Y}{\phi}}{} where
  $\mu \supseteq \mathtt{idl}_{X}$ and $\phi \subseteq \mathtt{idl}_{Y}$ is continuous.
\end{subsubsection}

\begin{subsubsection}{}
  \label{sssec:extreme-nbd}{}
  Given a context $\CAL{A} = (\Bb{A}, \mathsf{E}, \mathsf{M})$, for
  any $c \in \mathtt{EGM}(X)$:
  \begin{align}
    \label{eq:far-for-egm}
    \mathtt{Far}_{\Phi(c)}p
    & = 
      \bigl\{
      x \in \Sub{X}{\mathsf{M}}: c(x) \wedge p = \sigma_{X}
      \bigr\}, \\[1.2ex]
    \label{eq:closure-for-egm}
    \cls{p}{\Phi(c)}
    & = \bigvee
      \bigl\{
      x \not> p: c(x) \wedge p \neq \sigma_{X}
      \bigr\}, \\[1.2ex]
    \label{eq:closed-for-egm}
    \mc{\Phi(c)}
    & = 
      \Biggl\{
      p \in \Sub{X}{\mathsf{M}}: x \in \Sub{X}{\mathsf{M}} \Rightarrow 
      \biggl(
      (x > p) \text{ or } 
      \bigl(
      c(x) \wedge p \neq \sigma_{X} \Rightarrow x \leq p
      \bigr)
      \biggr)
      \Biggr\},     
  \end{align}
  while $\mo{\Phi(c)} = \Fix{c}$,
  $\intr{p}{\Phi(c)} = \bigvee \bigl\{ x \in \Fix{c}: x \leq p \bigr\}$.

  In particular, \mc{\nabla_X}{} is a chain and each element of \mc{\nabla_X}{}
  \emph{cuts} the lattice \Sub{X}{\mathsf{M}}{} in two parts, one
  above and the other below. In presence of good conditions on
  \Sub{X}{\mathsf{M}}{} (e.g., if each element is complemented)
  \cls{p}{\nabla_X}{} is closed, i.e., \cls{}{\nabla_X}{} is
  idempotent. Further, since \mc{\nabla_X}{} is closed under joins,
  \cls{}{\nabla_X}{} is additive, irrespective of presence of
  distributivity. On the other hand, every atom is in \mc{\uparrow_X}{,} the
  $p \in \mc{\uparrow_X}$ is nearly like an atom: for every
  $x\in\Sub{X}{\mathsf{M}}$ either $x \wedge p = \sigma_{X}$, $p < x$ or
  $x \leq p$, the last condition being the exception to an atom.
\end{subsubsection}{}



\begin{subsubsection}{}
  The functions \arrp{C}{D}{\Sub{X}{\mathsf{M}}}{\Sub{X}{\mathsf{M}}}
  provide competitors for a closure operation induced from a
  preneighbourhood system $\mu$ on $X$:
  \begin{align}
    \label{df:closure-alt}{}
    C(p) & = \bigvee
           \bigl\{
           x \in \Sub{X}{\mathsf{M}}: x \neq \id{X}\text{ and }x \notin \mathtt{Far}_{\mu}p
           \bigr\}, \\ \intertext{and}
    D(p) & = \bigvee
           \bigl\{
           x \in \Sub{X}{\mathsf{M}}: \mu(x) \varsubsetneqq \mu(p)\text{ and }x \notin \mathtt{Far}_{\mu}p
           \bigr\}
  \end{align}

  Evidently, $C, D$ are closure operations and
  $D \leq \cls{p}{\mu} \leq C$ and hence $\Fix{C} \leq \mc{\cls{}{\mu}} \leq \Fix{D}$.

  However, the function $C$ is trivial, i.e., either
  $C(p) = \sigma_{X}$ or
  $C(p) = \bigvee \bigl\{ x \in \Sub{X}{\mathsf{M}}: x \neq \id{X}
  \bigr\}$\footnote{I am thankful to the anonymous referee who pointed
    the triviality of this function.}. To see this, let
  $C(p) \neq \sigma_{X}$ and choose a $x \neq \id{X}$ with
  $x \notin \mathtt{Far}_{\mu}p$. Since $x \neq 1$ there exists a
  $y \neq\id{X}$ such that $x \vee y \neq \id{X}$ and for such a $y$,
  $x \vee y \notin \mathtt{Far}_{\mu}p$. Consequently:
  \begin{equation*}
    C(p) = \bigvee
    \bigl\{
    x \in \Sub{X}{\mathsf{M}}: x \neq \id{X}
    \bigr\} \leq \bigvee
    \bigl\{
    x \vee y: y\neq \id{X}
    \bigr\} \leq C(p).
  \end{equation*}

  On the other hand, the function $D$ uses a partial order relation on
  \Sub{X}{\mathsf{M}}{} larger than the usual partial order. Since the
  order structure on \Sub{X}{\mathsf{M}}{} make centre stage of this
  paper, the closure operation $D$ is left for a future detailed
  investigation.

\end{subsubsection}


\section{Closed morphisms}
\label{sec:closed-morphisms}

Having defined a closure operation induced from a preneighbourhood system, this section describe morphisms which preserve the closure operation.

\subsection{}
\label{ssec:closed-morphisms}
Let $\CAL{A} = (\Bb{A}, \mathsf{E}, \mathsf{M})$ be a context.

\begin{Df}
  \label{df:closed-morphisms}{}
  Given the internal preneighbourhood spaces \opair{X}{\mu} and
  \opair{Y}{\phi}, a morphism \Arr{f}{X}{Y}{} is \emph{$\mu$-$\phi$ closed},
  or simply \emph{closed} when the preneighbourhood systems are
  evident, if:
  \begin{equation}
    \label{eq:closed-morphisms}
    p \in \mc{\mu} \Rightarrow \img{f}{p} \in \mc{\phi}.
  \end{equation}
  The (possibly large) set of closed morphisms is
  \closedmorphism{\Bb{A}}.
\end{Df}

\begin{Thm}
  \label{thm:closed-morphisms}{}
  Given the internal preneighbourhood spaces \opair{X}{\mu}{,}
  \opair{Y}{\phi}{,} \opair{Z}{\psi}{} and the morphisms
  \Arr{f}{X}{\Arr{g}{Y}{Z}}{} the following statements are true.
  \begin{enumerate}[label=(\alph*{})]
  \item \label{item:closed-morphisms-alt}{} The morphism $f$ is a
    closed morphism if and only if for every
    $p \in \Sub{X}{\mathsf{M}}$:
    \begin{equation}
      \label{eq:closed-morphisms-alt}
      \clsm{\img{f}{p}}{\phi} \leq \img{f}{\clsm{p}{\mu}}.
    \end{equation}

  \item \label{item:closed-morphisms-alt-for-cont}{} If $f$ is
    continuous then $f$ is closed if and only if for every
    $p \in \Sub{X}{\mathsf{M}}$,
    $\img{f}{\clsm{p}{\mu}} = \clsm{\img{f}{p}}{\phi}$. In particular,
    $m \in \Sub{X}{\mathsf{M}}$ is a closed map if and only if
    $m \in \mc{\mu}$.

  \item \label{item:iso-closed}{} The set \closedmorphism{\Bb{A}}{}
    contain all isomorphisms.

  \item \label{item:closed-comp-closed}{} The set
    \closedmorphism{\Bb{A}}{} is closed under compositions.

  \item \label{item:closed-morphisms-rt-cancel}{} If \comp{g}{f}{} is
    a closed morphism and $f$ is formally surjective and continuous
    then $g$ is a closed morphism.

  \item \label{item:closed-cont-morphisms-pb-semistable}{} If $f$ is a
    closed continuous morphism then for each $m \in \mc{\phi}$ the
    corestriction $f_{m}$ of $f$ along $m$ is closed and continuous.
  \end{enumerate}
\end{Thm}

\begin{proof}
  The statement in \ref{item:iso-closed}{} is immediate from
  definition, while \ref{item:closed-comp-closed}{} is immediate from
  $\img{\comp{g}{f}}{} = \comp{\img{g}{}}{\img{f}{}}$. If $f$ is
  closed, then for any $p \in \Sub{X}{\mathsf{M}}$,
  $\clsm{p}{\mu} \in \mc{\mu}$ implies
  $\img{f}{\clsm{p}{\mu}} \in \mc{\phi}$ and hence
  $\clsm{\img{f}{p}}{\phi}\leq\img{f}{\clsm{p}{\mu}}$
  ($\because, \img{f}{p} \leq \img{f}{\clsm{p}{\mu}}$); on the other
  hand, if \eqref{eq:closed-morphisms-alt}{} is true then for
  $p \in \mc{\mu}$, $\clsm{\img{f}{p}}{\phi}\leq\img{f}{p}$ proving
  $\img{f}{p} \in \mc{\phi}$. This proves
  \ref{item:closed-morphisms-alt}{.} The first part of statement in
  \ref{item:closed-morphisms-alt-for-cont}{} is immediate from
  \ref{item:closed-morphisms-alt}{} and continuity (Proposition
  \ref{prop:cont-wrt-cl-idem-hull}{\ref{item:cont-wrt-idem-hull}});
  the second part is immediate from definition and Remark
  \ref{rem:closed-embed}{.} If \comp{g}{f}{} is a closed morphism, $f$
  is formally surjective and continuous then for any
  $y \in \Sub{Y}{\mathsf{M}}$:
  \begin{equation*}
    \begin{array}{rlr}
      \clsm{\img{g}{y}}{\psi}
      & = \clsm{\img{g}{\img{f}{\finv{f}{y}}}}{\psi}
      & \text{ (since $f$ is formally surjective)} \\[1.2ex]
      & = \clsm{\img{\comp{g}{f}}{\finv{f}{y}}}{\psi} \\[1.2ex]
      & \leq \img{\comp{g}{f}}{\clsm{\finv{f}{y}}{\mu}}
      & \text{ (since \comp{g}{f}{} is closed)} \\[1.2ex]
      & = \img{g}{\img{f}{\clsm{\finv{f}{y}}{\mu}}} \\[1.2ex]
      & \leq \img{g}{\clsm{\img{f}{\finv{f}{y}}}{\phi}}
      & \text{ (since $f$ is continuous)} \\[1.2ex]
      & = \img{g}{\clsm{y}{\phi}}
      & \text{ (since $f$ is formally surjective)},
    \end{array}    
  \end{equation*}
  proving \ref{item:closed-morphisms-rt-cancel}{.} Finally, given
  $\xymatrix{ {P} \ar@{>->}[]!<2.4ex,0ex>;[r]^-p & {\finv{f}{M}}
    \ar@{>->}[]!<0ex,-2.4ex>;[d]_-{\finv{f}{m}} \ar[r]^-{f_m} & {M}
    \ar@{>->}[]!<0ex,-2.4ex>;[d]^-m \\ & {X} \ar[r]_-f & {Y} }$, where
  $f$ is closed continuous, the square is the pullback of
  $m \in \mc{\phi}$ along $f$, $\finv{f}{m} \in \mc{\mu}$ (Proposition
  \ref{prop:cont-wrt-cl-idem-hull}\ref{item:preimage-pres-closed}{}),
  yielding:
  \begin{equation*}
    \begin{array}{rlr}
      \comp{m}{\clsm{\img{f_m}{p}}{\rest{\phi}{m}}}
      & = \clsm{(\comp{m}{(\img{f_m}{p})})}{\phi}
      & \text{ (using \ref{item:closed-morphisms-alt-for-cont}{,} second part)} \\[1.2ex]
      & = \clsm{(\img{f}{(\comp{(\finv{f}{m})}{p})})}{\phi} \\[1.2ex]
      & = \img{f}{\clsm{(\comp{(\finv{f}{m})}{p})}{\mu}}
      & \text{ (since $f$ is closed and continuous, \ref{item:closed-morphisms-alt-for-cont}{})} \\[1.2ex]
      & = \img{f}{\comp{(\finv{f}{m})}{\clsm{p}{\rest{\mu}{\finv{f}{m}}}}}
      & \text{ (since $\finv{f}{m} \in \mc{\mu}$, \ref{item:closed-morphisms-alt-for-cont}{,} second part)} \\[1.2ex]
      & = \comp{m}{\img{f_m}{\clsm{p}{\rest{\mu}{\finv{f}{m}}}}} \\[1.2ex]
      \Rightarrow \clsm{\img{f_m}{p}}{\rest{\phi}{m}}
      & = \img{f_m}{\clsm{p}{\rest{\mu}{\finv{f}{m}}}},
    \end{array}
  \end{equation*}
  completing the proof of
  \ref{item:closed-cont-morphisms-pb-semistable}{.}
\end{proof}

\begin{rem}
  In view of the second part of Theorem
  \ref{thm:closed-morphisms}{\ref{item:closed-morphisms-alt-for-cont}},
  for any internal preneighbourhood space \opair{X}{\mu}{,}
  $m \in \mc{\mu}$ is called a \emph{closed embedding} and
  \closedembed{\Bb{A}}{} is the (possibly large) set of closed
  embeddings.
\end{rem}


\subsection{}
\label{sssec:fr+r0+E-closed}

In this section some examples of closed morphisms are
provided.

Given the preneighbourhood spaces \opair{X}{\mu} and \opair{Y}{\phi}{} let:
\begin{equation*}
  \mathtt{Cl}\opair{\mu}{\phi} = 
  \bigl\{
  \Arr{f}{X}{Y}: \mu \supseteq \invfil{f}{\phi\img{f}{}}\text{ and }(\forall p \in \Sub{X}{\mathsf{M}})(\clsm{\img{f}{p}}{\phi} \leq \img{f}{\clsm{p}{\mu}})
  \bigr\}.
\end{equation*}
Evidently,
\Arr{\mathtt{Cl}}{\opp{\pnhd{X}}\times\pnhd{Y}}{2^{\Homo{\Bb{A}}{X}{Y}}}{}
is an order preserving map.

\begin{Thm}
  \label{thm:reflzero-fs-imply-closed}{}
  If \Arr{f}{X}{Y}{} is a reflecting zero Frobenius
  \textsf{E}-morphism and $\phi \in \pnhd{Y}$ then
  $f \in \mathtt{Cl}(\invfil{f}{\phi\img{f}{}},\phi)$.
\end{Thm}

\begin{proof}
  Since a Frobenius \textsf{E}-morphism is formally surjective
  (Proposition \ref{prop:fr-fs-connection}{}) for each
  $y \in \Sub{Y}{\mathsf{M}}$, $y = \img{f}{\finv{f}{y}}$. Hence using Lemma \ref{lem:far-of-dom-cod}{,} for any
  $p \in \mc{\invfil{f}{\phi\img{f}{}}}$,
  $y \notin \mathtt{Far}_{\phi}(\img{f}{p}) \Leftrightarrow \finv{f}{y} \notin
  \mathtt{Far}_{\invfil{f}{\phi\img{f}{}}}p \Rightarrow p \leq \finv{f}{y} \text{ or
  }p \geq \finv{f}{y}$; further
  $\img{f}{p} \not\leq y \Leftrightarrow p \not\leq \finv{f}{y}$ (from
  \adjt{\img{f}{}}{\finv{f}{}}{}). Hence if
  $y \notin \mathtt{Far}_{\phi}(\img{f}{p})$ and
  $\img{f}{p} \not\leq y$, then
  $\finv{f}{y} \leq p \Rightarrow y = \img{f}{\finv{f}{y}} \leq \img{f}{p}$, proving
  $\img{f}{p} \in \mc{\phi}$.
 
\end{proof}

Regarding examples in some other specific contexts:
\begin{enumerate}[label=(\roman*{})]
\item In the context $(\Loc, \mathsf{Epi}, \mathsf{RegMon})$, if $X$
  and $Y$ are equipped with their $T$- neighbourhood system then a
  localic map \Arr{f}{X}{Y}{} is a closed morphism if and only if $f$
  is a closed morphism in the usual localic sense.

\item In the context $(\Grp,\mathsf{RegEpi},\mathsf{Mono})$ (respectively, $(\CRng, \mathsf{RegEpi}, \mathsf{Mono})$) if \Arr{f}{X}{Y} with $f \in \mathsf{RegEpi}$ then $f$ is $\Phi(\mathtt{ncl}_X)$-$\Phi(\mathtt{ncl}_X)$ closed (respectively, $\Phi(\mathtt{idl}_X)$-$\Phi(\mathtt{idl}_X)$) closed.
\end{enumerate}


\section{Dense morphisms}
\label{sec:dense-morphisms}

In this section the notion of \emph{dense morphisms} shall be introduced, the \emph{dense}-(closed embedding) factorisation system exhibited.

\subsection{}
\label{ssec:dense-mor}
Let $\CAL{A} = (\Bb{A}, \mathsf{E}, \mathsf{M})$ be a context.

\begin{Df}
  \label{df:dense-mor}
  Given the internal preneighbourhood spaces \opair{X}{\mu}{,}
  \opair{Y}{\phi}{,} a morphism \Arr{f}{X}{Y}{} is \emph{$\mu$-$\phi$ dense}
  morphism (or in short \emph{dense} morphism, if $\mu$ and $\phi$ are
  evident) if $f = \comp{m}{h}$ for some $m \in \mc{\phi}$ implies $m$ is
  an isomorphism. The (possibly large) set of all dense morphisms is
  denoted by \densemorphism{\Bb{A}}{.}
\end{Df}

\begin{rem}\label{rem:dense-is-closure-dense}{}
  Consider the commutative diagram:
  \begin{equation*}
    \xymatrixcolsep{4.8em}
    \xymatrixrowsep{3.6em}
    \xymatrix{
      {X} \ar[r]^-f \ar@{->>}[d]_-{f^{\mathsf{E}}} & {Y} \\
      {\Img{f}} \ar@{>->}[]!<3.6ex,0ex>;[r]_-{j_f}
      \ar@/_{4.8ex}/@{>->}[]!<2.4ex,-2.4ex>;[rr]_-{u_f}
      \ar@{>->}[]!<2.4ex,2.4ex>;[ur]|-{f^{\mathsf{M}}} &
      {\overline{\Img{f}}} \ar@{>->}[]!<2.4ex,0ex>;[r]_-{k_f}
      \ar@{>->}[]!<0ex,3.06ex>;[u]|(0.36){\cls{f^{\mathsf{M}}}{\phi}} &
      {\widehat{\Img{f}}} 
      \ar@{>->}[]!<-2.4ex,2.4ex>;[ul]_-{\clsm{f^{\mathsf{M}}}{\phi}}
    }
  \end{equation*}
  where $j_{f}, k_{f}, u_{f} = \comp{k_f}{j_f}$ are the comparisons
  between the respective admissible subobjects.  Evidently,
  $f = \comp{(\clsm{f^{\mathsf{M}}}{\phi})}{\comp{u_f}{f^{\mathsf{E}}}}$
  is dense if and only if $\clsm{f^{\mathsf{M}}}{\phi} = \id{Y}$. For a
  general $f$, since
  $ \left( \comp{u_f}{f^{\mathsf{E}}} \right)^{\mathsf{E}} =
  f^{\mathsf{E}}$,
  $ \left( \comp{u_f}{f^{\mathsf{E}}} \right)^{\mathsf{M}} = u_{f}$,
  an use of \eqref{eq:closure-idemp-hereditary} shows on taking
  $m = \clsm{f^{\mathsf{M}}}{\phi}$:
  \begin{equation*}
    \clsm{(\comp{u_f}{f^{\mathsf{E}}})^{\mathsf{M}}}{\rest{\phi}{m}} =
    \clsm{u_f}{\rest{\phi}{m}} = \finv{m}{\clsm{(\comp{m}{u_f})}{\phi}} =
    \finv{m}{m} = \id{\widehat{\Img{f}}},
  \end{equation*}
  i.e., \comp{u_f}{f^{\mathsf{E}}}{} is a dense morphism. In
  particular,
  $f =
  \comp{(\clsm{f^{\mathsf{M}}}{\phi})}{(\comp{u_f}{f^{\mathsf{E}}})}$
  shows every morphism factor as a dense morphism followed by a closed
  embedding.
\end{rem}

\begin{Thm}
  \label{thm:dense-mor-prop}
  Given the internal preneighbourhood spaces \opair{X}{\mu}{,}
  \opair{Y}{\phi}{,} \opair{Z}{\psi}{} and the morphisms
  \Arr{f}{X}{\Arr{g}{Y}{Z}}{,} the following statements hold.
  \begin{enumerate}[label=(\alph*{})]

  \item \label{item:dense-iff-image-is-dense}{} The morphism $f$ is a
    dense morphism if and only if $\clsm{f^{\mathsf{M}}}{\phi} = \id{Y}$.

  \item \label{item:dense-has-E-and-iso-meet-with-clemb}{}
    $\mathsf{E} \subseteq \densemorphism{\Bb{A}}$ and
    $\densemorphism{\Bb{A}} \cap \closedembed{\Bb{A}} = \Iso{\Bb{A}}$.

  \item \label{item:dense-cont-and-dense-comp-is-dense}{} If $g$ is
    dense continuous and $f$ is dense then \comp{g}{f}{} is dense.
      
  \item \label{item:dense-mor-rt-cancel} If
    $\comp{g}{f} \in \densemorphism{\Bb{A}}$ then
    $g \in \densemorphism{\Bb{A}}$.

  \item \label{item:dense-closed-emb-almost-pfs}
    $\dn{\densemorphism{\Bb{A}}} \subseteq \closedembed{\Bb{A}}$ and
    $\up{\closedembed{\Bb{A}}} \subseteq \densemorphism{\Bb{A}}$.

  \item \label{item:dense-clemb-pfs} If every preneighbourhood
    morphism is continuous then
    $\densemorphism{\Bb{A}} \subseteq \up{\closedembed{\Bb{A}}}$ and
    \opair{\densemorphism{\Bb{A}}}{\closedembed{\Bb{A}}}{} is a
    factorisation system on \Bb{A}{.}

  \item \label{item:dense-po-stab} If all preneighbourhood morphisms
    are continuous then dense morphisms are pushout stable.
      
  \item \label{item:dense-colimit-closed}{} If all preneighbourhood
    morphisms are continuous and
    $\xymatrix{ {\Bb{X}} \ar@<2ex>[r]^-F="u" \ar@<-2ex>[r]_-G="d" &
      {\Bb{A}} \ar@2{->} <0pt,-2pt>+"u"; <0pt,2pt>+"d"^-{\alpha} }$
    $\alpha_X \in \densemorphism{\Bb{A}}$ ($X \in \Bb{X}_0$), \colimit{F}{}{}
    and \colimit{G}{}{} exists then
    $\colimit{\alpha}{} \in \densemorphism{\Bb{A}}$.
  \end{enumerate}
\end{Thm}

\begin{proof}
  Since
  $f = \comp{(\clsm{f^{\mathsf{M}}}{\phi})}{\comp{u_f}{f^{\mathsf{E}}}}$
  (Remark \ref{rem:dense-is-closure-dense}), the proof of
  \ref{item:dense-iff-image-is-dense} follows. If $f$ is dense and
  $f \in \mc{\phi}$ then
  $f = f^{\mathsf{M}} = \clsm{f^{\mathsf{M}}}{\phi}$, proving
  $\closedembed{\Bb{A}} \cap \densemorphism{\Bb{A}} \subseteq \Iso{\Bb{A}}$; the
  other inclusion is obvious. Also, if $f \in \mathsf{E}$ and
  $f = \comp{m}{h}$ with $m \in \mc{\phi}$ then
  $m \in \mathsf{E} \cap \mathsf{M} = \Iso{\Bb{A}}$. This proves
  \ref{item:dense-has-E-and-iso-meet-with-clemb}{.} If $g$ is dense
  continuous, $f$ is dense and $\comp{g}{f} = \comp{m}{h}$ for some
  $m \in \mc{\psi}$, then from the pullback of $m$ along $g$:
  \begin{equation*}
    \xymatrixcolsep{4.8em}
    \xymatrixrowsep{2.4em}
    \xymatrix{
      {X} \ar@/^{6ex}/[rrd]^-h \ar@/_{6ex}/[ddr]_-f \ar@{.>}[dr]|-{!\,u} \\
      & {\cdot} \ar[r]^-{g_m} \ar@{>->}[]!<0ex,-2.4ex>;[d]_-{\finv{g}{m}}
      & {\cdot} \ar@{>->}[]!<0ex,-2.4ex>;[d]^-m \\
      & {Y} \ar[r]_-g & {Z}
    }
  \end{equation*}
  $\finv{g}{m} \in \mc{\phi}$ (Proposition
  \ref{prop:cont-wrt-cl-idem-hull}\ref{item:preimage-pres-closed}{}),
  there exists a unique morphism $u$ such that
  $f = \comp{(\finv{g}{m})}{u}$; hence from density of $f$,
  \finv{g}{m}{} is an isomorphism. Hence
  $g = \comp{m}{\comp{g_m}{\inv{(\finv{g}{m})}}}$ forces $m$ an
  isomorphism from density of $g$. This proves
  \ref{item:dense-cont-and-dense-comp-is-dense}{.} If
  $\comp{g}{f} \in \densemorphism{\Bb{A}}$, $g = \comp{m}{h}$ for some
  $m \in \mc{\psi}$, then $\comp{g}{f} = \comp{m}{\comp{h}{f}}$ forces from
  density of \comp{g}{f}{,} $m$ is an isomorphism; hence
  $g \in \densemorphism{\Bb{A}}$, proving
  \ref{item:dense-mor-rt-cancel}{.} If
  $f \in \dn{\densemorphism{\Bb{A}}}$, then
  \down{(\comp{u_f}{f^{\mathsf{E}}})}{f} forces $f^{\mathsf{E}}$ an
  isomorphism proving $f \in \mathsf{M}$. Consequently,
  $f = \comp{(\clsm{f}{\phi})}{u_f}$ and \down{u_f}{f} forces
  $\clsm{f}{\phi} \leq f$, i.e., $f \in \closedembed{\Bb{A}}$, proving
  $\dn{\densemorphism{\Bb{A}}} \subseteq \closedembed{\Bb{A}}$; the inequality
  $\up{\closedembed{\Bb{A}}} \subseteq \densemorphism{\Bb{A}}$ is trivial,
  proving \ref{item:dense-closed-emb-almost-pfs}{.} If every
  preneighbourhood morphism is continuous, $f$ is dense and
  $\comp{v}{f} = \comp{m}{u}$ for some $m \in \closedembed{\Bb{A}}$,
  then \finv{v}{m}{} is a closed embedding (Proposition
  \ref{prop:cont-wrt-cl-idem-hull}{\ref{item:preimage-pres-closed}{}})
  and there exists a unique morphism $w$ such that
  $f = \comp{(\finv{v}{m})}{w}$, implying \finv{v}{m}{} an
  isomorphism; hence \comp{v_m}{\inv{(\finv{v}{m})}}{} is the unique
  diagonal fill-in, i.e., \down{f}{m}{,} proving
  $\densemorphism{\Bb{A}} \subseteq \up{\closedembed{\Bb{A}}}$. Hence if every
  preneighbourhood morphism is continuous,
  $\densemorphism{\Bb{A}} = \up{\closedembed{\Bb{A}}}$,
  $\dn{\densemorphism{\Bb{A}}} = \dn{\up{\closedembed{\Bb{A}}}} \subseteq
  \closedembed{\Bb{A}} \subseteq \up{\dn{\closedembed{\Bb{A}}}}$ proving
  $\dn{\densemorphism{\Bb{A}}} = \closedembed{\Bb{A}}$. This completes
  the proof of \ref{item:dense-clemb-pfs}{.} The rest of the
  properties in \ref{item:dense-po-stab},
  \ref{item:dense-colimit-closed}{} follow from
  \ref{item:dense-clemb-pfs}{} and the properties of a
  prefactorisation system \cite[see][Proposition 2.2]{Janel1997b}.
\end{proof}

\begin{rem}
  \label{rem:dense-closed-emb}
  If every preneighbourhood morphism is continuous then from Theorem
  \ref{thm:dense-mor-prop}\ref{item:dense-clemb-pfs},
  \opair{\densemorphism{\Bb{A}}\cap\pNHD{\Bb{A}}_1}
  {\closedembed{\Bb{A}}\cap\pNHD{\Bb{A}}_1}{} is a factorisation
  structure for \pNHD{\Bb{A}}{.} The dense-(closed embedding)
  factorisation of preneighbourhood morphisms is well known for
  topological spaces, i.e., for neighbourhood systems of the context
  $(\Set, \mathsf{Surjections}, \mathsf{Injections})$, as well as for
  locales with $T$-neighbourhood systems
  \cite[see][\S{XV.2.2}{}]{PicadoPultr2012}. Theorem
  \ref{thm:dense-mor-prop}\ref{item:dense-clemb-pfs}{} generalise
  these results to larger subcategories of preneighbourhood spaces.
\end{rem}

\begin{rem}
  \label{rem:dense-closed-emb-proper}{}
  The factorisation system \opair{\densemorphism{\Bb{A}} \cap
    \pNHD{\Bb{A}}_1{}}{\closedembed{\Bb{A}} \cap\pNHD{\Bb{A}}_1{}}{} is
  proper if and only if for each internal neighbourhood space
  \opair{X}{\mu},
  $\xymatrix{ {\opair{X}{\mu}} \ar@{>->}[]!<24pt,0pt>;[r]^-{d_{X}} &
    {\opair{X \times X}{\mu \times \mu}} }$ is a closed embedding, i.e.,
  \opair{X}{\mu}{} is an internal \emph{Hausdorff space} (see Definition
  \ref{df:hausdorff}{} and Theorem \ref{thm:ihs-alt}{}) --- note the
  assumption of \emph{continuity of each preneighbourhood morphism} is
  already embedded by Theorem
  \ref{thm:dense-mor-prop}\ref{item:dense-clemb-pfs}{,} relaxing the
  \emph{proper} condition in Definition \ref{df:hausdorff} (see also
  Theorem
  \ref{thm:proper-mor-prop}\ref{item:proper-mor-have-closed-emb}{}){.}
\end{rem}


\subsection{}
\label{ssec:densemaps-ex}{}
This section exhibit examples of dense subobjects.

\begin{ex}
  \label{ex:dense-morphisms-set}{}
  In the context $(\Set, \mathsf{Surjections}, \mathsf{Injections})$
  the dense morphisms between neighbourhood spaces are precisely the
  usual continuous maps with dense image.

  In the context $(\Top, \mathsf{Epi}, \mathsf{ExtMon})$, the dense
  morphisms between neighbourhood spaces are precisely the
  bicontinuous functions between the bitopological spaces, which have
  dense image with respect to the second topologies.
\end{ex}

\begin{ex}
  \label{ex:dense-morphisms-loc}
  In the context $(\Loc, \mathsf{Epi}, \mathsf{RegMon})$, a localic
  map \Arr{f}{X}{Y}{} is a dense morphism with respect to the
  $T$-neighbourhood systems on $X$ and $Y$, if and only if,
  $f(0) = 0$, i.e., is a dense localic map in the usual sense
  \cite[see][\S8.2]{PicadoPultr2012}.
\end{ex}

\begin{ex}
  \label{ex:dense-morphisms-grp}{}
  In the context $(\Grp, \mathsf{RegEpi}, \mathsf{Mono})$ a group
  homomorphism \Arr{f}{X}{Y} is
  $\Phi(\mathtt{ncl}_X)$-$\Phi(\mathtt{ncl}_Y)$ dense if and only if the
  image $f(X)$ non-trivially meets every non-trivial normal subgroup
  of $Y$, i.e., $f(X)$ is \emph{essential with respect to normal
    subgroups}.

  \begin{Df}
    \label{df:essential}{}
    An admissible subobject $m \in \Sub{X}{\mathsf{M}}$ is
    \emph{essential with respect to \eulf{a}} or
    \emph{\eulf{a}{-}essential}
    ($\eulf{a} \subseteq \Sub{X}{\mathsf{M}}$) if
    $a \in \eulf{a} \Rightarrow a \wedge m \neq \sigma_{X}$; \Sub{X}{\mathsf{M}}{-}essential is
    shortened to \emph{essential}.\footnote{The idea of essential
      abelian groups can be obtained from \cite[][page
      19]{Griffith1970}.}
  \end{Df}
\end{ex}

\begin{ex}
  \label{ex:dense-morphisms-rng}{}
  In the context $(\CRng, \mathsf{RegEpi}, \mathsf{Mono})$ a ring
  homomorphism \Arr{f}{X}{Y} is
  $\Phi(\mathtt{idl}_X)$-$\Phi(\mathtt{idl}_Y)$ dense if and only if the
  image $f(X)$ non-trivially meets every non-trivial ideal $Y$, i.e.,
  the image $f(X)$ is essential with respect to ideals.
\end{ex}


\section{Stably closed morphisms}
\label{sec:stably-clos-morph}

This section discuss pullback stable closed morphisms and \emph{compact preneighbourhood} spaces as a special case.

\subsection{}
\label{ssec:stably-clos-morph}
Let $\CAL{A} = (\Bb{A}, \mathsf{E}, \mathsf{M})$ be a context.

\begin{Df}
  \label{df:proper-mor}

  A preneighbourhood morphism \Arr{f}{(X, \mu)}{(Y, \phi)} is said to be
  {\em proper} if it is stably in \closedmorphism{\Bb{A}}, i.e., for
  every preneighbourhood morphism \Arr{h}{(T, \tau)}{(Y, \phi)}, the
  pullback $f_{h}$ of $f$ along $h$ is a closed morphism. The symbol
  \propermorphism{\Bb{A}} denotes the (possibly large) set of proper
  morphisms in \Bb{A}.
\end{Df}

\begin{Thm}
  \label{thm:proper-mor-prop}
  Let \Arr{f}{\opair{X}{\mu}}{\Arr{g}{\opair{Y}{\phi}}{\opair{Z}{\psi}}},
  \Arr{h}{\opair{T}{\tau}}{\opair{Y}{\phi}}{} be preneighbourhood morphisms.
  \begin{enumerate}[label=(\alph*{})]
  \item \label{item:proper-mor-alt} The preneighbourhood morphism $f$
    is a proper morphism if and only if for any internal
    preneighbourhood space \opair{T}{\tau} every corestriction of
    \Arr{f\times\id{T}}{X\times{T}}{Y\times{T}}{} along a section of the second
    product projection $p_{2}$ is a closed morphism.
  
  \item \label{item:proper-mor-pb-stable-comp-closed}{} The set
    \propermorphism{\Bb{A}} is a pullback stable set, is closed under
    compositions.

  \item \label{item:proper-mor-have-closed-emb}{} If all
    preneighbourhood morphisms are continuous then
    $\closedembed{\Bb{A}} \subseteq \propermorphism{\Bb{A}}$.

  \item \label{item:proper-mor-rt-can-by-prop-mor-stably-in-E} If
    $\comp{g}{f} \in \propermorphism{\Bb{A}}$ and $f$ is stably in
    \textsf{E} and is stably continuous, i.e., for any morphism $h$
    the pullback $f_{h}$ of $f$ along $h$ is in \textsf{E} and is
    continuous, then $g \in \propermorphism{\Bb{A}}$.

  \item \label{item:proper-mor-lt-can-by-mono} If
    $\comp{g}{f} \in \propermorphism{\Bb{A}}$ and
    $g \in \Mono{\Bb{A}}$ then $f \in \propermorphism{\Bb{A}}$.

  \item \label{item:prod-prop-is-prop}{} If $f$ is a proper morphism
    then $f \times f$ is a proper morphism.

  \end{enumerate}
\end{Thm}

\begin{proof}
  Firstly, sections of the second product projection
  \Arr{p_2}{Y\times{T}}{T}{} are precisely determined by
  \Arr{h}{\opair{T}{\tau}}{\opair{Y}{\phi}}{,} namely \opair{h}{\id{T}}{.}
  Towards the proof of \ref{item:proper-mor-alt}, consider the
  commutative diagram:
  \begin{equation*}
    \xymatrixcolsep{1.8cm}
    \xymatrix{
      {P} \ar[rr]^-{f_h} \ar[dd]_-{h_f} \ar[dr]^-{(h_f, f_h)} & &
      {T} \ar[dd]^(0.66)h \ar[dr]^-{(h, \id{T})} \\
      & {X \times T} \ar[rr]|(0.492){\hole}|(0.66){f \times \id{T}}
      \ar[dl]^-{p_1} & &
      {Y \times T}
      \ar[dl]^-{p_1} \\
      {X} \ar[rr]_-f & & {Y}
    }
  \end{equation*}
  in which $p_1$'s are product projections, the horizontal square is
  the pullback of $p_1$ along $f$. From properties of pullbacks,
  $f_{h}$ is the pullback of $f$ along $h$ if and only if $f_{h}$ is
  the corestriction of $f \times \id{T}$ along the section
  \opair{h}{\id{T}} of $p_{2}${.} This completes the proof of
  \ref{item:proper-mor-alt}{.}  Since for any $m \in \mc{\phi}$,
  $\finv{f}{m} \in \mc{\mu}$ when $f$ is continuous (Proposition
  \ref{prop:cont-wrt-cl-idem-hull}), closed embeddings are proper when
  every morphism is continuous, proving
  \ref{item:proper-mor-have-closed-emb}{.} For the first part of
  \ref{item:proper-mor-pb-stable-comp-closed}{,}
  \propermorphism{\Bb{A}} is the largest pullback stable subset of
  \closedmorphism{\Bb{A}}. If $f$ and $g$ are proper preneighbourhood
  morphisms then from the diagram
  \begin{equation*}
    \xymatrixcolsep{1.2cm}
    \xymatrix{
      {R} \ar[r]^-{f_{w_g}} \ar[d]_-{w_{\comp{g}{f}}} & {S} \ar[r]^-{g_w}
      \ar[d]|-{w_g} & {W} \ar[d]^-{w} \\
      {X} \ar[r]_-f & {Y} \ar[r]_-g & {Z}
    }
  \end{equation*}
  where \Arr{w}{(W, \omega)}{(Z, \psi)} is a preneighbourhood
  morphism, the right hand square is the pullback of $w$ along $g$ and
  the left hand square is the pullback of $w_g$ along $f$, the outer
  square is the pullback of $w$ along \comp{g}{f}.  If $g$ and $f$ are
  proper morphisms, $g_w$ and $f_{w_g}$ are both closed morphisms and
  hence their composite \comp{g_w}{f_{w_g}} is closed (Theorem
  \ref{thm:closed-morphisms}\ref{item:closed-comp-closed}), proving
  \comp{g}{f} is a proper morphism. This proves the second part of
  \ref{item:proper-mor-pb-stable-comp-closed}{.} On the other hand if
  the composite \comp{g}{f}{} is a proper morphism then
  \comp{g_w}{f_{w_g}}{} is a closed morphism. Further if $f$ is a
  morphism stably continuous and stably in $\mathsf{E}$, $f_{w_{g}}$
  is a continuous morphism stably in $\mathsf{E}$; hence $f_{w_{g}}$
  is a formally surjective continuous morphism. Hence $g_{w}$ is a
  closed morphism (Theorem \ref{thm:closed-morphisms}
  \ref{item:closed-morphisms-rt-cancel}{}) proving $g$ to be a proper
  morphism. This proves
  \ref{item:proper-mor-rt-can-by-prop-mor-stably-in-E}{.}  If
  \comp{g}{f} is a proper morphism and $g$ is a monomorphism then from
  the commutative diagram:
  \begin{equation*}
    \xymatrixcolsep{1.2cm}
    \xymatrix{
      {T'} \ar[r]^-v \ar[d]_-{(u, v)} & {T} \ar[d]|-{(h, \id{T})} \ar@2{-}[r] & {T}
      \ar[d]^-{(\comp{g}{h}, \id{T})} \\
      {X \times T} \ar[r]_-{f \times \id{T}} & {Y \times T} \ar[r]_-{g \times \id{T}}
      & {Z \times T}
    }
  \end{equation*}
  in which the left hand square is the pullback of $f \times \id{T}$
  along \opair{h}{\id{T}}, since $g$ is a monomorphism the right hand
  square is a pullback square, implying the outer square is the
  pullback of \opair{\comp{g}{h}}{\id{T}} along
  $(\comp{g}{f}) \times \id{T}$.  Since \comp{g}{f} is proper, using
  \ref{item:proper-mor-alt} on the outer pullback square, $v$ is a
  closed morphism. Hence using \ref{item:proper-mor-alt}{} again, $f$
  is a proper morphism. This proves
  \ref{item:proper-mor-lt-can-by-mono}.

  Finally, towards the proof of \ref{item:prod-prop-is-prop}, given
  \Arr{\opair{t}{s}}{\opair{T}{\tau}}{\opair{Y\times{Y}}{\phi\times\phi}{}}
  consider the diagrams:
  \begin{equation*}
    \label{eq:product-proper}\tag{$\star\star$}
    \begin{aligned}
      \xymatrix{ {Q} \ar[r]^-{v} \ar[d]_-u & {Q'} \ar[r]^-{s_f}
                                             \ar[d]|-{f_s}
      & {X} \ar[d]^-f \\
      {P} \ar[r]|-{f_t} \ar[d]_-{t_f} & {T} \ar[d]^-t \ar[r]_-s
      & {Y} \\
      {X} \ar[r]_-f & {Y} }
    \end{aligned}\quad\text{and}\quad
    \begin{aligned}
      \xymatrixcolsep{1.56cm} \xymatrixrowsep{1.2cm} \xymatrix{ {Q}
      \ar[r]^-{\comp{f_t}{u}}
      \ar[d]_-{\opair{\comp{t_f}{u}}{\comp{s_f}{v}}}
      & {T} \ar[d]^-{\opair{t}{s}} \\
      {X \times X} \ar[r]_-{f\times{f}} & {Y \times Y} }
    \end{aligned}
  \end{equation*}
  where in the left hand diagram each square is a pullback
  square. Since $f$ is proper then $f_{t}$, $f_{s}$, $u$ and $v$ are
  each closed morphisms. Clearly, the right hand square of
  \eqref{eq:product-proper}{} is a pullback square. Hence the
  composite \comp{f_t}{u}, which is the pullback of $f \times f$ along
  \opair{t}{s}{,} is a closed morphism (Theorem
  \ref{thm:closed-morphisms}\ref{item:closed-comp-closed}{}). Therefore
  $f \times f$ is a proper morphism, proving
  \ref{item:prod-prop-is-prop}.

\end{proof}


\subsection{}
\label{ssec:propermaps-ex}

This section exhibit examples of proper maps.

\begin{ex}
  \label{ex:proper-maps-top-loc}
  In the context $(\Set, \mathsf{Surjections}, \mathsf{Injections})$
  the proper maps of internal neighbourhood spaces are precisely the
  proper maps of topological spaces. In the context
  $(\Top, \mathsf{Epi}, \mathsf{ExtMon})$ the proper maps are
  precisely the proper maps between the second topological spaces.  In
  the context $(\Loc, \mathsf{Epi}, \mathsf{RegMon})$ the proper maps
  between the internal $T$-neighbourhood spaces are precisely the
  usual localic proper maps \cite[see][]{Vermeulen2001,Vermeulen1994,PultrTozzi1992}.

\end{ex}

\begin{ex}
  \label{ex:proper-maps-grp-rng}{}
  In the context $(\Grp,\mathsf{RegEpi},\mathsf{Mono})$ (respectively, $(\CRng, \mathsf{RegEpi}, \mathsf{Mono})$) if \Arr{f}{\opair{X}{\Phi(\mathtt{ncl}_X)}}{\opair{Y}{\Phi(\mathtt{ncl}_Y)}} (respectively, \Arr{f}{\opair{X}{\Phi(\mathtt{idl}_X)}}{\opair{Y}{\Phi(\mathtt{idl}_Y)}})  with $f \in \mathsf{RegEpi}$ then $f$ is proper.
\end{ex}


\subsection{}
\label{ssec:compact-obj}

This section introduce the \emph{compact preneighbourhood spaces}.

\begin{df}
  \label{df:compact-obj}
  An internal preneighbourhood space \opair{X}{\mu}{} is \emph{compact}
  if the unique morphism
  \Arr{\terma{X}}{\opair{X}{\mu}}{\opair{\termo}{\nabla_{\termo}}} is
  proper. The full subcategory of all compact objects is denoted by
  \Comp{\pNHD{\Bb{A}}}{}.
\end{df}

\begin{rem}
  \label{rem:compact}{} Immediately from Theorem
    \ref{thm:proper-mor-prop}{\ref{item:proper-mor-alt}}: an internal
    preneighbourhood space \opair{X}{\mu}{} is compact if and only if
    for every internal preneighbourhood space \opair{Y}{\phi}{,} the
    projection \Arr{p_2}{\opair{X \times Y}{\mu \times \phi}}{\opair{Y}{\phi}{}} is a
    closed morphism, in fact a proper morphism.
  \end{rem}

\begin{Thm}
  \label{thm:compact-prop}
  \begin{enumerate}[label=(\alph*)]
  \item \label{item:dom-prop-cpt-cod-is-prop} If \opair{Y}{\phi}{} is
    compact and \Arr{f}{\opair{X}{\mu}}{\opair{Y}{\phi}}{} is a proper
    morphism then \opair{X}{\mu}{} is compact.

  \item \label{item:stably-E-img-cpt-is-cpt} If \opair{X}{\mu}{} is
    compact and \Arr{f}{\opair{X}{\mu}}{\opair{Y}{\phi}} is a
    preneighbourhood morphism with $f$ stably continuous and stably in \textsf{E} then
    \opair{Y}{\phi}{} is compact.

  \item \label{item:cpt-fin-productive-closed-her}{} The category
    \Comp{\pNHD{\Bb{A}}} is finitely productive; if every preneighbourhood morphism is continuous then \Comp{\pNHD{\Bb{A}}}{} is closed hereditary, i.e.,
    if \opair{X}{\mu}{} is compact and $m \in \mc{\mu}$ then
    \opair{M}{\rest{\mu}{M}}{} is compact.
  \end{enumerate}
\end{Thm}

\begin{proof}
  Since $\terma{X} = \comp{\terma{Y}}{f}$,
  \ref{item:dom-prop-cpt-cod-is-prop} \&
  \ref{item:stably-E-img-cpt-is-cpt}{} follow from Definition and
  Theorem \ref{thm:proper-mor-prop}(\ref{item:proper-mor-pb-stable-comp-closed} \& \ref{item:proper-mor-rt-can-by-prop-mor-stably-in-E}). Since
  every isomorphism is a proper morphism, \opair{\termo}{\nabla_{\termo}}{}
  is compact. Further, binary products of compact objects is proper
  from \ref{item:dom-prop-cpt-cod-is-prop} and Definition \ref{df:compact-obj}{.} Hence \Comp{\pNHD{\Bb{A}}}{} is
  finitely productive. Since every closed embedding is a proper
  morphism if every preneighbourhood morphism is continuous (Theorem
  \ref{thm:proper-mor-prop}\ref{item:proper-mor-have-closed-emb}{}) the closed
  heredity of \Comp{\pNHD{\Bb{A}}}{} follows from \ref{item:dom-prop-cpt-cod-is-prop}{}.
\end{proof}

A detailed treatment of \Comp{\pNHD{\Bb{A}}}{} shall be done in a later
paper.


\section{Separated morphisms}
\label{sec:separated-mor}

Given the context $(\Bb{A}, \mathsf{E}, \mathsf{M})$, let \Arr{f}{\opair{X}{\mu{}}}{\opair{Y}{\phi{}}} be a preneighbourhood
morphism with \arrp{f_1}{f_2}{\Kerp{f}}{X} its kernel
pair. Since $\xymatrixcolsep{6em}\xymatrix{ {\opair{X}{\rest{(\mu\times_{\phi}\mu)}{d_f}}} \ar@{>->}[]!<11.4ex,0ex>;[r]^-{d_f = \opair{\id{X}}{\id{X}}} & {\opair{\Kerp{f}}{\mu\times_{\phi}\mu}} }$, $X$ is an admissible
subobject of \Kerp{f}.  Evidently, any
morphism \Arr{t}{T}{\Kerp{f}{}} is determined by the pair
\arrp{t_1}{t_2}{T}{X} of morphisms such that $\comp{f_{i}}{t} = t_{i}$
($i = 1, 2$) and
$t_{i} = \comp{t_{i}^{\mathsf{E}}}{t_{i}^{\mathsf{M}}}$ is the
\fact{\textsf{E}}{\textsf{M}} of $t_{i}$ ($i = 1,2$).
If $t = \opair{t_{1}}{t_{2}} \in \Sub{\Kerp{f}}{\mathsf{M}}$
  and
  $\xymatrixcolsep{1.2cm} \xymatrix{ {\levset{t_{1}}{t_{2}}{=}\,\,}
    \ar@{>->}[r]^-{m_{t}} & {T} \ar@<0.6ex>[r]^-{t_{1}}
    \ar@<-0.6ex>[r]_-{t_{2}} & {X} }$ is the equaliser of the pair
  \opair{t_1}{t_2}, then:
  \begin{align}
    \label{eq:diag-inv}
    \finv{d_{f}}{t} & = \comp{t_1}{m_t} = \comp{t_{2}}{m_{t}}, \\ 
    \label{eq:diag-meet}
    d_{f} \wedge t & = \opair{\comp{t_{1}}{m_{t}}}{\comp{t_{2}}{m_{t}}}, \\
    \intertext{and}
    \label{eq:pnbd-diag-inv}
    \mu(\finv{d_{f}}{t}) & \supseteq \mu(t_{1}^{\mathsf{M}}) \vee \mu(t_{2}^{\mathsf{M}}).
  \end{align}

  The statements in \eqref{eq:diag-inv}\&\eqref{eq:diag-meet}{} are trivial computations; for  \eqref{eq:pnbd-diag-inv}{,}
  $\comp{t_i}{m_t} =
  \comp{t_i^{\mathsf{M}}}{\comp{t_i^{\mathsf{E}}}{m_t}}$ implies
  $\comp{t_i}{m_t} \leq t_{i}^{\mathsf{M}}$, ($i = 1, 2$) yielding the
  result from \eqref{eq:diag-inv}{.} An use of
  \eqref{eq:pnbd-diag-inv}{} shows
  $\rest{(\mu\times_{\phi}\mu)}{d_f} \leq \mu$. Since for every $v \in \mu(u)$, $v = \finv{d_f}{(\finv{f_1}{v}\wedge\finv{f_2}{v})} \Leftrightarrow d_{f} \wedge \finv{f_1}{v}\wedge\finv{f_2}{v} \leq \comp{d_f}{v} = \opair{v}{v}$, $v \in \rest{(\mu\times_{\phi}\mu)}{d_f}(u)$; hence $\mu = \rest{(\mu\times_{\phi}\mu)}{d_f}$.


\subsection{}
\label{ssec:separated-morphisms}

This section discuss the notion of \emph{separated morphisms}.

\begin{Df}
  \label{df:separated-morphism}
  A preneighbourhood morphism \Arr{f}{(X, \mu)}{(Y, \phi)} is said to be a
  {\em separated morphism} if $d_f$ is a proper morphism.  The symbol
  \sepmorphism{\Bb{A}} denotes the (possibly large) set of all
  separated morphisms of \Bb{A}.
\end{Df}

\begin{rem}
  In case when every preneighbourhood morphism is continuous an use of
  Theorem \ref{thm:proper-mor-prop}{} shows, a preneighbourhood
  morphism $f$ is separated if and only of $d_{f}$ is a closed
  embedding.
\end{rem}

\begin{Ex}\label{ex:sep-map}
  In the context $(\Set, \mathsf{Surjections}, \mathsf{Injections})$
  the separated morphisms between internal neighbourhood spaces are
  precisely those continuous maps in whose fibres distinct points are
  separated by disjoint neighbourhoods. In the context
  $(\Top, \mathsf{Epi}, \mathsf{ExtMon})$, the separated morphisms
  between the internal neighbourhood spaces are precisely the
  separated maps with respect to the second topologies.
\end{Ex}

\begin{Thm}
  \label{thm:sep-map-prop}
  Let \Arr{f}{\opair{X}{\mu}}{\Arr{g}{\opair{Y}{\phi}}{\opair{Z}{\psi}}}{} be
  preneighbourhood morphisms.
  \begin{enumerate}[label=(\alph*{})]
  \item \label{item:sep-mor-contain-mono}{} The set
    \sepmorphism{\Bb{A}}{} contain all monomorphisms.
    
  \item \label{item:sep-mor-pb-stable}{} The set
    \sepmorphism{\Bb{A}}{} is pullback stable.

  \item \label{item:sep-mor-comp-closed}{} If
    $g, f \in \sepmorphism{\Bb{A}}$ then
    $\comp{g}{f} \in \sepmorphism{\Bb{A}}$.

  \item \label{item:sep-mor-rt-can} If
    $\comp{g}{f} \in \sepmorphism{\Bb{A}}$ and $f$ is a proper morphism
    stably continuous and stably in \textsf{E} then
    $g \in \sepmorphism{\Bb{A}}$.
    
  \item \label{item:sep-mor-lt-can} If
    $\comp{g}{f} \in \sepmorphism{\Bb{A}}$ then
    $f \in \sepmorphism{\Bb{A}}$.

  \end{enumerate}

\end{Thm}

\begin{proof}
  Since the kernel pair of a monomorphism $f$ is trivial, $d_f$ is an
  isomorphism. Hence every monomorphism is separated, proving
  \ref{item:sep-mor-contain-mono}{.}  For the rest of the proof, let
  \arrp{f_1}{f_2}{\Kerp{f}}{X}{,} \arrp{g_1}{g_2}{\Kerp{g}}{Y},
  \arrp{h_1}{h_2}{\Kerp{h}}{X}{} the kernel pairs, where
  $h = \comp{g}{f}$. Evidently,
  $\xymatrixcolsep{2.4cm} \xymatrix{ {\Kerp{f}
      \ar[r]^-{\opair{f_1}{f_2}}} & {\Kerp{h}}
    \ar[r]^-{\opair{\comp{f}{h_1}}{\comp{f}{h_2}}} & {\Kerp{g}}
  }$. Consider the commutative the diagram
  \begin{equation*}\label{eq:pb-dg-along-dh}\tag{$\star$}
    \xymatrixcolsep{9.6em}
    \xymatrixrowsep{2.4em}
    \xymatrix{
      {X\,\,} \ar@{>->}[r]^-{d_f} \ar@2{-}[d] & {\Kerp{f}}
      \ar[r]^-{\comp{f}{f_2}} \ar[d]|-{\opair{f_1}{f_2}}
      & {Y} \ar@{>->}[]!<0pt,-12pt>;[d]^-{d_g} \\
      {X\,\,} \ar@{>->}[r]_-{ d_{h}} & {\Kerp{h}}
      \ar[r]|-{\opair{\comp{f}{h_1}}{\comp{f}{h_2}}}
      \ar@{>->}[]!<0pt,-12pt>;[d]_-{\opair{h_{1}}{h_{2}}}
      & {\Kerp{g}} \ar@{>->}[]!<0pt,-12pt>;[d]^-{\opair{g_{1}}{g_{2}}} \\
      & {X \times X} \ar[r]_-{f\times f} & {Y \times Y} 
    }.
  \end{equation*}
  The top right hand square is a pullback square --- if \Arr{p}{P}{Y}
  and \arrp{q}{r}{P}{X} be morphisms such that
  $\comp{h}{q} = \comp{h}{r}$ and
  $\opair{p}{p} = \comp{d_g}{p} =
  \comp{\opair{\comp{f}{h_1}}{\comp{f}{h_2}}}{(q, r)}$ then
  $\comp{f}{q} = p = \comp{f}{r}$. Hence,
  \Arr{\opair{q}{r}}{P}{\Kerp{f}{}} is the unique morphism such that
  $\comp{\opair{f_1}{f_2}}{\opair{q}{r}} = \opair{q}{r}$ and
  $\comp{f}{\comp{f_2}{\opair{q}{r}}} = \comp{f}{r} = p$, proving the
  assertion.  On the other hand the top outer square is trivially a
  pullback square. Hence using properties of pullback squares the
  pullback of \opair{f_1}{f_2}{} along $d_{h}$ is \id{X}.  Finally the
  bottom right hand square is trivially a pullback square.  If $f$ is
  a proper morphism then so also is $f \times f$ (Theorem
  \ref{thm:proper-mor-prop}\ref{item:prod-prop-is-prop}); from the
  right hand bottom pullback square of \eqref{eq:pb-dg-along-dh},
  \opair{\comp{f}{h_1}}{\comp{f}{h_2}}{} is a proper morphism. From
  the top outer square in \eqref{eq:pb-dg-along-dh}{}, since
  $\comp{d_g}{f} = \comp{\opair{\comp{f}{h_1}}{\comp{f}{h_2}}}{d_h}$,
  if $h$ is a separated morphism then \comp{d_g}{f}{} is a proper
  morphism (Theorem
  \ref{thm:proper-mor-prop}\ref{item:proper-mor-pb-stable-comp-closed}{}).
  Hence, if $f$ is a proper morphism stably continuous and stably in
  $\mathsf{E}$, from Theorem \ref{thm:proper-mor-prop}
  \ref{item:proper-mor-rt-can-by-prop-mor-stably-in-E}, $d_{g}$ is
  proper morphism, i.e., $g$ is separated, proving
  \ref{item:sep-mor-rt-can}.  If $h$ is separated then $d_{h}$ is a
  proper morphism. Hence from the left hand pullback square, $d_f$
  being the pullback of $d_{h}$ along \opair{f_1}{f_2}{} is a proper
  morphism. Hence $f$ is a separated, proving
  \ref{item:sep-mor-lt-can}.  If $g$ is separated, $d_{g}$ is a proper
  morphism and hence \opair{f_1}{f_2}{} is proper. Further if $f$ is
  separated, $d_{h} = \comp{\opair{f_{1}}{f_{2}}}{d_f}$ is a proper
  morphism.  Hence h is a separated morphism proving
  \ref{item:sep-mor-comp-closed}{.}  Towards
  \ref{item:sep-mor-pb-stable}{,} consider the diagram in
  \eqref{eq:pb-stab-sep}{} with $f$ a proper morphism.  In
  \eqref{eq:pb-stab-sep}{} the front vertical and the right hand
  vertical squares depict the pullback of $f$ along
  \Arr{p}{\opair{Z}{\psi}}{\opair{Y}{\phi}}{,} the base horizontal square is
  the kernel pair of $f$ and the top horizontal square is the kernel
  pair of $f_p$ the pullback of $f$ along $p$.  Since
  $\comp{f}{\comp{p_f}{(f_p)_2}} = \comp{p}{\comp{f_p}{(f_p)_2}} =
  \comp{p}{\comp{f_p}{(f_p)_1}} = \comp{f}{\comp{p_f}{(f_p)_1}}$,
  there exists the unique morphism \Arr{v}{\Kerp{f_p}}{\Kerp{f}} such
  that $\comp{f_1}{v} = \comp{p_f}{(f_p)_1}$ and
  $\comp{f_2}{v} = \comp{p_f}{(f_p)_2}$.  Furthermore, using
  properties of pullbacks squares, all the faces of the cube are
  pullback squares.
      
  \begin{equation*}
    \label{eq:pb-stab-sep}\tag{$\dagger$}
    \xymatrixcolsep{6em}
    \xymatrixrowsep{2.4em}
    \xymatrix{
      {P} \ar@{>->}[]!<12pt,0pt>;[rr]^-{d_{f_p}}
      \ar@{-->}[dd]_-{!\, w } \ar@2{-}[dr]
      & & {\Kerp{f_p}} \ar[rr]^-{(f_p)_2} \ar[dl]_-{(f_p)_1}
      \ar@{.>}[dd]|-{\hole}|(0.66){!\, v} & &
      {P} \ar[dl]|-{f_p}
      \ar[dd]^-{p_f} \\
      & {P} \ar[rr]|(0.66){f_p} \ar[dd]_-(0.66){p_f} & & {Z}
      \ar[dd]^-(0.66){p} & \\
      {X} \ar@{>->}[]!<12pt,0pt>;[rr]|-{\hole}|(0.66){d_f}
      \ar@2{-}[dr] & &
      {\Kerp{f}}
      \ar[rr]|-{\hole}|(0.66){f_2}
      \ar[dl]^-{f_1} & & {X} \ar[dl]^-f \\
      & {X} \ar[rr]_-f & & {Y} &
    }
  \end{equation*}

  Since $d_f$ is the equaliser of \opair{f_1}{f_2} and
  $\comp{f_1}{\comp{v}{d_{f_p}}} = \comp{f_2}{\comp{v}{d_{f_p}}}$,
  there exists a unique morphism \Arr{w}{P}{X} such that
  $\comp{d_f}{w} = \comp{v}{d_{f_p}}$. Hence
  $w = \comp{f_1}{\comp{d_f}{w}} = \comp{f_1}{\comp{v}{d_{f_p}}} =
  \comp{p_f}{\comp{(f_p)_1}{d_{f_p}}} = p_f$.  Furthermore, since the
  left most square is trivially a pullback square it follows from
  properties of pullbacks that $p_f$ is the pullback of $v$ along
  $d_f$.  Consequently, if $f$ is separated, then $d_{f_p}$ being the
  pullback of $d_f$ along $v$ is also a proper morphism. This proves
  $f_p$ a separated morphism whenever $f$ is a separated
  morphism. Hence \sepmorphism{\Bb{A}} is stable under pullbacks.
  This completes the proof.
\end{proof}


\subsection{}
\label{ssec:hausdorff}

Given \Arr{f}{\opair{X}{\mu}}{\opair{Y}{\phi}}, $f$ is separated if and
only if in the context \slice{\CAL{A}}{Y}, the internal
preneighbourhood space \opair{\opair{X}{f}}{\slice{\mu}{Y}}{} has the
property: the unique morphism $f = \terma{\opair{X}{f}}$ is separated,
i.e., the diagonal morphism from \opair{X}{f}{} to
$\opair{X}{f} \times \opair{X}{f}$ is a proper morphism.

\begin{Df}
  \label{df:hausdorff}
  An internal preneighbourhood space \opair{X}{\mu}{} is
  \emph{Hausdorff} if the unique morphism
  \Arr{\terma{X}}{\opair{X}{\mu}}{\opair{\termo}{\nabla_{\termo}}}{} is
  separated.
\end{Df}

Evidently, if \opair{X}{\mu}{} is Hausdorff and $\mu' \geq \mu$, then
\opair{X}{\mu'}{} is also Hausdorff (Theorem
\ref{thm:sep-map-prop}\ref{item:sep-mor-contain-mono}{}).

\begin{Thm}
  \label{thm:ihs-alt}
  The following are equivalent for any internal preneighbourhood space
  \opair{X}{\mu} of \Bb{A}:
  \begin{enumerate}[label=(\alph*),series=ihs]
  \item \label{item:hausdorffspace} \opair{X}{\mu} is an internal
    Hausdorff space.

  \item \label{item:closeddiagonal} The diagonal morphism
    $\xymatrix{ {\opair{X}{\mu}} \ar@{>->}[]!<20pt,0pt>;[r]^-{d_X} &
      {\opair{X\times{X}}{\mu\times\mu}} }$ is a proper morphism.

  \item \label{item:everymapwithhausdorffdomisseparated} Every
    preneighbourhood morphism with \opair{X}{\mu} as domain is
    separated.

  \item \label{item:existenceofsepmorphismwithhausdorffcodomain} There
    exists a separated preneighbourhood morphism from \opair{X}{\mu} to
    an internal Hausdorff space.

  \item \label{item:closedunderformalepiimage} For every proper
    morphism \Arr{f}{(X, \mu)}{(Y, \phi)} with $f$ stably continuous
    and stably in $\mathsf{E}$, \opair{Y}{\phi} is an internal
    Hausdorff space.

  \item \label{item:productprojectiontoanyissep} The product
    projection \Arr{p_Y}{(X \times Y, \mu \times \phi)}{(Y, \phi)} is
    a separated morphism for every internal preneighbourhood space
    \opair{Y}{\phi}.

  \item \label{item:binproductsofhausdoff} For every internal
    Hausdorff space \opair{Y}{\phi} the product \opair{X \times Y}{\mu
      \times \phi} is an internal Hausdorff space.

  \item \label{item:equalisersofhausdorffvaluedareclosed} If
    $\xymatrixcolsep{1.2cm}\xymatrix{ {(E, \rest{\psi}{E})}
      \ar[r]^-{e} & {(Z, \psi)} \ar@<1ex>[r]^-f \ar@<-1ex>[r]_-g &
      {(X, \mu)} }$ be the equaliser diagram for $f$ and $g$ then $e$
    is a proper morphism.
    
  \end{enumerate}

\end{Thm}

\begin{proof}
  
  Evidently \ref{item:hausdorffspace} and \ref{item:closeddiagonal}
  are equivalent by definition.  Given any preneighbourhood morphism
  \Arr{f}{(X, \mu)}{(Y, \phi)} since
  $\terma{X} = \comp{\terma{Y}}{f}$, an use of Theorem
  \ref{thm:sep-map-prop}\ref{item:sep-mor-lt-can}, shows
  \ref{item:hausdorffspace} implies
  \ref{item:everymapwithhausdorffdomisseparated}. On the other hand,
  \ref{item:everymapwithhausdorffdomisseparated} evidently implies
  \ref{item:hausdorffspace}.  Since \opair{\termo}{\nabla_{\termo}} is
  already an internal Hausdorff space, \ref{item:hausdorffspace}
  automatically implies
  \ref{item:existenceofsepmorphismwithhausdorffcodomain}. On the other
  hand if \opair{Y}{\phi} be an internal Hausdorff space and \Arr{f}{(X,
    \mu)}{(Y, \phi)} is a separated preneighbourhood morphism then
  $\terma{X} = \comp{\terma{Y}}{f}$ implies from Theorem
  \ref{thm:sep-map-prop}\ref{item:sep-mor-comp-closed}{}, $X$ is an
  internal Hausdorff space. Hence
  \ref{item:existenceofsepmorphismwithhausdorffcodomain} implies
  \ref{item:hausdorffspace}.  Given any proper morphism \Arr{f}{(X,
    \mu)}{(Y, \phi)} with $f$ stably continuous and stably in
  $\mathsf{E}$ an use of Theorem
  \ref{thm:sep-map-prop}\ref{item:sep-mor-rt-can-by-prop-mor-stably-in-E}
  prove from $\terma{X} = \comp{\terma{Y}}{f}$ the implication of
  \ref{item:closedunderformalepiimage} from
  \ref{item:hausdorffspace}. On the contrary, assuming
  \ref{item:closedunderformalepiimage} and considering $Y = X$,
  $f = \id{X}$, \ref{item:hausdorffspace} follows.  Since the product
  projection \Arr{p_Y}{(X \times Y, \mu \times \phi)}{(Y, \phi)} is the pullback of
  \terma{X} along \terma{Y}, \ref{item:hausdorffspace} implies
  \ref{item:productprojectiontoanyissep} from pullback stability of
  separated morphisms (Theorem
  \ref{thm:sep-map-prop}\ref{item:sep-mor-pb-stable}{}).  If
  \opair{Y}{\phi} is an internal Hausdorff space and \Arr{p_Y}{(X \times Y, \mu
    \times \phi)}{(Y, \phi)} is the product projection, then
  $\terma{X \times Y} = \comp{\terma{Y}}{p_Y}$ implies
  \ref{item:binproductsofhausdoff} from
  \ref{item:productprojectiontoanyissep} \& Theorem
  \ref{thm:sep-map-prop}\ref{item:sep-mor-comp-closed}{}.  Since any
  internal preneighbourhood space isomorphic to an internal Hausdorff
  space is also an internal Hausdorff space,
  \ref{item:binproductsofhausdoff} evidently implies
  \ref{item:hausdorffspace}.  Since
  $\xymatrix{ {E} \ar[r]^-{e} & {Z} \ar@<1ex>[r]^-f \ar@<-1ex>[r]_-g &
    {X} }$ is an equaliser diagram if and only if \comp{f}{e}{} is the
  pullback of \Arr{\opair{f}{g}}{Z}{X\times{X}}{} along $d_{X}$,
  \ref{item:closeddiagonal} implies
  \ref{item:equalisersofhausdorffvaluedareclosed} from the pullback
  stability of separated morphisms proved in Theorem
  \ref{thm:sep-map-prop}\ref{item:sep-mor-pb-stable}{}. On the other
  hand, since $d_{X}$ is the equaliser of the product projections,
  \ref{item:equalisersofhausdorffvaluedareclosed}{} implies
  \ref{item:closeddiagonal}{.} This completes the proof of Theorem.

\end{proof}

\begin{Cor}
  \label{cor:hausisfincomplete}
  The category \Int{\Haus}{\pNHD{\Bb{A}}} is a finitely complete
  subcategory of \pNHD{\Bb{A}} closed under subobjects and images of
  preneighbourhood morphisms stably continuous and stably in
  $\mathsf{E}$.
\end{Cor}

\begin{proof}
  Since \opair{\termo}{\nabla_{\termo}} is an internal Hausdorff space,
  from \ref{item:binproductsofhausdoff} \&
  \ref{item:equalisersofhausdorffvaluedareclosed} of Theorem the
  category \Int{\Haus}{\Bb{A}} is closed under finite products and
  regular subobjects. Hence \Int{\Haus}{\Bb{A}} is finitely complete.
  Let \opair{X}{\mu} be an internal Hausdorff space and \Arr{f}{Y}{X} be
  a monomorphism. Let $\mu_f$ be the smallest preneighbourhood system on
  $Y$ making $f$ a preneighbourhood morphism. Then \Arr{f}{(Y,
    \mu_f)}{(X, \mu)} is separated morphism (Theorem
  \ref{thm:sep-map-prop}\ref{item:sep-mor-contain-mono}{}) and hence
  \opair{Y}{\mu_f} is an internal Hausdorff space from
  \ref{item:existenceofsepmorphismwithhausdorffcodomain} of Theorem.
  Finally from \ref{item:closedunderformalepiimage} of Theorem if
  \opair{X}{\mu} is an internal Hausdorff space and \opair{Y}{\phi} be an
  internal preneighbourhood space such that $Y = \img{f}{X}$ for some
  preneighbourhood morphism $f$ stably continuous and stably in
  $\mathsf{E}$ then \opair{Y}{\phi} is also an internal Hausdorff space.
\end{proof}


\section{Perfect morphisms}
\label{sec:perfect-mor}
This section discuss \emph{perfect morphisms}, i.e., morphisms which are both proper and separated.

\subsection{}
\label{ssec:perfect-morphisms}

Firstly, some results for preneighbourhood morphisms between compact
and Hausdorff preneighbourhood spaces.

\begin{Thm}
  \label{thm:dom-cpt-cod-haus}
  \begin{enumerate}[label=(\alph*{})]
  \item \label{item:dom-cpt-cod-haus}{} Every preneighbourhood
    morphism from a compact preneighbourhood space to a Hausdorff
    preneighbourhood space is proper.

  \item \label{item:cod-cpt-haus-proper-alt} A preneighbourhood
    morphism with a compact Hausdorff codomain is proper if and only
    if the domain is compact.

  \item \label{item:cpt-sub-haus-is-closed}{} Every compact admissible
    subobject of a Hausdorff preneighbourhood space is closed.
  \end{enumerate}
\end{Thm}

\begin{proof}
  The statement in \ref{item:cod-cpt-haus-proper-alt}{} follows from
  \ref{item:dom-cpt-cod-haus}{} and composition closed property of
  proper morphisms (Theorem
  \ref{thm:proper-mor-prop}\ref{item:proper-mor-pb-stable-comp-closed}). Since
  an admissible subobject is proper if and only if it is closed, the
  statement in \ref{item:cpt-sub-haus-is-closed}{} follows from
  \ref{item:dom-cpt-cod-haus}{.} Towards the proof of
  \ref{item:dom-cpt-cod-haus}{,} if
  \Arr{f}{\opair{X}{\mu}}{\opair{Y}{\phi}}{} is a preneighbourhood morphism
  from a compact preneighbourhood space \opair{X}{\mu}{} to a Hausdorff
  preneighbourhood space \opair{Y}{\phi}{,} and \Arr{p_Y}{X \times Y}{Y} is
  the product projection, then from the pullback square
  $\xymatrixcolsep{4.8em}\xymatrix{ {X}
    \ar@{>->}[]!<3.6ex,0ex>;[r]^-{\opair{\id{X}}{f}} \ar[d]_-f & {X \times
      Y} \ar[d]^-{f\times\id{Y}} \\ {Y} \ar@{>->}[]!<3.6ex,0ex>;[r]_-{d_Y}
    & {Y \times Y} }$, since \opair{Y}{\phi}{} is Hausdorff the diagonal
  $d_{Y}$ is proper (Theorem
  \ref{thm:ihs-alt}\ref{item:closeddiagonal}{}) implying
  \opair{\id{X}}{f}{} is proper. Since \opair{X}{\mu}{} is compact,
  $p_{Y}$ is proper (Remark \ref{rem:compact}). Hence the composite
  $f = \comp{p_Y}{\opair{\id{X}}{f}}$ is proper (Theorem
  \ref{thm:closed-morphisms}\ref{item:closed-comp-closed}), proving
  \ref{item:dom-cpt-cod-haus}{.}
\end{proof}

\begin{Df}
  \label{df:perfect-mor}
  A preneighbourhood morphism \Arr{f}{\opair{X}{\mu}}{\opair{Y}{\phi}}{} is
  a \emph{perfect morphism} if it is both proper and separated. The
  symbol
  $\perfectmorphism{\Bb{A}} = \propermorphism{\Bb{A}} \cap
  \sepmorphism{\Bb{A}}$ is the (possibly large) set of all perfect
  morphisms of \Bb{A}{.}
\end{Df}


\subsection{}
\label{ssec:perfectmaps-prop}

As an immediate consequence of Theorem \ref{thm:proper-mor-prop} \&
\ref{thm:sep-map-prop}:

\begin{Thm}
  \label{thm:perfect-map-prop}
  The set \perfectmorphism{\Bb{A}} of all perfect morphisms of \Bb{A}
  is a pullback stable set, is closed under composition and satisfies
  the properties:
  \begin{enumerate}[label=(\alph*{})]
  \item \label{item"closed-emb-perfect-under-cont}{} If every
    preneighbourhood morphism is continuous then
    $\closedembed{\Bb{A}} \subseteqq \perfectmorphism{\Bb{A}}$.
    
  \item \label{item:sep-mor-rt-can-by-prop-mor-stably-in-E} If
    $\comp{g}{f} \in \perfectmorphism{\Bb{A}}$ and $f$ is a proper
    morphism, stably continuous and stably in \textsf{E} then
    $g \in \perfectmorphism{\Bb{A}}$.
    
  \item \label{item:sep-maps-lt-can-for-prop}{} If
    $\comp{g}{f} \in \propermorphism{\Bb{A}}$ and
    $g \in \sepmorphism{\Bb{A}}$ then $f \in \propermorphism{\Bb{A}}$.

  \item \label{item:sep-mor-lt-can-for-sep}{} If
    $\comp{g}{f} \in \perfectmorphism{\Bb{A}}$ and
    $g \in \sepmorphism{\Bb{A}}$ the $f \in \perfectmorphism{\Bb{A}}$.
  \end{enumerate}

\end{Thm}

\begin{proof}
  It is enough to prove \ref{item:sep-maps-lt-can-for-prop}. Let
  \Arr{f}{\opair{X}{\mu}} {\Arr{g}{\opair{Y}{\phi}}{\opair{Z}{\psi}}}{} be
  preneighbourhood morphisms such that \comp{g}{f}{} is proper and $g$
  is separated. Evidently, in the context \slice{\CAL{A}}{Z},
  \Arr{f}{\opair{\opair{X}{\comp{g}{f}}}
    {\slice{\mu}{Z}}}{\opair{\opair{Y}{g}}{\slice{\phi}{Z}}}{} is a
  preneighbourhood morphism from a compact preneighbourhood space to a
  Hausdorff preneighbourhood space. Hence from Theorem
  \ref{thm:dom-cpt-cod-haus}{}\ref{item:dom-cpt-cod-haus}{,} $f$ is
  proper in \slice{\CAL{A}}{Z}{,} and hence proper in \CAL{A}{.} The
  statement in \ref{item:sep-mor-lt-can-for-sep}{} follows from
  \ref{item:sep-maps-lt-can-for-prop}{} and Theorem
  \ref{thm:sep-map-prop}\ref{item:sep-mor-lt-can}.

\end{proof}


\subsection{Three types of internal preneighbourhood spaces}
\label{ssec:some-spaces}

In conclusion, three important types of internal preneighbourhood spaces are defined. Detailed investigation of these spaces shall be done in later papers.

\begin{Df}
  \label{df:some-spaces}
  An internal preneighbourhood space \opair{X}{\mu}{} is:
  \begin{enumerate}[label=(\alph*)]
  \item \label{item:cpt-hausdorff}
    \emph{compact Hausdorff} if \Arr{\terma{X}}{\opair{X}{\mu}}{\opair{\termo}{\nabla_{\termo}}}{} is a perfect morphism.

  \item \label{item:tychonoff}{}
    \emph{Tychonoff} if there exists a morphism $\xymatrix{ {\opair{X}{\mu}} \ar@{>->}[]!<21.6pt,0pt>;[r]^-m & {\opair{Y}{\phi}} }$, where \opair{Y}{\phi}{} is a compact Hausdorff internal preneighbourhood space and $m \in \mathsf{M}$.

  \item \label{item:abs-closed}{}
    \emph{absolutely closed} if for every morphism $\xymatrix{ {\opair{X}{\mu}} \ar@{>->}[]!<21.6pt,0pt>;[r]^-m & {\opair{Y}{\phi}} }$ where \opair{Y}{\phi}{} is a Hausdorff internal preneighbourhood space and $m \in \mathsf{M}$ the morphism $m \in \mc{\phi}$.
  \end{enumerate}

  The symbols \Int{\KHaus}{\pNHD{\Bb{A}}}{,} \Int{\Tych}{\pNHD{\Bb{A}}}{,} \Int{\AbCl}{\pNHD{\Bb{A}}}{} respectively denote the full subcategories of compact Hausdorff, Tychonoff, absolutely closed internal preneighbourhood spaces.
\end{Df}



\section{Concluding Remarks}
\label{sec:concluding-remarks}

Let $\CAL{A} = (\Bb{A}, \mathsf{E}, \mathsf{M})$ be a context.

Given (possibly large) sets \eulf{a}, \eulf{b}{} of morphisms of \Bb{A}, the phrases \emph{\eulf{b}{} is composition closed} or \emph{\eulf{b}{} is (pullback) stable} is well known; the set \eulf{b} shall be said to be \emph{left \eulf{a}{} cancellative} (respectively, \emph{right \eulf{a}{} cancellative}) if $\comp{g}{f} \in \eulf{b}$ and $g \in \eulf{a}$ (respectively, $f \in \eulf{a}$) implies $f \in \eulf{b}$ (respectively, $g \in \eulf{b}$). The set \eulf{b}{} is \emph{stably in $\mathsf{E}$} if every pullback $f_{g}$ of $f$ along any morphism $g$ is in $\mathsf{E}$. If \eulf{b}{} is a set of preneighbourhood morphisms then it is said to be \emph{stably continuous} if for any $\mu$-$\phi$ continuous preneighbourhood morphism \Arr{f}{\opair{X}{\mu}}{\opair{Y}{\phi}}{} in \eulf{b}, and for any preneighbourhood morphism \Arr{g}{\opair{Z}{\psi}}{\opair{Y}{\phi}}, the pullback \Arr{f_g}{\opair{X\times_YZ}{\mu\times_{\phi}\psi}}{\opair{Z}{\psi}}{} of $f$ along $g$ is $(\mu\times_{\phi}\psi)$-$\psi$ continuous and is also in \eulf{b}. Table \ref{tab:mor-types-comp-table}{} summarise the properties
deduced in this paper. The cells highlighted in {\color{col} this colour} are the properties where the \emph{continuity} condition is required; the others do not require \emph{continuity} of the involved preneighbourhood morphism, and hence are purely consequences of the preneighbourhood morphism property.

\begin{table}
  \centering
  \begin{threeparttable}
  \begin{tabular}{|m{0.06\textwidth}|m{0.084\textwidth}|m{0.264\textwidth}|
    m{0.24\textwidth}|m{0.264\textwidth}|} \hline
    & contains & stability & closed under composition & cancellation properties
    \\[2ex]\hline

    \closedmorphism{\Bb{A}} &
    \Iso{\Bb{A}}
    & \cellcolor{col}{
      $m \in \closedembed{\Bb{A}}, f\in\Bb{A}_{\mathtt{clc}}$ {}
      $\Rightarrow {f_m\in\closedmorphism{\Bb{A}}}$ \newline{}
      $m \in \closedembed{\Bb{A}}, f\in\Bb{A}_{\mathtt{c}}$ {}
      $\Rightarrow \finv{f}{m}\in\closedembed{\Bb{A}}$}
    & composition closed
    & \cellcolor{col}{right $\Bb{A}_{\mathtt{fsc}}$ cancellative} \\[2ex] \hline

    \densemorphism{\Bb{A}} & $\mathsf{E}$
         & & \cellcolor{col}{
             $g \in\Bb{A}_{\mathtt{dc}},f\in\densemorphism{\Bb{A}}$ {}
             $\Rightarrow\comp{g}{f}\in\densemorphism{\Bb{A}}$}
                                                     & right $\Bb{A}_1$ cancellative \\[2ex] \hline

    {\propermorphism{\Bb{A}}} &
    {\cellcolor{col}{\clemb{\Bb{A}}, {} in
                              $\mathtt{RZC}$}} & 
    {pullback stable} &
    {composition closed}
               & \cellcolor{col}{right
                 $\Bb{A}_{\mathtt{st(E,c)}}$ cancellative} \\\cline{5-5}
    & & & & left \Mono{\Bb{A}} cancellative \\[2ex]\hline

    \sepmorphism{\Bb{A}}{} & \Mono{\Bb{A}}{} & pullback stable & composition closed & \cellcolor{col}{right $\Bb{A}_{\mathtt{st(E,c,cl)}}$ cancellative} \\\cline{5-5}
    & & & & left $\Bb{A}_1$ cancellative \\[2ex]\hline

    \perfectmorphism{\Bb{A}} & \cellcolor{col}{\clemb{\Bb{A}}, {} in $\mathtt{RZC}$}
               & pullback stable & composition closed &
                                                        \cellcolor{col}{right $\Bb{A}_{\mathtt{st(E,c,cl)}}$ cancellative} \\\cline{5-5}
    & & & & left \perfectmorphism{\Bb{A}} {} cancellative \\[2ex]\hline 
                                                                 
  \end{tabular}
  
  \begin{tablenotes}
  \item[1] $\Bb{A}_1$ is the (possibly large) set of all morphisms
  \item[2] \closedmorphism{\Bb{A}} is the (possibly large) set of all closed morphisms
  \item[3] \clemb{\Bb{A}}{} is the (possibly large) set of all closed embeddings
  \item[4] \densemorphism{\Bb{A}} is the (possibly large) set of all dense preneighbourhood morphisms
  \item[5] \propermorphism{\Bb{A}} is the (possibly large) set of all proper preneighbourhood morphisms
  \item[6] \sepmorphism{\Bb{A}} is the (possibly large) set of all separated preneighbourhood morphisms
  \item[7] \perfectmorphism{\Bb{A}} is the (possibly large) set of all perfect preneighbourhood morphisms
  \item [8] $\Bb{A}_{\mathtt{c}}$ is the (possibly large) set of all continuous preneighbourhood morphisms
  \item[9] $\Bb{A}_{\mathtt{dc}}$ is the (possibly large) set of all dense and continuous preneighbourhood morphisms
  \item[10] $\Bb{A}_{\mathtt{fsc}}$ is the (possibly large) set of all formally surjective and continuous preneighbourhood morphisms
  \item[11] $\Bb{A}_{\mathtt{clc}}$ is the (possibly large) set of all closed and continuous preneighbourhood morphisms
  \item[11] $\Bb{A}_{\mathtt{st(E,c)}}$ is the (possibly large) set of all preneighbourhood morphisms which are stably continuous and stably in $\mathsf{E}$
  \item[12] $\Bb{A}_{\mathtt{st(E,c,cl)}}$ is the (possibly large) set of all preneighbourhood morphisms which are stably continuous, stably in $\mathsf{E}$ and stably closed
  \item[13] $\mathtt{RZC}$ abbreviates \emph{reflecting zero context}
  \item[14] the cells in \emph{{\color{col}{this colour}}} indicate the presence of \emph{continuity} in the assertion
  \item[15] additionally, every $\mathtt{RZC}$ with continuous preneighbourhood morphisms has \opair{\densemorphism{\Bb{A}}}{\clemb{\Bb{A}}}{} factorisation structure
  \end{tablenotes}

  \caption{Comparative list of properties}
  \label{tab:mor-types-comp-table}

\end{threeparttable}
\end{table}

The following definition appears in \cite[][\S{2}]{HofmannTholen2012}: 
\begin{Df}
  A pullback stable (possibly large) set \eulf{a}{} of morphisms of \Bb{A}{} is called a \emph{topology} if it contains isomorphisms and is closed under compositions.

  If \eulf{a} is a topology and right \eulf{a}{} cancellative, a topology \eulf{b}{} is called a \emph{\eulf{a}-topology} if it is right \eulf{a} cancellative.
\end{Df}
Drawing inspiration from \cite[][]{ClementinoGiuliTholen2004}, it is observed in \cite[see][\S{2}]{HofmannTholen2012} that in case when a finitely complete category \Bb{A}{} with a proper \fact{\mathsf{E}}{\mathsf{M}}{} system has a set \closedmorphism{\Bb{A}}{} of closed morphisms described by axioms \cite[see][Axioms (F3)-(F5)]{ClementinoGiuliTholen2004}, then the set of proper morphisms (i.e., morphisms stably in \closedmorphism{\Bb{A}}{}) is a \eulf{s}{-}topology, where \eulf{s}{} is the set of morphisms stably in $\mathsf{E}$. 

In terms of Definition, Table \ref{tab:mor-types-comp-table}{} shows the set  $\Bb{A}_{\mathtt{st(E,c,cl)}}$ is a right $\Bb{A}_{\mathtt{st(E,c,cl)}}$ cancellative topology and each of the sets \propermorphism{\Bb{A}}{,} \sepmorphism{\Bb{A}}{,} \perfectmorphism{\Bb{A}} are $\Bb{A}_{\mathtt{st(E,c,cl)}}$-topologies. The difference between the two approaches arises from the fact that in the present case \closedmorphism{\Bb{A}}{} is right $\Bb{A}_{\mathtt{fsc}}$ ($\subset \mathsf{E}$) cancellative, while the axioms of \cite{ClementinoGiuliTholen2004} assert \closedmorphism{\Bb{A}}{} is right $\mathsf{E}$ cancellative. In case when \CAL{A}{} is $\mathtt{RZC}$ and $\mathsf{E}$ is pullback stable the present case reduces to the situation considered in \cite{ClementinoGiuliTholen2004}.


\printbibliography


\end{document}
